\documentclass[12pt]{article}
\usepackage{epsfig}
\usepackage{latexsym}
\usepackage{amssymb}
\title{\bf Transcendence and CM on Borcea-Voisin towers of Calabi-Yau manifolds}
\author{\bf Paula Tretkoff\footnote{The previous name of the author was Paula B. Cohen.}}
\begin{document}
\maketitle
\begin{abstract} \noindent This paper is a sequel to \cite{TCY}, in which we showed the validity of a special case of a conjecture of Green, Griffiths and Kerr \cite{GGK} for certain families of Calabi-Yau manifolds over Hermitian symmetric domains. These results are analogues of a celebrated theorem of Th. Schneider \cite{Sch} on the transcendence of values of the elliptic modular function, and its generalization to the context of abelian varieties in \cite{Co1}, \cite{SW}. In the present paper, we apply related techniques to many of the examples of families of Calabi-Yau varieties with dense sets of CM fibers in the work of Rohde \cite{Roh}, and in particular to Borcea-Voisin towers. Our results fit into a broader context of transcendence theory for variations of Hodge structure of higher weight.\footnote{Research supported by NSF grant DMS--0800311 and NSA grant 1100362.

\medskip

\noindent Department of Mathematics, Texas A\&M University, College Station, TX 77842--3368, USA, and,\\
CNRS, UMR 8524, Universit\'e Lille 1, Cit\'e Scientifique 59655 Villeneuve d'Ascq C\'edex,  FRANCE.}
\end{abstract}

Calabi-Yau manifolds

\medskip

14C30 (Secondary)

\newpage

\section{Introduction}\label{s:intro}

\noindent Transcendence and linear independence properties of periods and of quasi-periods of algebraic $1$-forms on group varieties defined over number fields are well understood due to results of many researchers, all having their origins in Baker's method. Since the mid-1980's, these results can all be viewed as applications of W\"ustholz's Analytic Subgroup Theorem. In his 1986 lecture at the ICM in Berkeley, W\"ustholz proposed some results on the transcendence of periods of algebraic forms of higher order and recognized the importance of Hodge Theory for such problems. Nonetheless, relatively little is known in this direction. One can also study the transcendence properties of normalized period matrices in Griffiths' period domains.

\noindent The theory of Mumford-Tate domains, developed by Green, Griffiths and Kerr \cite{GGK}, aims at generalizing the theory of Shimura varieties (moduli spaces for PEL level 1 polarized Hodge structures) to certain analogous moduli spaces for Hodge structures of higher level with some additional structure. Although their main interest is in the arithmetic, geometry and representation theory of Mumford-Tate domains, they also propose some problems in transcendental number theory. This motivated both a previous publication of the author with M.D. Tretkoff \cite{TCY} and the present paper, even though the results obtained so far can be viewed in the context of much older work of Griffiths. In \cite{GGK} a conjecture is formulated for geometric variations of Hodge structure whose affirmation generalizes a classical result of Th. Schneider on the transcendence of values of the classical elliptic modular function at algebraic arguments, and its generalization to the Siegel modular case, which is a joint result of the author (formerly Paula B. Cohen) with Shiga and Wolfart \cite{Co1}, \cite{SW}. One way to formulate Schneider's result is to say that an elliptic curve is defined over a number field, and the ratio of its periods is algebraic, if and only if the elliptic curve has complex multiplication. The complex multiplication property for an elliptic curve ${\mathcal E}$ means that the endomorphism algebra of ${\mathcal E}$ is an imaginary quadratic number field, or, equivalently, that the Mumford-Tate group of the usual rational Hodge structure on $H^1({\mathcal E},{\mathbb Q})$ is abelian. There are several natural generalizations of this problem to Hodge structures of higher level, and their geometric variations, which we discuss in \S\S\ref{s:Schneider},\ref{s:VHS} of this paper. Let $X$ be a smooth projective variety of dimension $n$ defined over a number field, and for any $k\le n$, consider the usual rational Hodge structure $H^k(X,{\mathbb Q}_X)_{\rm prim}$, determined by the complex structure on $X$, on the primitive cohomology $H^k(X,{\mathbb Q})_{\rm prim}$. Associated to the Hodge decomposition of $H^k(X,{\mathbb C})_{\rm prim}$ into a direct sum of $(p,q)$ forms $H^{p,q}_{\rm prim}$, $p+q=k$, we have the Hodge filtration $F^{\ast,k}$, where $F^{p,k}=\oplus_{p'\ge p}H^{p',k-p'}_{\rm prim}$. We say that $F^{\ast,k}$ is defined over ${\overline{\mathbb Q}}$ if it is induced by extension of scalars to ${\mathbb C}$ from a filtration of $H^k(X,{\overline{\mathbb Q}})_{\rm prim}$ by ${\overline{\mathbb Q}}$-vector subspaces. One way to generalize the result of Th. Schneider is to ask whether, if $F^{\ast,k}$ is defined over ${\overline{\mathbb Q}}$, it follows that $H^k(X,{\mathbb Q}_X)$ has CM. By CM, we mean that the Mumford-Tate group of $H^k(X,{\mathbb Q}_X)$ is abelian (for details, see \S\ref{s:Schneider}). The level $k$ of interest depends on the variety $X$ that one is studying. For example, when $X$ is Calabi-Yau, then either the case $k=n$ or $k\le n$ is usually of interest, whereas, if $X$ is an abelian variety, the classical definition of CM corresponds to the case $k=1$. A natural restatement of this question for geometric variations of Hodge structure is given in \S\ref{s:VHS}. These problems can be restated in terms of normalized period matrices in Griffiths' domains of marked polarized Hodge structures, and this is more in keeping with the original result of Schneider and its generalization in \cite{Co1}, \cite{SW}, but we do not pursue this in the present paper.

\noindent In \cite{TCY}, with M.D. Tretkoff, we looked at such questions for certain families of Calabi-Yau varieties. The examples are taken from work of Rohde, although his interest is not in transcendence, but rather in finding families of projective Calabi-Yau varieties with Zariski-dense sets of CM fibers. In the present paper, we rework these examples and treat more cases from the book of Rohde.  A new ingredient is the treatment of a general Borcea-Voisin tower of Calabi-Yau varieties with involution in \S \ref{s:BVtower}, which allows for the construction of infinitely many examples of Calabi-Yau varieties of arbitrary dimension that satisfy an appropriate generalization of Schneider's Theorem. In \S\ref{s:expect1}, examples of families of Calabi-Yau varieties with involution are given for which a generalized Schneider's Theorem holds. They can therefore be used for the general inductive step of a Borcea-Voisin tower. These examples are the so-called Viehweg-Zuo towers, originating in \cite{VZ}, and are iterated cyclic covers of families of algebraic curves. They include the example already considered in \cite{TCY}. In \S\ref{s:examples}, we reconsider these examples from a more explicit point of view that involves computation of periods of forms of higher order. In \S\ref{s:exampleslow}, we consider examples of Borcea-Voisin towers in low dimension to show the necessity of some of the assumptions made in \S\ref{s:BVtower}. In \S\ref{s:lemmas}, we collect some basic lemmas from Hodge theory and multi-linear algebra which are used throughout the paper.

\noindent The basic transcendence techniques used in all the results of this paper are those of \cite{Co1}, \cite{Sch}, \cite{SW}, which can all be viewed as applications of W\"ustholz' Analytic Subgroup Theorem. Therefore, we must necessarily relate the family of Calabi-Yau varieties we are studying to some family of Hodge structures of level 1. It would be interesting to find a new method that relates directly to Hodge structures of higher level.

\noindent A problem posed by Schneider asks for a proof of the transcendence of the values of the elliptic modular function at non-CM points which does not use elliptic functions, that is, which only uses the intrinsic properties of the elliptic modular function and not the fact that its values parametrize complex isomorphism classes of elliptic curves. Similarly, in this paper, we do not in any way exploit the intrinsic properties of Griffiths' period domains.

\noindent Finally, I would like to thank my husband, Marvin D. Tretkoff, for many helpful conversations about the mathematics in this paper, as well as for his continuing love and support.

\section{Hodge Filtrations defined over ${\overline{\mathbb Q}}$ and CM}\label{s:Schneider}

\noindent We begin by recalling some definitions and results from the Hodge theory of smooth projective manifolds. We do not include proofs as there are many excellent references on this topic, for example \cite{Cat}, \cite{GrH}. Let $X$ be a smooth complex projective manifold of dimension $n$. We also denote by $X$ its set $X({\mathbb C})$ of complex points. Then $X$ is a compact K\"ahler manifold and has a natural complex structure inducing an action of ${\mathbb C}^\ast={\mathbb C}\setminus\{0\}$ on the bundle of $k$-forms $\Omega^k$. For non-negative integers $p,q$ with $p+q=k$, a $(p,q)$-form is a section of $\Omega^k$ on which $z\in {\mathbb C}^\ast$ acts by multiplication by $z^p\overline{z}^q$. A $(p,q)$-form $\omega$ can be written locally as
$$
\omega=\sum f_{i_1\ldots i_pj_1\ldots j_q}dz_{i_1}\wedge\ldots\wedge dz_{i_p}\wedge d{\overline{z}}_{j_1}\wedge \ldots\wedge d{\overline{z}}_{j_q},
$$
where $z_1,z_2,\ldots,z_n$ are local holomorphic coordinates.
We have a direct sum decomposition
$$
\Omega^k=\oplus_{p+q=k}\Omega^{p,q}.
$$
Using this, we can write the exterior differential $d$ as a direct sum $d=d'+d''$ where $d'$ takes a $(p,q)$-form to a $(p+1,q)$-form and $d''$ takes a $(p,q)$-form to a $(p,q+1)$-form. For $k=0,\ldots,n$, the de-Rham cohomology with complex coefficients is the complex vector space $H^k_\Omega(X,{\mathbb C})$ consisting of the closed $k$-forms modulo the exact $k$-forms. The Bott-Chern cohomology $H^{p,q}$ is defined to be the space of closed $(p,q)$-forms modulo the image under $d'd''$ of the $(p-1,q-1)$-forms. There is a natural injective map from $H^{p,q}(X)$ to $H_\Omega^{p+q}(X)$ and we denote the image also by $H^{p,q}(X)$. When $X$ is a compact K\"ahler manifold, there is an internal direct sum decomposition
$$
H^k_\Omega(X,{\mathbb C})=\oplus_{p+q=k}H^{p,q}(X).
$$
Moreover, the action of ${\mathbb C}^\ast$ on forms induces an action on cohomology. The above direct sum decomposition is called the Hodge decomposition, which is independent of the choice of K\"ahler metric on $X$, and the $h^{p,q}=\dim(H^{p,q}(X))$ are called the Hodge numbers. There is an associated filtration $F^{\ast,k}_\Omega$ of $H^k_\Omega(X,{\mathbb C})$ by ${\mathbb C}$-vector spaces given by
$$
\{0\}\subset F^{k,k}_\Omega\subset F^{k-1,k}_\Omega\subset\ldots\subset F^{1,k}_\Omega\subset F^{0,k}_\Omega=H^k_\Omega(X,{\mathbb C})
$$
where
$$
F^{p,k}_\Omega=\oplus_{p'\ge p}\,H^{p',k-p'}(X).
$$
This is called the Hodge filtration of $H^k_\Omega(X,{\mathbb C})$. Denoting $f^{p,k}=\sum_{p'\ge p}h^{p',k-p'}$, $p=0,\ldots n$, we have $\dim(F^{p,k}_\Omega)=f^{p,k}$. The decreasing sequence $$f^{\ast,k}=(f^{0,k},\ldots,f^{k,k})$$ is called the signature of the filtration. Notice that
$$
F^{p,k}_\Omega\oplus {\overline{F}}^{k-p+1}_\Omega=H^k_\Omega(X,{\mathbb C}).
$$

\noindent Consider now the singular cohomology group $H^k(X,{\mathbb Z})$, $k=0,\ldots,\dim(X)$, and let $V^k=V^k_{\mathbb Q}$ be the ${\mathbb Q}$-vector space $H^k(X,{\mathbb Q})=H^k(X,{\mathbb Z})\otimes{\mathbb Q}$. For a field $K$ containing ${\mathbb Q}$, and any ${\mathbb Q}$-vector space $U$, let $U_K=U\otimes_{\mathbb Q}K$. Then $V_K^k=H^k(X,K)=H^k(X,{\mathbb Q})\otimes_{\mathbb Q}K$, and we denote ${\rm{GL}}(V^k)_K={\rm{GL}}(V_K^k)$. By the de-Rham Theorem with complex coefficients, there is a canonical isomorphism
$$
\iota_k:H^k(X,{\mathbb C})\simeq H^k_\Omega(X,{\mathbb C}).
$$
We denote by $F^{\ast,k}=F^{\ast,k}(X)$ the filtration of $H^k(X,{\mathbb C})$ by the complex vector spaces $F^{p,k}=\iota_k^{-1}(F^{p,k}_\Omega)$, $p=0,\ldots,k$. In particular $f^{p,k}=\dim_{\mathbb C}(F^{p,k})$, $p=0,\ldots,k$, and
$$
F^{p,k}\oplus {\overline{F}}^{k-p+1}=H^k(X,{\mathbb C}).
$$
We denote by $H^k(X,{\mathbb Q}_X)$ the vector space $V^k$ together with this filtration $F^{\ast,k}$, and call it the rational Hodge structure on $H^k(X,{\mathbb Q})$. We call the filtration $F^{\ast,k}=F^{\ast,k}(X)$ the
Hodge filtration of $H^k(X,{\mathbb Q}_X)_{\mathbb C}$ to emphasize its dependence on the complex manifold structure of $X$.

\medskip

\noindent{\bf Definition 1:} Let $W$ be a ${\mathbb Q}$-vector space, and $g^\ast=\{g^p\}_{p=0}^k$ a decreasing sequence of $k+1$ positive integers with $g^0=\dim_{\mathbb Q}(W)$. A $K$-filtration of $W_K=W\otimes_{\mathbb Q}K$ of signature $g^\ast$ is a filtration of $W_K$ by $K$-vector subspaces of the form
$$
\{0\}\subset G^k_K\subset G^{k-1}_K\subset\ldots\subset G^1_K\subset G^0_K=W_K,
$$
with $\dim_KG^p_K=g^p$, $p=0,\ldots,k$.

\medskip

\noindent We now introduce the notion of a Hodge filtration defined over $K$. It is equivalent to the condition that the filtration correspond to a $K$-rational point of the compact dual of an appropriate Griffiths period domain. In the polarized case, this equivalent notion can be found in \cite{GGK}, and we describe it briefly in \S\ref{s:VHS} of this paper. Let $L/K$ be a field extension. A $K$-vector space $U$ can be extended to an $L$-vector space $L\otimes_K U$, and any $K$-basis of $U$ yields an $L$-basis of $L\otimes_K U$. On the other hand, for an $L$-vector space $W$, a $K$-subspace $U$ such that a $K$-basis of $U$ is an $L$-basis of $W$ is called a $K$-form of $W$.

\medskip

\noindent {\bf Definition 2:} Let $X$ be a compact K\"ahler manifold of dimension $n$. For $k=0,\ldots,n$, let $f^{\ast,k}=\{f^{p,k}\}_{p=0}^k$, where $f^{p,k}=\sum_{p'\ge p}h^{p',k-p'}$. The Hodge filtration of $H^k_\Omega(X,{\mathbb C})$ is defined over $K$ if there is a $K$-filtration of $H^k(X,K)$, of signature $f^{\ast,k}$,
$$
\{0\}\subset F^{k,k}_K\subset F^{k-1,k}_K\subset\ldots\subset F^{1,k}_K\subset F^{0,k}_K=H^k(X,K),
$$
such that $\iota_k(F^{p,k}_K\otimes_K{\mathbb C})=F^{p,k}_\Omega$. In other words, the $K$-vector subspace $F^{p,k}_K$ of $H^k(X,K)$ is a $K$-form for the ${\mathbb C}$-vector subspace $F^{p,k}=\iota_k^{-1}(F^{p,k}_\Omega)$ of $H^k(X,{\mathbb C})$. We also say that the Hodge filtration $F^{\ast,k}=F^{\ast,k}(X)$ of $H^k(X,{\mathbb Q}_X)_{\mathbb C}$ is defined over $K$.

\medskip

 \noindent Whereas, for trival reasons, filtrations fulfilling Definition 1 always exist, the filtrations described in Definition 2 do not always exist. Like any ${\mathbb C}$-vector space, the vector spaces $F^{p,k}$ have $K$-forms, but there is no reason \emph{a priori} for these $K$-forms to be subspaces of $H^k(X,K)$, since this requires there to be a ${\mathbb C}$-basis of $F^{p,k}$ with elements in $H^k(X,K)$.

\medskip

\noindent We can associate to the Hodge filtration of $H^k(X,{\mathbb Q}_X)_{\mathbb C}$ a homomorphism of ${\mathbb R}$-algebraic groups
$$\widetilde{\varphi}:{\mathbb S}({\mathbb R})\rightarrow {\rm{GL}}(V^k)_{\mathbb R}$$
with $\widetilde{\varphi}({\rm{diag}}(r,r))=r^k{\rm{Id}}_{V^k}$, for $r\in {\mathbb R}$, $r\not=0$. As before $V^k=H^k(X,{\mathbb Q})$. Here ${\mathbb S}$ is the group with $K$-points
$$
{\mathbb S}(K)=\left\{\begin{pmatrix}{a&-b\cr b&a}\end{pmatrix}: a^2+b^2\not=0,\, a,b\in K\right\},
$$
which is often called the Deligne torus. Any $z\in {\mathbb C}$, $z\not=0$, is viewed as the matrix $\begin{pmatrix}{\Re(z)&-\Im(z)\cr\Im(z)&\Re(z)}\end{pmatrix}$ in ${\mathbb S}({\mathbb R})$. Indeed, this defines an isomorphism of ${\mathbb R}$-algebraic groups between ${\mathbb S}({\mathbb R})$ and the restriction of scalars ${\rm{Res}}_{{\mathbb C}/{\mathbb R}}({\mathbb C}^\ast)$. Given the Hodge filtration $F^{\ast,k}=F^{\ast,k}(X)$, we define $\widetilde{\varphi}$ by
$$
\widetilde{\varphi}\left(\begin{pmatrix}{a&-b\cr b&a}\end{pmatrix}\right)(v_{p,q})=(a+ib)^p(a-ib)^qv_{p,q},\qquad {\rm for\; all}\;\; v_{p,q}\in F^{p,k}\cap{\overline{{F}}}^{q,k}.
$$
Given ${\widetilde{\varphi}}$, we can recover the Hodge decomposition of $H^k(X,{\mathbb C})$ by setting $H^{p,q}$ equal to the vectors $v_{p,q}$ satisfying the above formula, and then in turn recover the Hodge filtration $F^{\ast,k}$. In this way, we see that the Hodge filtration and the associated ${\widetilde{\varphi}}$ are equivalent.
We call $(H^k(X,{\mathbb Q}),{\widetilde{\varphi}})=(H^k(X,{\mathbb Q}),F^{\ast,k})$ the Hodge structure of weight $k$ on $H^k(X,{\mathbb Q})$, and again denote it by $H^k(X,{\mathbb Q}_X)$. As already remarked, it is the natural rational Hodge structure induced by the complex structure on $X$.

\noindent More abstractly, we call any pair $(V^k,\widetilde{\varphi})$, where $V^k$ is a ${\mathbb Q}$-vector space, $k$ is an integer, and $\widetilde{\varphi}:{\mathbb S}({\mathbb R})\rightarrow {\rm GL}(V^k)_{\mathbb R}$ is a homomorphism of ${\mathbb R}$-algebraic groups with the above properties, a (pure rational) Hodge structure of weight $k$. Two Hodge structures $(V_1^k,\widetilde{\varphi}_1)$ and $(V_2^k,\widetilde{\varphi}_2)$ of level $k$ are isomorphic if there is a ${\mathbb Q}$-vector space isomorphism from $V_1^k$ to $V_2^k$ commuting with the actions of ${\mathbb S}({\mathbb R})$ defined by $\widetilde{\varphi}_1$ and $\widetilde{\varphi}_2$. If the two rational Hodge structures arise from extension of scalars to ${\mathbb Q}$ of two Hodge strucures over ${\mathbb Z}$, such an isomorphism is sometimes called an isogeny of the Hodge structures over ${\mathbb Z}$. A rational Hodge structure is equivalent to giving the real structure of $V_{\mathbb C}$ together with a filtration $F^{\ast,k}$ of $V_{\mathbb C}$ satisfying $F^{p,k}\oplus {\overline F}^{k-p+1,k}=V_{\mathbb C}$, $p=0,\ldots,k$. For more details, see \cite{GGK}, and also \cite{TCY} for a discussion of some of the transcendence questions raised in the present paper from the more abstract viewpoint.

\noindent We again let $V^k=H^k(X,{\mathbb Q})$, where $X$ is a smooth projective variety of dimension $n$, and, for $k\le n$, let $F^{\ast,k}=F^{\ast,k}(X)$ be the Hodge filtration of $V^k_{\mathbb C}$ with associated ${\mathbb R}$-algebra homomorphism ${\widetilde{\varphi}}$ from ${\mathbb S}({\mathbb R})$ to ${\rm GL}(V^k)_{\mathbb R}$. We define $M_{\widetilde{\varphi}}$ to be the smallest ${\mathbb Q}$-algebraic subgroup of ${\rm GL}(V^k)$ whose real points contain $\widetilde{\varphi}(S({\mathbb R}))$. The group $M_{\widetilde{\varphi}}$ is known as the Mumford-Tate group of $\widetilde{\varphi}$, and we denote its set of $K$-points by $M_{\widetilde{\varphi}}(K)$, where $K$ is any field extension of ${\mathbb Q}$. As $\widetilde{\varphi}$ is equivalent to the Hodge filtration $F^{\ast,k}=F^{\ast,k}(X)$ of $H^k(X,{\mathbb Q}_X)_{\mathbb C}$, we also call $M_{\widetilde{\varphi}}$ the Mumford-Tate group of $F^{\ast,k}$ and denote it by $M_{F^{\ast,k}}$.

 \medskip

 \noindent {\bf Definition 3:} Let $X$ be a smooth projective variety of dimension $n$ and let $k$ be an integer with $0\le k\le n$. Let $F^{\ast,k}=F^{\ast,k}(X)$ be the Hodge filtration of $H^k(X,{\mathbb Q}_X)_{\mathbb C}$. Let $M_{F^{\ast,k}}$ be the Mumford-Tate group of $F^{\ast,k}$ and $H_{F^{\ast,k}}$ be the stabilizer of $F^{\ast,k}$ in ${\rm GL}(H^k(X,{\mathbb Q}))_{\mathbb R}$. Then $F^{\ast,k}$ has complex multiplication (CM) if and only if $M_{F^{\ast,k}}({\mathbb R})$  is contained in  $H_{F^{\ast,k}}$. We also say that $H^k(X,{\mathbb Q}_X)$ has CM.

 \medskip

 \noindent It is well-known that if the Hodge filtration has CM, then it is defined over ${\overline{\mathbb Q}}$ (for a proof, see \cite{UllY}). It is also well-known that the Hodge filtration has CM if and only if its Mumford-Tate group is an algebraic torus (see for example \cite{GGK}, V.4). With notation as in the paragraph before Definition 3, let ${\rm End}_0(V^k_{\mathbb Q},{\widetilde{\varphi}})$ be the ${\mathbb Q}$-algebra of ${\mathbb Q}$-linear endomorphisms of $V^k_{\mathbb Q}=H^k(X,{\mathbb Q})$ commuting with the elements in the image of $\widetilde{\varphi}$. Then ${\rm End}_0(V^k_{\mathbb Q},{\widetilde{\varphi}})$ consists of the ${\mathbb Q}$-linear endomorphisms of $V^k_{\mathbb Q}$ commuting with the elements of $M_{\widetilde{\varphi}}$, since the Mumford-Tate group is generated by the elements in the image of $\widetilde{\varphi}$ and their conjugates over ${\mathbb Q}$. If $M_{{\widetilde{\varphi}}}$ is contained in the stabilizer of ${\widetilde{\varphi}}$ in ${\rm GL}(V^k)$, it commutes with itself and is therefore abelian. Conversely, if $M_{\widetilde{\varphi}}$ is a torus algebraic group, it is diagonalizable in ${\rm{GL}}(V^k)_{\mathbb C}$. It follows that its commutator in ${\rm{GL}}(V^k)$ contains maximal commutative semi-simple subalgebras $R$ with $[R:{\mathbb Q}]=\dim_{\mathbb Q}V^k$.  Such $R$ give rise to non-trivial endomorphisms in ${\rm{End}}_0(V^k,\widetilde{\varphi})$ corresponding to ``complex multiplications''. For example, let $k=1$, $V^1={\mathbb Q}^{2g}$, and let $\widetilde{\varphi}$ be determined by the complex structure associated to an element $\tau$ in the Siegel upper half space of genus $g$. Let $A$ be the abelian variety with complex points isomorphic to the complex torus ${\mathbb C}^g/{\mathcal {L}}_{\tau}$ where ${\mathcal {L}}_{\tau}={\mathbb Z}^g+\tau{\mathbb Z}^g$. We have ${\rm {End}}_0(A)\simeq {\rm{End}}(V^1,\widetilde{\varphi})$ and the existence of $R$ is equivalent to the usual definition of CM for abelian varieties. For more precise details, see \cite{GGK} (Chapter V), \cite{Mu2}, \cite{ShT}.

 \noindent Consider for example a complex elliptic curve $X$ and generators $\gamma_1$, $\gamma_2$ of the ${\mathbb Z}$-module $H_1(X,{\mathbb Z})$ of rank 2. Then $H^1(X,{\mathbb Q})$ is the dual vector space to $H_1(X,{\mathbb Z})\otimes{\mathbb Q}$, and we let $\gamma_1^\ast$, $\gamma_2^\ast$ be the basis of $H^1(X,{\mathbb Q})$ dual to $\gamma_1$, $\gamma_2$. The Hodge filtration is
 $$
 \{0\}\subset H^{1,0}(X)\subset H^1_\Omega(X,{\mathbb C})
 $$
 and $\dim_{\mathbb C}(H^{1,0}(X))=1$. Let $\omega'$ be a non-zero holomorphic $1$-form on $X$ and denote also by $\omega'$ its class in $H^{1,0}(X)$. Then $H^{1,0}(X)={\mathbb C}\omega'$ and, with respect to $\gamma_1^\ast$, $\gamma_2^\ast$, now viewed as a basis of $H^1(X,{\mathbb C})$, we have
 $$
 \iota_1^{-1}(\omega')=\left(\int_{\gamma_1}\omega'\right)\gamma_1^\ast+\left(\int_{\gamma_2}\omega'\right)\gamma_2^\ast,
 $$
where $\int_{\gamma_i}\omega'\not=0$, $i=1,2$. Replacing $\omega'$ by $\omega=(\int_{\gamma_1}\omega')^{-1}\omega'$ we may assume that the generator of $F^{1,1}=\iota_1^{-1}(H^{1,0}(X))$ is of the form
 $$
 \gamma_1^\ast+\tau\gamma_2^\ast,
 $$
 where $\Im(\tau)\not=0$. We may even suppose $\Im(\tau)>0$, on replacing $\gamma_2$ by $-\gamma_2$ if necessary. We therefore have
 $$
 \{0\}\subset F^{1,1}={\mathbb C}(\gamma_1^\ast+\tau\gamma_2^\ast)\subset H^1(X,{\mathbb C})={\mathbb C}\gamma_1^\ast+{\mathbb C}\gamma_2^\ast.
 $$
 Now, any ${\overline{\mathbb Q}}$-filtration of $H^1(X,{\overline{\mathbb Q}})$ with signature $(2,1)$ is of the form
 $$
 \{0\}\subset F^{1,1}_{\overline{\mathbb Q}}={\overline{\mathbb Q}}(\gamma_1^\ast+\alpha\gamma_2^\ast)
 \subset H^1(X,{\overline{\mathbb Q}})={\overline{\mathbb Q}}\gamma_1^\ast+{\overline{\mathbb Q}}\gamma_2^\ast,
 $$
 for some $\alpha\in {\overline{\mathbb Q}}$. However, for $F^{1,1}_{\overline{\mathbb Q}}$ to be a ${\overline{\mathbb Q}}$-form of $F^{1,1}$ we must have $\tau\in {\overline{\mathbb Q}}$, which is a strong assumption. The ``obvious'' ${\overline{\mathbb Q}}$-form of $F^{1,1}$ is ${\overline{\mathbb Q}}(\gamma_1^\ast+\tau\gamma_2^\ast)$, which, again, is only a subspace of $H^1(X,{\overline{\mathbb Q}})$ if $\tau\in{\overline{\mathbb Q}}$. In other words, if $\tau\not\in{\overline{\mathbb Q}}$, there is no non-zero holomorphic $1$-form on $X$ whose image under $\iota_1^{-1}$ is in a ${\overline{\mathbb Q}}$-vector space $F^{1,1}_{\overline{\mathbb Q}}={\overline{\mathbb Q}}(\gamma_1^\ast+\alpha\gamma_2^\ast)$ with $\alpha\in{\overline{\mathbb Q}}$.

 \noindent We now let $\tau$ be an element of the Siegel upper half space of genus $g$, so that $\tau$ is a symmetric $g\times g$ matrix with positive definite imaginary part. We associate to $\tau$ the principally polarized abelian variety $A$ considered above, whose set of complex points is isomorphic to ${\mathbb C}^g/{\mathcal L}_\tau$, where ${\mathcal L}_\tau ={\mathbb Z}^g+\tau{\mathbb Z}^g$. It is well known that we can choose generators $\gamma_1,\ldots,\gamma_{2g}$ of the ${\mathbb Z}$-module $H_1(A,{\mathbb Z})$, and a basis $\omega_1,\ldots,\omega_g$ of the holomorphic 1-forms on $A$, such that
 $$
{\rm Id}_g=\left(\int_{\gamma_j}\omega_i\right)_{i,j=1,\ldots,g}, \qquad \tau=\left(\int_{\gamma_{g+j}}\omega_i\right)_{i,j=1,\ldots,g},
 $$
 where ${\rm Id}_g$ is the $g\times g$ identity matrix, and the rows of the above matrices are indexed by $i$ and the columns by $j$. Let $\gamma^\ast_1$,\ldots,$\gamma^\ast_{2g}$ be the basis of $H^1(A,{\mathbb Q})$ dual to $\gamma_1,\ldots,\gamma_{2g}$. For $i=1,\ldots,g$, the coordinates of $\omega_i$ with respect to the $\gamma^\ast_j$ are given by the $i$-th row of the $g\times 2g$ matrix $({\rm Id}_g,\tau)$. Explicitly, denoting also by $\omega_i$ the class of $\omega_i$ in $H^{1,0}(A)$, we have
 $$
 \iota_1^{-1}(\omega_i)=\gamma_i^\ast+\sum_{j=1}^g\left(\int_{\gamma_{g+j}}\omega_i\right)\gamma_{g+j}^\ast.
 $$
 Therefore, the Hodge filtration of $H^1(A,{\mathbb Q}_A)_{\mathbb C}$ is defined over ${\overline{\mathbb Q}}$ if and only if the entries in all the rows of $\tau$ are algebraic, that is, if and only if the matrix $\tau$ has all its entries algebraic. We can deal in a similar way with abelian varieties whose polarization is not principal.

 \noindent We assume throughout this paper that our complex manifolds are connected. We now assume in addition that $X$ has the underlying structure of a smooth projective variety $X_{\overline{\mathbb Q}}$ of dimension $n$ defined over ${\overline{\mathbb Q}}$, which we usually denote also by $X$ when the context is clear. We have an isomorphism
 $$\iota_2: H^k_\Omega(X,{\mathbb C})\simeq{\mathbb H}^k(X,\Omega_X^{\bullet}),$$
 where the right hand side is the hypercohomology of the algebraic de-Rham complex over ${\mathbb C}$, see \cite{GrH}, Chapter 3, and \cite{VoiAlg}, \S4. As $X$ is defined over ${\overline{\mathbb Q}}$, the ${\mathbb C}$-vector space ${\mathbb H}^k(X,\Omega_X^\bullet)$ has a natural ${\overline{\mathbb Q}}$-structure ${\mathbb H}_{\overline{\mathbb Q}}^k={\mathbb H}^k(X_{\overline{\mathbb Q}},\Omega_{X/{\overline{\mathbb Q}}}^\bullet)$. This ${\overline{\mathbb Q}}$-structure has nothing to do, however, with $H^k(X,{\overline{\mathbb Q}})$ as defined above. One can nonetheless define the Hodge filtration algebraically as the inverse image under $\iota_2$ of the filtration of ${\mathbb H}^k(X,\Omega_X^\bullet)$ given by
$$
F^p{\mathbb H}^k(X,\Omega_X^\bullet)={\rm Image}\left({\mathbb H}^k(X,\Omega_X^{\ge p})\rightarrow {\mathbb H}^k(X,\Omega_X^\bullet)\right).
$$
This ${\mathbb C}$-vector space clearly has a ${\overline{\mathbb Q}}$-form in ${\mathbb H}_{\overline{\mathbb Q}}$ (see \cite{Voi} for details).

\noindent Returning to the example of an elliptic curve $X$, assume now that the curve is defined over ${\overline{\mathbb Q}}$. Then it is the set of complex points that satisfy an equation of the form
$$
y^2=4x^3-g_2x-g_3,
$$
where $g_2$ and $g_3$ are algebraic numbers with $g_2^3-27g_3^2\not=0$, together with the ``point at infinity'' in ${\mathbb P}_2$.
There are differential forms $\omega$ and $\eta$ on $X$, defined over ${\overline{\mathbb Q}}$, with ${\mathbb H}_{\overline{\mathbb Q}}^1={\overline{\mathbb Q}}\omega+{\overline{\mathbb Q}}\eta$, such that $\omega$ is holomorphic and $\eta$ has poles with zero residues, see \cite{GrH}, Chapter 3. Explicitly $\omega=dx/y$ and $\eta=xdx/y$. By extension of scalars to ${\mathbb C}$, the following ${\overline{\mathbb Q}}$-filtration of ${\mathbb H}^1_{\overline{\mathbb Q}}$ of signature $(2,1)$ induces, via $\iota_2^{-1}$, the Hodge filtration of $H^1_\Omega(X,{\mathbb C})$:
$$
\{0\}\subset F^1{\mathbb H}^1_{\overline{\mathbb Q}}={\overline{\mathbb Q}}\omega\subset {\mathbb H}^1_{\overline{\mathbb Q}}={\overline{\mathbb Q}}\omega+{\overline{\mathbb Q}}\eta.
$$
The base change to the filtration
$$
\{0\}\subset F^1{\mathbb H}^1(X,\Omega^\bullet_X)={\mathbb C}\omega\subset F^0{\mathbb H}^1(X,\Omega_1^\bullet)={\mathbb C}\omega+{\mathbb C}\eta
$$
from the filtration (using $\iota_2\circ\iota_1$)
$$
\{0\}\subset F^{1,1}={\mathbb C}(\gamma_1^\ast+\tau\gamma_2^\ast)\subset F^{0,1}=H^1(X,{\mathbb C})={\mathbb C}\gamma_1^\ast+{\mathbb C}\gamma_2^\ast,
$$
is given by multiplication by $\int_{\gamma_1}\omega$ at the $F^1$-level and by applying the period matrix
$$
\begin{pmatrix}{\int_{\gamma_1}\omega&\int_{\gamma_2}\omega\cr\int_{\gamma_1}\eta&\int_{\gamma_2}\eta}\end{pmatrix}
$$
at the $F^0$-level. In the case of a principally polarized abelian variety $A$ of dimension $g$, defined over ${\overline{\mathbb Q}}$, it is known that there are $g$ algebraic differential forms $\omega_1,\ldots,\omega_g$, defined over ${\overline{\mathbb Q}}$, which also generate $H^{1,0}(A)$, and that the normalized period matrix $\tau$, referred to above, can be written as a matrix quotient $\Omega_1^{-1}\Omega_2$, where $\Omega_1=(\int_{\gamma_j}\omega_i)_{i,j=1,\ldots,g}$, $\Omega_2=(\int_{\gamma_{j+g}}\omega_i)_{i,j=1,\ldots,g}$ with respect to the homology basis introduced earlier.

\noindent The expectation is that, if $X$ is defined over ${\overline{\mathbb Q}}$, the existence of a filtration $F^{\ast,n}_{\overline{\mathbb Q}}$ of $H^n(X,{\overline{\mathbb Q}})$, $n=\dim(X)$, such that $F^{\ast,n}=F_{\overline{\mathbb Q}}\otimes_{\overline{\mathbb Q}}{\mathbb C}$ is the Hodge filtration of $H^n(X,{\mathbb Q}_X)_{\mathbb C}$, forces $F^{\ast,n}$ to have CM, that is forces $M_{F^{\ast,n}}({\mathbb R})$ to be a subgroup of the isotropy group $H_{F^{\ast,n}}$ (see \cite{GGK}). The oldest result in this direction is due to Th. Schneider \cite{Sch}, and is equivalent to the statement that, for an elliptic curve $X$ defined over ${\overline{\mathbb Q}}$, the Hodge filtration of $H^1_\Omega(X,{\mathbb C})$ is defined over ${\overline{\mathbb Q}}$ if and only if $X$ has complex multiplication. When $X({\mathbb C})$ has complex multiplication, it is isomorphic to a quotient ${\mathbb C}/{\mathcal L}_\tau$, with ${\mathcal L}_\tau={\mathbb Z}+\tau{\mathbb Z}$ acting additively and $\tau$ imaginary quadratic. Moreover, the Mumford-Tate group is a maximal torus $T$ in ${\rm GL}(2,{\mathbb Q})$ such that $T({\mathbb Q})\simeq{\mathbb Q}(\tau)^\ast$. Th. Schneider's result is equivalent to the statement that the elliptic modular function takes an algebraic value at an algebraic argument $\tau$ in the upper half plane if and only if $\tau$ is imaginary quadratic. An analogue for abelian varieties is due to P. Tretkoff (a.k.a. Paula B. Cohen) \cite{Co1} and Shiga-Wolfart \cite{SW}, and is joint work of all three authors published in two separate papers.

\noindent In general, however, it is important to stress that, even when $X$ is a smooth projective variety defined over ${\overline{\mathbb Q}}$, the assumption that the Hodge filtration be defined over ${\overline{\mathbb Q}}$ boils down to an intrinsic property of the vector subspaces $F^{p,k}(X)$ of $H^k(X,{\mathbb C})$ and does not necessarily translate in a naive manner to a transcendence statement about suitable quotients of periods of algebraic differential forms on $X$ defined over ${\overline{\mathbb Q}}$. Indeed, as already remarked, there is no direct relation between the underlying ${\overline{\mathbb Q}}$-structure of ${\mathbb H}^k(X,\Omega_X^{\bullet})$ and that of $H^k(X,{\mathbb C})$. Even in the case of an abelian variety discussed above, the relation of $\Omega_1$ and $\Omega_2$ to $\tau$ is via matrix inversion, which is a complicated operation. Moreover, for higher weight, the normalized period matrices can involve combinations of periods of forms of type $(p,q)$ and $(p',q')$ with $(p,q)\not=(p',q')$. Although we do not pursue this in the present paper, for smooth projective variations of Hodge structure, defined over ${\overline{\mathbb Q}}$ in a suitable sense, we  can use the algebraic de-Rham theorem to in general view any transcendence results in terms of suitable combinations of matrices of periods of algebraic differential forms defined over ${\overline{\mathbb Q}}$. For more details, see Proposition 4.4.1 of \cite{CMP}.  All the specific examples we consider are Calabi-Yau manifolds, with the main emphasis in practice being on the one dimensional complex vector space $F^{n,n}_\Omega(X)=H^{n,0}(X)$, $n=\dim(X)$. It is clear that we can, in this case, directly relate the existence or not of a generator of $\iota_1^{-1}(F^{n,n}_\Omega(X))$ in $H^n(X,{\overline{\mathbb Q}})$ to the algebraicity or not of the quotient of any two non-zero periods of a generator of $H^{n,0}(X)$, chosen to be an algebraic differential form on $X$ defined over ${\overline{\mathbb Q}}$. For more details, see \S\ref{s:examples}. For Calabi-Yau manifolds, the following question is of interest. We use the word ``Expectation'', as the conjectures formulated in transcendental number theory are usually very natural, and it is often the proofs that are more interesting.

\medskip

\noindent{\bf Expectation 1:} Let $X$ be a smooth projective variety of dimension $n$ defined over ${\overline{\mathbb Q}}$. If the Hodge filtration $F^{\ast,n}$ of $H^n(X,{\mathbb Q}_X)_{\mathbb C}$ is defined over ${\overline{\mathbb Q}}$ in the sense of Definition 2, then $F^{\ast,n}$ has CM in the sense of Definition 3.

\medskip

\noindent In certain situations, it is of interest to look at the CM properties  of the Hodge filtration of $H^k(X,{\mathbb Q}_X)_{\mathbb C}$ when $k<n$. For example, for abelian varieties of any dimension, the classical notion of complex multiplication in the sense of Shimura-Taniyama pertains to the level 1 Hodge structure. This corresponds to Part 2 of the Problem in \S\ref{s:VHS}, see also the comments in \cite{Roh}, Introduction, p.3. In this paper, we sometimes show that a variant of Expectation 1 holds (see the CMCY property in \cite{Roh}), which can be stated as follows.

\medskip

\noindent {\bf Expectation 2:} Let $X$ be a smooth projective variety of dimension $n$ defined over ${\overline{\mathbb Q}}$. Suppose that, for all $k=0,\ldots,n$, the Hodge filtration $F^{\ast,k}$ of $H^k(X,{\mathbb Q}_X)_{\mathbb C}$ is defined over ${\overline{\mathbb Q}}$ in the sense of Definition 2. Then $F^{\ast,k}$ has CM for all $k=0,\ldots,n$ in the sense of Definition 3.

\medskip

\noindent As in \S\ref{s:VHS} of this paper, the above predictions should be recast in terms of the primitive cohomology of $X$. Nonetheless, in our examples, it is often useful to work with the full cohomology group and establish Expectation 1, and then to deduce the appropriate results for primitive cohomology. Moreover, showing Expectation 2 for the full cohomology groups  implies a similar result for the primitive cohomology groups using the primitive Lefschetz decomposition. Expectations 1 and 2 can be viewed as predictions in transcendence theory. In particular, the proofs in this paper use transcendence techniques and are applications of W\"ustholz's Analytic Subgroup Theorem \cite{Wu3}, which has its classical origins in Baker's method.

\section{Variations of Hodge structure
 and CM}\label{s:VHS}

\noindent The examples we discuss in this paper all arise as fibers of families of smooth projective varieties, where the setting of variations of Hodge structure and period domains is appropriate. Moreover, for the algebraic and arithmetic properties of complex multiplication we should
 work with polarized Hodge structures, and therefore with primitive cohomology and polarized period domains.

\noindent Throughout this section, let $X$ be a smooth connected complex projective variety of complex dimension $n$. Assume that $X$ is embedded in a projective space and is not contained in any hyperplane. Restricting to $X$ the K\"ahler form associated to the Fubini-Study metric on the ambient projective space, we obtain a K\"ahler form $\omega$ whose class in $H^2(X,{\mathbb R})$, also denoted $\omega$, is in fact in $H^2(X,{\mathbb Q})$. For $k=0,\ldots,n$, the primitive cohomology, with coefficients in ${\mathbb Q}$, denoted $P^k=P^k_{\mathbb Q}=H^k(X,{\mathbb Q})_{\rm{prim}}$, is the kernel of the Lefschetz operator
$$
L^{n-k+1}: H^k(X,{\mathbb Q})\rightarrow H^{2n-k+2}(X,{\mathbb Q})
$$
$$
u\mapsto (\wedge^{n-k+1}\omega)\wedge u,
$$
where $L$ is the cup product with $\omega$, which corresponds to the wedge product with $\omega$ viewed as a representative of a class in $H^2_\Omega(X,{\mathbb R})$.
Let $P^k_K=P^k\otimes_{\mathbb Q}K$ for any field $K$ containing ${\mathbb Q}$ and for any $0\le k\le n$. As $L$ shifts the Hodge decomposition by $(1,1)$, the primitive cohomology $P^k$ has an induced Hodge decomposition
$$
P^k_{\mathbb C}=H^k(X,{\mathbb C})_{\rm prim}=\bigoplus_{p+q=k}P^{p,q}
$$
where
$$
P^{p,q}=H^{p,q}_{\rm prim}=H^{p,q}\cap H^k(X,{\mathbb C})_{\rm prim}
$$
We denote by $P^k(X)=H^k(X,{\mathbb Q}_X)_{\rm prim}$ the associated rational Hodge substructure of $H^k(X,{\mathbb Q}_X)$.

The filtration $F^{\ast,k}=\iota_1^{-1}(F^\ast_\Omega)$ given by the Hodge filtration induces, by restriction to $P^k_{\mathbb C}$, a filtration of $P^k_{\mathbb C}$ which we denote by $F_P^{\ast,k}= F_P^{\ast,k}(X)$. We call it the Hodge filtration of $H^k(X,{\mathbb Q}_X)_{{\rm prim},\mathbb C}$ to emphasize its dependence on the complex manifold structure of $X$. We denote by $f_P^{\ast,k}$ the signature of $F_P^{\ast,k}$. As in \S\ref{s:Schneider}, we say that $F_P^{\ast,k}$ is defined over ${\overline{\mathbb Q}}$ if it is induced by extension of scalars to ${\mathbb C}$ by a ${\overline{\mathbb Q}}$-filtration of $P^k_{\overline{\mathbb Q}}$. Let $h_P^{p,q}=\dim(F_P^p\cap {\overline{F_P}}^q)$, then $h_P^{p,q}=\dim(P^{p,q})$ and $f_P^{p,k}=\sum_{p'\ge p}h_P^{p',k-p'}$, $p=0,\ldots,k$, and $p+q=k$.

\noindent The homomorphism $\widetilde{\varphi}$ of \S\ref{s:Schneider}, equivalent to the Hodge filtration, preserves $P^k_{\mathbb R}$, and we denote by $\widetilde{\varphi}_P$ the induced map of ${\mathbb R}$-algebraic groups
$$
\widetilde{\varphi}_P:{\mathbb S}({\mathbb R})\rightarrow {\rm GL}(P^k)_{\mathbb R}
$$
with ${\widetilde{\varphi}}_P({\rm{diag}}(r,r))=r^k{\rm{Id}}_{P^k}$, for $r\in{\mathbb R}$, $r\not=0$. Let $\varphi_P$ be the restriction of $\widetilde{\varphi}_P$ to ${\mathbb U}({\mathbb R})$, where ${\mathbb U}$ is the group with $K$-points
$$
{\mathbb U}(K)=\left\{\begin{pmatrix}{a&-b\cr b&a}\end{pmatrix}:\, a^2+b^2=1,\, a,b\in K\right\}.
$$
Then $\varphi_P$ has image in ${\rm{SL}}(P^k)_{\mathbb R}$, and giving $\varphi_P$ together with the weight $k$ is equivalent to giving ${\widetilde{\varphi}}_P$ and also equivalent to the Hodge filtration $F^{\ast,k}_P$.
The endomorphism $C=\varphi_P(\sqrt{-1})$ is called the Weil operator. It acts on $P^{p,q}$ by multiplication by $\sqrt{-1}^{\;p-q}$, preserving the real vector space $(V^{p,q}\oplus V^{q,p})\cap V_{\mathbb R}$.

\medskip

\noindent We drop the subscript $P$ in all the above notations when it is clear we are working with primitive cohomology.

\medskip

\noindent The polarization $Q=Q(\;\cdot\;,\;\cdot\;)$ of $(P^k,\varphi_P)=(P^k,F_P^{\ast,k})$, is defined to be the bilinear non-degenerate map defined by the \emph{Hodge-Riemann} (HR) pairing
$$
Q(u,v):=\,(-1)^{\frac{k(k-1)}2}\,\int_X(\wedge^{n-k}\omega)\wedge u\wedge v,\qquad u,v\in P^k.
$$
It satisfies
\begin{equation}\label{eq:sym} Q(u,v)=(-1)^kQ(v,u)\end{equation}
together with the Hodge-Riemann relations
$$
Q(F^p,F^{k-p+1})=0, \qquad ({\rm{HR1}}),
$$
$$
Q(u,C\,\overline{u})>0, \quad  u\not=0,\quad u\in P^k_{\mathbb C},\qquad ({\rm{HR2}}).
$$
We refer to $(P^k,Q,\varphi_P)=(P^k,Q,F^{\ast,k}_P)$ as a polarized Hodge structure. An isomorphism of Hodge structures that preserves the given polarizations is known as a Hodge isometry. The Hodge-Riemann relations imply that $Q(P^{p,q},P^{p',q'})=0$ for $(p,q)\not=(q',p')$ and that $\sqrt{-1}^{p-q}Q(v,{\overline{v}})>0$, for $v\in P^{p,q}$. The above definition of $Q$ extends as is to all of $H^k(X,{\mathbb Q})$ where it gives a non-degenerate bilinear map, also denoted by $Q$. For every $k=0,\ldots,n$, we have the \emph{primitive Lefschetz decomposition},
$$
H^k(X,{\mathbb Q})=\,H^k(X,{\mathbb Q})_{\rm prim}\bigoplus_{2\le 2r\le k} L^rH^{k-2r}(X,{\mathbb Q})_{\rm prim},
$$
where the decomposition is orthogonal with respect $Q$, see for example \cite{Cat}, Theorem 7.3.4.  To prove the above decomposition, one shows that $\sum_{2r\le k}L^r$ is an isomorphism of ${\mathbb Q}$-vector spaces from $\oplus_{2r\le k}P^{k-2r}$ to $H^k(X,{\mathbb Q})$. In particular, this shows that the Hodge filtration of $H^k(X,{\mathbb Q}_X)_{\mathbb C}$ is defined over ${\overline{\mathbb Q}}$ if and only if the Hodge filtrations of $H^{k-2r}(X, {\mathbb Q}_X)_{\rm prim, {\mathbb C}}$, for $2r\le k$, are all defined over ${\overline{\mathbb Q}}$.

\noindent On $H^k(X,{\mathbb C})$, the hermitian form (sesquilinear and with $h(u,v)={\overline{h(v,u)}}$) given by $h(u,v)=(\sqrt{-1})^kQ(u,\overline{v})$ is definite ($h(u,u)\not=0$) of sign $(-1)^q$ on $P^{p,q}=H^{p,q}_{\rm prim}$, $p+q=k$, and of sign $(-1)^{q+r}$ on $L^rH^{p,q}_{\rm prim}$, where $2r+p+q=k$, $2\le 2r\le k$. Moreover, the Hodge decomposition $H^k(X,{\mathbb C})=\oplus_{p+q=k}H^{p,q}$ is orthogonal with respect to $h$, so that $h(H^{p,q},H^{p',q'})=0$ when $(p,q)\not=(p',q')$.

\noindent As $H^k(X,{\overline{\mathbb Q}})$ is closed under complex conjugation and $h(u,v)\in{\overline{\mathbb Q}}$ for $u,v$ in $H^k(X,{\overline{\mathbb Q}})$, this hermitian form restricts to a ${\overline{\mathbb Q}}$-valued hermitian form on $H^k(X,{\overline{\mathbb Q}})$. By \S\ref{s:lemmas}, Lemma 4, if the Hodge filtration $F^{\ast,k}$ of $H^k(X,{\mathbb Q}_X)_{\mathbb C}$ is defined over ${\overline{\mathbb Q}}$, then so is the Hodge filtration of $H^k(X,{\mathbb Q}_X)_{{\rm prim},{\mathbb C}}$ and $L^rH^{k-2r}(X,{\mathbb Q}_X)_{{\rm prim},{\mathbb C}}$, $2\le 2r\le k$. From \S\ref{s:lemmas}, Lemma 10, applied to these Hodge structures, it then follows that each $L^rH^{p,q}(X)_{\rm prim}$, $p+q=k$, has a basis in $L^rH^{p,q}(X)_{\rm prim}\cap H^k(X,{\overline{\mathbb Q}})$, $0\le 2r\le k$, and therefore that each $H^{p,q}(X)$, $p+q=k$, has a basis in $H^{p,q}(X)\cap H^k(X,{\overline{\mathbb Q}})$.

\noindent For $i=1,2$, let $(V_i,{\widetilde{\varphi}}_i)$ be a rational Hodge structure of weight $k_i$ and $Q_i$ a non-degenerate bilinear form on $V_i$ such that $Q_i(v_i,v_i')=(-1)^{k_i}Q_i(v_i',v_i)$, $v_i, v_i'\in V_i$. Let $h_i$ be the hermitian form on $V_{i,{\mathbb C}}$ given by
$$
h_i(u,v)=(\sqrt{-1})^{k_i}Q_i(u,{\overline{v}})
$$
and suppose that the Hodge decomposition of $V_{i,{\mathbb C}}$ is orthogonal with respect to $h_i$. Suppose in addition that each $(V_i,{\widetilde{\varphi}}_i)$ is a direct sum, orthogonal with respect to $Q_i$, of the form
$$
(V_i,{\widetilde{\varphi}}_i)=\oplus_{r_i}(V_{r_i},{\widetilde{\varphi}}_{r_i}),
$$
where each $(V_{r_i},{\widetilde{\varphi}}_{r_i})$ is a rational Hodge substructure of weight $k_i$ such that $h_i$ is definite on each summand $(V_{r_i})^{p,q}=(V_{r_i}\cap V_i^{p,q})$, $p+q=k_i$, of the Hodge decomposition of $(V_{r_i})_{\mathbb C}$.
Let
$$
(V,{\widetilde{\varphi}})=(V_1\otimes_{\mathbb Q}V_2,{\widetilde{\varphi_1}}\otimes{\widetilde{\varphi_2}})
$$
be the rational tensor product Hodge structure of weight $k=k_1+k_2$. Let $Q$ be the bilinear form on $V$ defined on elementary tensors by
$$
Q(v_1\otimes_{\mathbb Q}v_2,v_1'\otimes_{\mathbb Q}v_2')=Q_1(v_1,v_1')Q_2(v_2,v_2'),\qquad v_1,v_1'\in V_1,\;v_2,v_2'\in V_2,
$$
and extended to all of $V$ by ${\mathbb Q}$-bilinearity. Then $Q$ is non-degenerate and we have $Q(u,v)=(-1)^kQ(v,u)$, $u,v\in V$. Let $h$ be the hermitian form on $V_{\mathbb C}$ given by $h(u,v)=(\sqrt{-1})^kQ(u,{\overline{v}})$. The $(p,q)$-part, $p+q=k$, of the Hodge decomposition of $V_{\mathbb C}$ is of the form
$$
V_{\mathbb C}^{p,q}=\oplus_{i+j=p}\left(V_1^{i,k_1-i}\otimes_{\mathbb C}V_2^{j,k_2-j}\right).
$$
By first checking on elementary tensors, it is easy to see that
$$
h(V_1^{i,k_1-i}\otimes_{\mathbb C}V_2^{j,k_2-j},V_1^{i',k_1-i'}\otimes_{\mathbb C}V_2^{j',k_2-j'})=0,
$$
unless $i=i'$ and $j=j'$. Therefore, the $V_1^{i,k_1-i}\otimes_{\mathbb C}V_2^{j,k_2-j}$ are mutually orthogonal with respect to $h$. We have the decomposition of rational Hodge structures,
$$
(V,{\widetilde{\varphi}})=\oplus_{r_1,r_2}(V_{r_1}\otimes_{\mathbb Q}V_{r_2},{\widetilde{\varphi}}_{r_1}\otimes {\widetilde{\varphi}}_{r_2}),
$$
orthogonal with respect to $Q$, and $h$ is definite on each summand
$$
(V_1^{i,k_1-i}\cap V_{r_1})\otimes_{\mathbb C}(V_2^{j,k_2-j}\cap V_{r_2}),\qquad i+j=p
$$
of $((V_{r_1})_{\mathbb C}\otimes_{\mathbb C}(V_{r_2})_{\mathbb C})^{p,q}$, $p+q=k$.

\noindent Returning to our discussion of primitive cohomology, the Mumford-Tate group $M_\varphi$ of the Hodge structure $(P^k,\varphi)$ is the smallest ${\mathbb Q}$-algebraic subgroup of ${\rm{SL}}(V)$ whose real points contain $\varphi({\mathbb U}({\mathbb R}))$. Here, we use the terminology of \cite{GGK}, rather than calling this the Hodge group or special Mumford-Tate group. The Mumford-Tate group $M_{\widetilde{\varphi}}$ as defined in \S\ref{s:Schneider}, but for $(P^k,{\widetilde{\varphi}})$, is abelian if and only if $M_{\varphi}$ is abelian, in which case they are both algebraic tori. The associated filtration $F_P^{\ast,k}$ of $P^k_{\mathbb C}$ then has CM in the sense that $M_{F^{\ast,k}_P}:=M_{\varphi}$ is contained in the stabilizer $H_{F^{\ast,k}_P}$ of $F^{\ast,k}_P$ in $G({\mathbb R})$, where $G={\rm Aut}(P^k,Q)$ is the ${\mathbb Q}$-algebraic subgroup of ${\rm SL}(P^k)$ whose elements preserve $Q$ . This leads to Definition 4 below.

\noindent First, define a Calabi-Yau variety $X$ of dimension $n$ to be a smooth connected projective variety of dimension $n$ with $h^{j,0}(X)=0$, for $j=1,\ldots, n-1$, and with a nowhere vanishing holomorphic $n$-form, so that $h^{n,0}=1$.

\medskip

\noindent {\bf Definition 4:} Let $X$ be a smooth projective variety of dimension $n$ and $k$ an integer with $0\le k\le n$. Let $F^{\ast,k}_P$ be the Hodge filtration of $H^k(X,{\mathbb Q}_X)_{\rm prim,{\mathbb C}}$. Let $M_{F^{\ast,k}_P}$ be the Mumford-Tate group of $F^{\ast,k}_P$ and $H_{F^{\ast,k}_P}$ be the stabilizer of $F^{\ast,k}_P$ in $G({\mathbb R})$. Then $F^{\ast,k}_P$ has complex multiplication (CM) if and only if $M_{F^{\ast,k}_P}({\mathbb R})$  is contained in  $H_{F^{\ast,k}_P}$. We also say that $H^k(X,{\mathbb Q}_X)_{\rm prim}$ has CM. If $X$ is a Calabi-Yau variety, we say that $X$ has CM if $F^{\ast,n}_P$ has CM and that $X$ has strong CM if $F^{\ast,k}_P$ has CM for $k=0,\ldots,n$.

\medskip

\noindent By our remarks above, for an abelian variety $X=A$, the usual definition of CM, in the accepted sense of Shimura and Taniyama, is that $A$ has CM if and only if $H^1(A,{\mathbb Q}_A)=H^1(A,{\mathbb Q}_A)_{\rm prim}$ has CM, even when $\dim(A)>1$. Therefore, to avoid ambiguity, we should always mention the rational Hodge structure, and not just the variety, for the CM property. Definition 4 is sufficient for the examples of Calabi-Yau varieties in this paper, see \S\ref{s:BVtower} for more details.

\medskip

\noindent When $X$ is a curve, we clearly have $H^1(X,{\mathbb Q})=H^1(X,{\mathbb Q})_{\rm prim}$. When $X$ is a surface, the Hodge structures $H^2(X,{\mathbb Q}_X)$ and $H^2(X,{\mathbb Q}_X)_{\rm prim}$ can only differ by a $(1,1)$-part, and the Mumford-Tate group of a Hodge structure concentrated in bidegree $(p,p)$, for any $p$, is trivial. Indeed, when the surface $X$ is connected, if $\omega$ is the K\"ahler form used to define the polarization, we have $H^2(X,{\mathbb Q}_X)=H^2(X,{\mathbb Q}_X)_{\rm prim}\oplus{\mathbb Q}\,\omega$. For a Calabi--Yau 3-fold, we have $H^1(X,{\mathbb Q})=\{0\}$, so that $H^3(X,{\mathbb Q})_{\rm{prim}}=H^3(X,{\mathbb Q})$. We see by the primitive Lefschetz decomposition that the rational Hodge structures $H^k(X,{\mathbb Q}_X)_{\rm prim}$  have CM for all $k=0,\ldots,n$ if and only if the $H^k(X,{\mathbb Q}_X)$ have CM for all $k=0,\ldots,n$. In \S\ref{s:BVtower}, we are interested in showing that certain Calabi-Yau varieties have strong CM, so we can work equally well with the full cohomology groups. In \S\ref{s:expect1}, our examples are either curves, surfaces, or Calabi-Yau 3-folds, so that, again, we may work with the full cohomology group to show that the varieties have CM.

\medskip

\noindent All our examples naturally occur in families of smooth projective varieties parametrized by a quasi-projective variety defined over ${\overline{\mathbb Q}}$. We now fix a smooth projective variety $X_b$, where $b$ stands for base point, of dimension $n$, with primitive cohomology $P^k=H^k(X_b,{\mathbb Q})_{\rm prim}$, $k=0,\ldots,n$. Consider the associated Hodge numbers $h_P^{p,q}=\dim_{\mathbb C}P^{p,q}$, $p+q=k$, and the polarization $Q$, as above. Therefore $h_P^{p,q}=h_P^{q,p}$ and $\sum_{p=0}^kh^{p,k-p}=\dim_{\mathbb C}P^k_{\mathbb C}$. Denote by $F^{\ast,k}_P=F^{\ast,k}_P(X_b)$ the Hodge filtration of $H^k(X_b,{\mathbb Q}_{X_b})_{{\rm prim},{\mathbb C}}$. In particular $F_P^{\ast,k}$ satisfies (HR1) and (HR2), so that $F^{p,k}_P\oplus{\overline{F}}^{k-p+1}_P=P^k_{\mathbb C}$. Moreover $F_P^{\ast,k}$ has signature $f^{\ast,k}_P$, where $f^{p,k}_P=\sum_{p'\ge p}h^{p',k-p'}$, $p=0,\ldots,k$. We define the period domain $D^k=D(P^k,Q,f^{\ast,k}_P)$ to be the set of all the ${\mathbb C}$-filtrations $F^{\ast,k}$ of $P^k_{\mathbb C}$,
$$F^{k,k}\subset F^{k-1,k}\subset\ldots\subset F^{1,k}\subset F^{0,k}=P^k_{\mathbb C}$$
of signature $f^{\ast,k}_P$, which satisfy the Hodge-Riemann relations (HR1) and (HR2) for $Q$. Let $G={\rm{Aut}}(P^k,Q)$, and denote by $G(K)$ the $K$-points of $G$, for any field $K$ containing ${\mathbb Q}$. The period domain is a homogeneous space. If we fix a ${\mathbb C}$-filtration $F^{\ast,k}_0$ of $P^k_{\mathbb C}$ in $D^k$ with isotropy group $H_0=H_0({\mathbb R})$ in $G({\mathbb R})$, then $D^k=G({\mathbb R})/H_0$.

 \noindent  The compact dual ${\check{D}}^k$ of $D^k$ is defined to be the set of all the ${\mathbb C}$-filtrations
 $$
 F^{k,k}\subset F^{k-1,k}\subset \ldots \subset F^{0,k}=P^k_{\mathbb C},
 $$
of signature $f^{\ast,k}_P$ satisfying (HR1). Therefore $D^k\subset{\check{D}}^k$. The complex Lie group $G({\mathbb C})$ acts transitively on ${\check{D}}^k$. The isotropy group in $G({\mathbb C})$ of any filtration $F_0^{\ast,k}$ in ${\check{D}}^k$ is a parabolic subgroup $P_0({\mathbb C})$ of $G({\mathbb C})$, and
$$
{\check{D}}^k=G({\mathbb C})/P_0({\mathbb C}).
$$
It is well-known that the compact dual can be embedded in a projective space as a subvariety defined over ${\mathbb Q}$ (see for example \cite{GGK}), and that $F^{\ast,k}\in{\check{D}}^k({\overline{\mathbb Q}})$ if and only if $F^{\ast,k}$ is defined over ${\overline{\mathbb Q}}$ in the sense we have already defined, namely that it is induced by extension of scalars to ${\mathbb C}$ from a ${\overline{\mathbb Q}}$-filtration of $P^k_{\overline{\mathbb Q}}$. In other words, there is a basis of each vector space $F^{p,k}$, $p=0,\ldots,k$, with elements in $P^k_{\overline{\mathbb Q}}$ (see \cite{Sha1}, Chapter I, \S4, Example 1). The same proof that shows $G({\mathbb C})$ acts transitively on ${\check{D}}^k({\mathbb C})$ can be adapted to show that $G({\overline{\mathbb Q}})$ acts transitively on ${\check{D}}^k({\overline{\mathbb Q}})$.

\noindent Consider a smooth complex projective family ${\mathcal X}\subset {\mathbb P}_N$, for some $N$,
$$
\pi:{\mathcal X}\rightarrow S.
$$
In particular, the map $\pi$ is a surjective, proper, holomorphic submersion of complex manifolds. For $s\in S$, let ${\mathcal X}_s$ be the corresponding fiber of ${\mathcal X}$ at $s$, and let $\dim({\mathcal X}_s)=n$ (independent of $s$). Assume the family has an underlying structure over ${\overline{\mathbb Q}}$. In particular, the base $S$ is a smooth quasi-projective variety defined over ${\overline{\mathbb Q}}$ such that ${\mathcal X}_s$, $s\in S({\overline{\mathbb Q}})$, is a smooth projective variety defined over ${\overline{\mathbb Q}}$. Fix a base point $b$ in $S$ and let  $P^k=H^k({\mathcal X}_b,{\mathbb Q})_{\rm{prim}}$ with the polarization $Q$ defined above. Let $f^{\ast,k}_P=f^{\ast,k}_P({\mathcal X}_b)$ be the signature of the Hodge filtration of $P^k_{\mathbb C}$ associated to the usual complex structure on ${\mathcal X}_b$. There is a lattice $P^k_{\mathbb Z}$ in $P^k$, namely a ${\mathbb Z}$-module of rank $\dim_{\mathbb Q}P^k$, which we fix, such that $P^k=P^k_{\mathbb Z}\otimes_{\mathbb Z}{\mathbb Q}$. The lattice $P^k_{\mathbb Z}$ comes from the torsion free part of the singular cohomology $H^k({\mathcal X}_b,{\mathbb Z})$. Let $G({\mathbb Z})$ be the set of $Q$-isometries of $P^k_{\mathbb Z}$ (this will involve a choice of ${\mathbb Z}$-basis of $P_{\mathbb Z}^k$ if we wish to represent $G({\mathbb Z})$ as a ring of matrices.) Then $G({\mathbb Z})$ is an arithmetic subgroup of $G({\mathbb Q})$ with a proper discontinuous action on $D^k$ and ${\check{D}}^k$. For $s\in S$, the Hodge filtration of $H^k({\mathcal X}_s,{\mathbb Q}_{{\mathcal X}_s})_{{\rm prim},{\mathbb C}}$, given by the complex structure on ${\mathcal X}_s$, depends on $s$, but nonetheless can be pulled back to a filtration of $P^k_{\mathbb C}$ with signature $f^{\ast,k}_P=f^{\ast,k}_P({\mathcal X}_b)$ (independent of $s$). Let $(P^k,Q,F^{\ast,k}_s)$ be
the corresponding polarized Hodge structure, or filtration, on $P^k$. For $k\le n$, the induced map from $S$ to the corresponding period domain $D^k$ is multi-valued when $S$ has non-trivial fundamental group, but its image in $\Gamma\backslash D^k$ is well-defined, where $\Gamma\subseteq G({\mathbb Z})$ is the image of the monodromy representation (for a rigorous treatment and full details, see \cite{CMP}, Chapter 4). Therefore, we have a well-defined period map
$$
\Phi_k: S\rightarrow \Gamma\backslash D^k.
$$
Let $\rho_k:D^k\rightarrow \Gamma\backslash D^k$ be the natural projection.

\medskip

\noindent We now state three variants of a Problem whose solution would yield analogues of Th. Schneider's result on the transcendence of the values of the elliptic modular function at algebraic arguments which are not quadratic imaginary. Part (a) is a special case of Conjecture (VIII.A.8) of \cite{GGK}. For the Hodge filtration of primitive cohomology, part 1 corresponds to the notion of CM and part 2 to the notion of strong CM.

\medskip

\goodbreak

\noindent {\bf Problem:} Let $s\in S({\overline{\mathbb Q}})$.

\begin{enumerate}

\item Suppose the filtration $F^{\ast,n}\in D^n$ satisfies $\rho_n(F^{\ast,n})=\Phi_n(s)$. Show that $F^{\ast,n}\in{\check{D}}^n({\overline{\mathbb Q}})$ if and only if $F^{\ast,n}$ has CM.

\item Suppose the filtration $F^{\ast,k}\in D^k$ satisfies $\rho_k(F^{\ast,k})=\Phi_k(s)$, for all $k\le n$. Show that $F^{\ast,k}\in{\check{D}}^k({\overline{\mathbb Q}})$, for all $k\le n$, if and only if $F^{\ast,k}$ has CM, for all $k\le n$.

\item Suppose the filtration $F^{\ast,k}\in D^k$ satisfies $\rho_k(F^{\ast,k})=\Phi_k(s)$, for some $k\le n$. Show that $F^{\ast,k}\in{\check{D}}^k({\overline{\mathbb Q}})$ if and only if $F^{\ast,k}$ has CM.

\end{enumerate}
\bigskip

\noindent The ``if'' part of the above statement is immediate in all the examples we consider, the only work being in the ``only if'' part.

\noindent Notice that, by the results of \cite{Co1}, \cite{SW}, and the remarks made in \S\ref{s:Schneider}, when $A$ is an abelian variety defined over ${\overline{\mathbb Q}}$, and $F^{\ast,1}$ is the Hodge filtration of $H^1(A,{\mathbb C})$, then $A$ will have CM, in the classical sense, if and only if $F^{\ast,1}$ is defined over ${\overline{\mathbb Q}}$, even when $\dim(A)>1$. Therefore, it can happen that we are interested in $F^{\ast,k}$ for one $k<n$, which corresponds to the Problem, Part 3. Indeed, the definition of CM for an abelian variety refers to the rational Hodge structure $H^1(A,{\mathbb Q}_A)$. For the examples of this paper, we shall mainly be interested in Parts 1 and 2 of the Problem. As a consequence of the result on abelian varieties, and using the well-known description of the Siegel upper half space ${\mathcal H}_g$ of genus $g$ in terms of complex structures on ${\mathbb R}^{2g}$ (see for example \cite{Run}, \S3), we have by \cite{Co1}, \cite{SW}:

\bigskip

\noindent {\bf Proposition 1:} Let $\pi:{\mathcal X}\rightarrow S$ be a family of smooth projective algebraic curves of genus $g$ satisfying the above assumptions. Then, we may take $D={\mathcal H}_g$ and $\Gamma\subseteq {\rm{PSp}}(2g,{\mathbb Z})$, and the statement of the Problem is true.

\bigskip

\noindent In \cite{TCY}, we proved the statement of the Problem, part 2, for two families of Calabi-Yau 3-folds with dense sets of CM fibers, studied respectively by Borcea \cite{Bor}, Viehweg-Zuo \cite{VZ}, and Rohde \cite{Roh}, and the first step in a tower of Calabi--Yau manifolds that starts with these examples. In the present paper, we extend these results to more of the examples of Calabi-Yau families with a dense set of CM fibers appearing in \cite{Roh}. In \S\ref{s:BVtower}, we show that the statement of the Problem, part 2, holds for all families of Calabi-Yau varieties with involution in a Borcea-Voisin tower, as long as any new Calabi-Yau family introduced at any step, as well as the ramification loci of its involution, already satisfy the Problem, part 2. In \S\ref{s:expect1}, we give examples of Viehweg-Zuo towers, or families of iterated cyclic covers of algebraic curves, for which the statement of the Problem, part 1, holds. In \S\ref{s:examples}, we approach the same questions by the explicit computation of points in period domains.

\bigskip

\section{The general step up a Borcea-Voisin tower}\label{s:BVtower}

\bigskip

\noindent In this section, we consider the general inductive step in the construction of a Borcea-Voisin tower of Calabi-Yau manifolds with involution, as in \cite{BorK3}, \cite{VoiK3}, and as further developed in \cite{Roh}. We define a Calabi-Yau manifold $X$ of dimension $n$ to be a complex connected K\"ahler manifold of dimension $n$ with $h^{j,0}(X)=0$, for $j=1,\ldots, n-1$, and with a nowhere vanishing holomorphic $n$-form, so that $h^{n,0}=1$.  Note that this is a more restrictive definition of Calabi-Yau than that found in some other references. An abelian surface, for example, is often considered to be a Calabi-Yau manifold, but has $h^{1,0}=2$. For a connected complex manifold, we always have $h^{0,0}=1$.

\noindent In the rest of the paper, the complex manifolds we consider are compact K\"ahler manifolds and all maps between them are assumed holomorphic. We assume in addition that these manifolds have the underlying structure of smooth complex projective varieties, so that such maps are also algebraic morphisms by Chow's Theorem (see \cite{GrH}, Chapter 1, \S3). It may happen that our Calabi-Yau varieties and the algebraic morphisms between them can be defined over ${\overline{\mathbb Q}}$. We always assume that ${\overline{\mathbb Q}}$ is a subfield of ${\mathbb C}$, by the choice of a suitable embedding. We work with the strong CM definition of \cite{Roh}, which is restated in \S\ref{s:VHS}, Definition 4, of this paper. Correspondingly, we construct varieties $X$ such that the Hodge filtration of $H^k(X,{\mathbb Q}_X)_{\mathbb C}$ is defined over ${\overline{\mathbb Q}}$ for all $k\le \dim(X)$. As remarked already in \S\ref{s:VHS}, this is equivalent to the Hodge filtration of $H^k(X,{\mathbb Q}_X)_{\rm prim, {\mathbb C}}$ being defined over ${\overline{\mathbb Q}}$ for all $k\le \dim(X)$. In the course of the proofs of this section we do need to work with the polarization on primitive cohomology defined in \S\ref{s:VHS}, but mostly we use its extension to the full cohomology. We show that, if the statement of the Problem, part 2, of \S\ref{s:VHS} holds for two initial families of Calabi-Yau varieties with involution, and for the family of ramification loci of that involution, then it holds for the family of Calabi-Yau varieties with involution in the next step of the tower. This amounts to showing Expectation 2 of \S\ref{s:Schneider} for the fibers over the algebraic points of the parameter space.

\noindent We first describe the general step in a Borcea-Voisin tower of Calabi-Yau varieties with involution. For $i=1,2$, let $A_i$ be a Calabi-Yau variety of dimension $d_i$ over the field $K= {\overline{\mathbb Q}}$ or ${\mathbb C}$. Suppose that each $A_i$ carries an involution $I_i$, that is a smooth holomorphic map of order 2 of the associated complex variety $A_i^{an}({\mathbb C})$ into itself. Suppose that as an algebraic morphism this involution can be defined over $K$, giving a map from $X(K)$ to itself of order 2. Suppose that the ramification locus $R_i$ of the map $A_i\rightarrow A_i/I_i$ consists of a union of smooth disjoint hypersurfaces of $A_i$ defined over $K$. Consider the projective variety $\widehat{Y}$ given by the blow-up of the fixed locus of the product involution $I_{1,2}=I_1\times I_2$ on the product variety $Y=A_1\times A_2$. By the universal property of blowing up, the variety $\widehat{Y}$ carries an involution $\widehat{I}_{1,2}$, induced by $I_{1,2}$, whose ramification locus is the exceptional set of the blow-up (see \cite{Ha}, II, Corollary 7.15). Therefore the variety $B=\widehat{Y}/\widehat{I}_{1,2}$ is a smooth projective variety of dimension $d=d_1+d_2$ over $K$, and by \cite{Roh}, Lemma 7.2.4 and Proposition 7.2.5, pp.148-149, the associated complex manifold is Calabi-Yau. The Calabi-Yau variety $B$ carries an involution $I$ such that the ramification locus $R$ of the map $B\rightarrow B/I$ consists of a union of smooth disjoint hypersurfaces of $B$. Equivalently, we obtain a variety isomorphic to $B$ by blowing up the singular locus of $(A_1\times A_2)/(I_1\times I_2)$. On this last quotient, the involutions $I_1\times {\rm{Id}}$ and ${\rm{Id}}\times I_2$ are equivalent, and induce the involution $I$.

\medskip

\noindent {\bf Proposition 2:} In the above situation, suppose that, for all $k\le\dim(B)$, the Hodge filtration of $H^k(B,{\mathbb Q}_B)_{\mathbb C}$ is defined over ${\overline{\mathbb Q}}$.  Then, for $i=1,2$, the Hodge filtration of $H^k(A_i,{\mathbb Q}_{A_i})_{\mathbb C}$  is defined over ${\overline{\mathbb Q}}$ for all $k\le\dim(A_i)$ and the Hodge filtration of $H^k(R_i,{\mathbb Q}_{R_i})_{\mathbb C}$ is defined over ${\overline{\mathbb Q}}$ for all $k\le\dim(R_i)$.

\medskip

\noindent {\bf Proof:} For $X$ a smooth complex projective variety, recall from \S\ref{s:Schneider} that we denote by $F^{\ast,\;k}_\Omega(X)$ the Hodge filtration of $H^k_\Omega(X,{\mathbb C})$ and by $F^{\ast,k}(X)$ the Hodge filtration of $H^k(X,{\mathbb Q}_X)_{\mathbb C}$, which is the inverse image of $F^{\ast,\;k}_\Omega(X)$ under $\iota_k$, the de-Rham isomorphism with complex coefficients. Moreover, by Definition 2, \S\ref{s:Schneider}, we say that $F^{\ast,k}(X)$ is defined over ${\overline{\mathbb Q}}$ if it is induced by extension of scalars to ${\mathbb C}$ by a ${\overline{\mathbb Q}}$-filtration of $H^k(X,{\overline{\mathbb Q}})$. Recall that this bears no relation to the field of definition of the projective variety $X$.

\noindent In what follows, we will encounter various finite dimensional ${\mathbb Q}$-vector spaces $V$ endowed with a ${\mathbb Q}$-linear involution $I$, so that $I^2=1$. Any such involution $I$ can be diagonalized over ${\mathbb Q}$ and has eigenvalues $\pm1$. We denote by $V^+$ the $+1$-eigenspace of $I$, or $I$-invariant subspace of $V$, and by $V^-$ the $-1$ eigenspace of $I$, or $I$-anti-invariant subspace of $V$, without necessarily  specifying the involution. When we do want to stress the involution giving rise to the $I$-invariant subspace of $V$, we will denote that subspace by $V^I$. Suppose $V$ carries a Hodge structure of weight $k$, say, so that $V_{\mathbb C}=\oplus_{p+q=k}V^{p,q}$ with $V^{p,q}={\overline{V}}^{q,p}$. Let $I_{\mathbb C}=I\otimes_{\mathbb Q}{\rm {Id}}_{\mathbb C}$ be the complex extension of $I$ to a ${\mathbb C}$-linear involution on $V_{\mathbb C}$. Then ${\overline{I_{\mathbb C}(v)}}={I_{\mathbb C}}({\overline{v}})$, for all $v\in V$. The involutions $I$ we consider will always respect the Hodge structure, in that $I_{\mathbb C}(v_{p,q})\in V^{p,q}$ for all $v_{p,q}\in V^{p,q}$ and all $p,q$ with $p+q=k$. Therefore, the involution $I_{\mathbb C}$ also leaves the vector spaces $F^{p,k}$ of the associated Hodge filtration invariant, and we again denote by $(F^{p,k})^+$ and $(F^{p,k})^-$ the $+1$ and $-1$ eigenspaces of the induced involution on $F^{p,k}$. The rational Hodge structures $(V^+,(F^{\ast,k})^+)$ and $(V^-,(F^{\ast,k})^-)$ are rational Hodge substructures of $(V^k,F^{\ast,k})$. Furthermore, the rational Hodge structure $(V^k,F^{\ast,k})$ is a direct sum of $(V^+,(F^{\ast,k})^+)$ and $(V^-,(F^{\ast,k})^-)$. We drop the subscript ${\mathbb C}$ of $I_{\mathbb C}$ in what follows.

\noindent Let $I$ be an involution of a compact K\"ahler manifold $X$, that is a holomorphic map $I:X\rightarrow X$ with $I^2={\rm Id}_X$. Let $n=\dim_{\mathbb C}(X)$. The pull-back $I^\ast$ by $I$ of complex differential forms sends a closed $(p,q)$-form to a closed $(p,q)$ form,  and therefore defines a map on $H^k_\Omega(X,{\mathbb C})$ preserving $H^{p,q}$. Moreover, by the functoriality of the de-Rham isomorphism, it corresponds to the pull-back induced by $I$ on singular cohomology, which maps $H^k(X,{\mathbb Z})$ to itself. The induced ${\mathbb Q}$-linear involution preserves $V^k=H^k(X,{\mathbb Q})$,  $0\le k\le n$, and we denote it also by $I$.  Indeed, if $\sigma=\sum_ia_i\sigma_i$ is a singular $k$-chain with coefficients $a_i\in{\mathbb Z}$, then so is $I_\ast(\sigma)=\sum_ia_i(I\circ\sigma_i)$, and $I_\ast$ commutes with the boundary operator. This induces a map $I$ on $H^k(X,{\mathbb Z})$ which extends by ${\mathbb Q}$-linearity to the map $I$ on $H^k(X,{\mathbb Q})$. Let $\omega\in H^k(X,{\mathbb R})$ and let $\varphi$ be a closed $k$-form on $X$ whose class $[\varphi]$ represents $\omega$ via the de-Rham isomorphism. Then, the functoriality of the de-Rham isomorphism amounts to the fact that $I[\varphi]=[I^\ast\varphi]$ (see \cite{GrH}, p.45). As we have seen, the $\pm1$ eigenspaces $(V^k)^\pm$ of $V^k$, with respect to $I$, carry rational Hodge substructures of $H^k(X,{\mathbb Q}_X)$, which we denote $H^k(X,{\mathbb Q}_X)^\pm$.

\noindent Returning to the Borcea-Voisin construction, for all $k\le n$, denote also by ${\widehat{I}}_{1,2}$  the involution on $H^k({\widehat{Y}},{\mathbb Q})$ induced by the involution ${\widehat{I}}_{1,2}$ on ${\widehat{Y}}$. As $B={\widehat{Y}}/{\widehat{I}}_{1,2}$, we can identify $H^k(B,{\mathbb Q}_B)$ with $H^k({\widehat{Y}},{\mathbb Q}_{\widehat{Y}})^{{\widehat{I}}_{1,2}}$ (see \cite{Gro}). Let $F^{\ast,k}(B)$ be the Hodge filtration of $H^k(B,{\mathbb Q}_B)_{\mathbb C}$.

\noindent By Lemma 3, \S\ref{s:lemmas}, we have
$$
F^{p,k}({\widehat{Y}})=F^{p,k}(A_1\times A_2)\oplus F^{p-1,k-2}(R_1\times R_2),
$$
where the second summand on the right hand side is empty for $k=0,1$, and for $p=0$. Considering the action of ${\widehat{I}}_{1,2}$ on the left hand side and $I_{1,2}=I_1\times I_2$ on the right hand side, we have
$$
F^{p,k}(B)=F^{p,k}({\widehat{Y}})^{{\widehat{I}}_{1,2}}=F^{p,k}(A_1\times A_2)^{I_{1,2}}\oplus F^{p-1,k-2}(R_1\times R_2).
$$

\noindent Let $V^k=H^k(A_1\times A_2,{\mathbb Q})^{I_{1,2}}$. For $i=1,2$, let $V^{r,\pm}_i=H^r(A_i,{\mathbb Q})^{\pm1}$, be the $\pm1$-eigenspaces with respect to the involution $I_i$. As explained above, these ${\mathbb Q}$-vector spaces carry the rational Hodge structures given by $(V^k,{\widetilde{\varphi}}^{(k)})=H^k(A_1\times A_2,{\mathbb Q}_{A_1\times A_2})^{I_{1,2}}$ and $(V^{r,\pm}_i,{\widetilde{\varphi}}^{(r,\pm)}_i)=H^r(A_i,{\mathbb Q}_{A_i})^\pm$, $i=1,2$. By the assumptions of the proposition $F^{\ast,k}(B)$ is defined over ${\overline{\mathbb Q}}$, so, by Lemma 4, \S\ref{s:lemmas}, the filtration $F^{\ast,k}(A_1\times A_2)^{I_{1,2}}$ is also defined over ${\overline{\mathbb Q}}$, in that it is obtained by extension of scalars to ${\mathbb C}$ from a ${\overline{\mathbb Q}}$-filtration of $V^k_{\overline{\mathbb Q}}$. We have by Lemma 2, \S\ref{s:lemmas}, that
$$
V^k=\;\oplus_{r+s=k}\left(V^{r,+}_1\otimes V^{s,+}_2\right)\oplus
\left(V^{r,-}_1\otimes V^{s,-}_2\right),
$$
where the rational Hodge structures on the vector spaces are related by
$$
{\widetilde{\varphi}}^{(k)}=\;\oplus_{r+s=k}\left({\widetilde{\varphi}}^{(r,+)}_1\otimes{\widetilde{\varphi}}^{(s,+)}_2\right)\oplus
\left({\widetilde{\varphi}}^{(r,-)}_1\otimes{\widetilde{\varphi}}^{(s,-)}_2\right).
$$
As the Hodge filtration associated to $(V^k,{\widetilde{\varphi}}^{(k)})$ is defined over ${\overline{\mathbb Q}}$, we can again use Lemma 4, \S\ref{s:lemmas}, to deduce that the Hodge filtrations associated to the direct summands on the right hand side of the above equations are all defined over ${\overline{\mathbb Q}}$. Therefore, for all $r,s$ with $r+s=k$, the Hodge filtration associated to
$$
(V^{r,s,\pm},{\widetilde{\varphi}}^{(r,s,\pm)})=\left(V^{r,\pm}_1\otimes V^{s,\pm}_2,{\widetilde{\varphi}}^{(r,\pm)}_1\otimes{\widetilde{\varphi}}^{(s,\pm)}_2\right)
$$
is defined over ${\overline{\mathbb Q}}$. By Lemma 2, \S\ref{s:lemmas}, this filtration is given by
$$
\left(F^{p,r,s}\right)^\pm=\sum_{i+j=p}\left(F^{i,r}(A_1)^\pm\otimes_{\mathbb C}F^{j,s}(A_2)^\pm\right),
$$
where the sum over the $i,j$ with $i+j=p$ is not direct, and $p=0,\ldots,k$. To see that the above sum is not direct, notice, for example, that we can have a $v_i\in F^{i,r}(A_1)^\pm$, $i\le p-1$, which also satisfies $v_i\in F^{i+1,r}(A_1)^\pm$. On the other hand, any $w_j\in F^{j,s}(A_2)^{\pm}$, $j\ge1$, is also in $F^{j-1,s}(A_2)^\pm$. Therefore, in that case, the vector $v_i\otimes_{\mathbb C}w_j$ is both in $F^{i,r}(A_1)^\pm\otimes_{\mathbb C}F^{j,s}(A_2)^\pm$ and in $F^{i+1,r}(A_1)^\pm\otimes_{\mathbb C}F^{j-1,s}(A_2)^\pm$.

\noindent We first consider the case $k=d$, $r=d_1$, and $s=d_2$, where, as before $d_1=\dim(A_1)$, $d_2=\dim(A_2)$ and $d=d_1+d_2=\dim(A_1\times A_2)=\dim(B)$. Recall that, for $i=1,2$, the variety $A_i$ is assumed to be Calabi-Yau and that, by \cite{Roh}, Lemma 7.2.4, we have $F^{d_i,d_i}(A_i)^+=\{0\}$ and $\dim(F^{d_i,d_i}(A_i)^-)=1$, the vector space $F^{d_i,d_i}(A_i)^-=H^{d_i,0}(A_i)$ being generated by the inverse image under the de-Rham isomorphism of a nowhere vanishing holomorphic $d_i$-form. Therefore,
$$
\left(F^{d,d_1,d_2}\right)^-= F^{d_1,d_1}(A_1)^-\otimes_{\mathbb C}F^{d_2,d_2}(A_2)^-,
$$
$$
=H^{d_1,0}(A_1)\otimes_{\mathbb C}H^{d_2,0}(A_2),
$$
which is a 1-dimensional ${\mathbb C}$-vector space, whereas $\left(F^{d,d_1,d_2}\right)^+=\{0\}$. Now, since $(F^{\ast,d_1,d_2})^-$ is defined over ${\overline{\mathbb Q}}$, there is an $\Omega\in V^{d_1,-}_{1,{\overline{\mathbb Q}}}\otimes_{\overline{\mathbb Q}} V^{d_2,-}_{2,{\overline{\mathbb Q}}}$ such that $(F^{d,d_1,d_2})^-={\mathbb C}(\Omega\otimes_{\overline{\mathbb Q}}1_{\mathbb C})$. In addition, we have $\Omega\otimes_{\overline{\mathbb Q}}1_{\mathbb C}=\Omega_1'\otimes_{\mathbb C}\Omega_2'$, for some $\Omega_i'\in F^{d_i,d_i}(A_i)^-$, $i=1,2$. By Lemma 5, \S\ref{s:lemmas}, we can choose $c_1, c_2$ such that $c_1\Omega_1'\in H^{d_1}(A_1,{\overline{\mathbb Q}})^-$ and $c_2\Omega_2'\in H^{d_2}(A_2,{\overline{\mathbb Q}})^-$. Therefore, for $i=1,2$, the 1-dimensional ${\mathbb C}$-vector space $F^{d_i,d_i}=H^{d_i,0}(A_i)$ has a generator $\Omega_i$ in $H^{d_i}(A_i,{\overline{\mathbb Q}})$. This gives the same result as that of Example 2, \S\ref{s:examples}, namely, if the Calabi-Yau varieties with involution $A_1$ and $A_2$ satisfy Expectation 4, \S\ref{s:examples}, then the Calabi-Yau variety with involution $B$ satisfies Expectation 4. This is enough for most the examples of this paper, namely those that exploit the Analytic Subgroup Theorem.

\noindent Recall from \S\ref{s:VHS} that, for $i=1,2$, and $0\le k_i\le d_i$, we have a non-degenerate bilinear form $Q_i$ on $H^{k_i}(A_i,{\mathbb Q})_{\rm prim}$ which extends to the full cohomology space $H^{k_i}(A_i,{\mathbb Q})$  and yields a non-degenerate bilinear form there with $Q_i(u,v)=(-1)^{k_i}Q_i(v,u)$. By the Lefschetz primitive decomposition, the rational Hodge structure $H^{k_i}(A_i,{\mathbb Q}_X)$ is a direct sum of rational Hodge substructures that are mutually orthogonal with respect to $Q_i$. Moreover, the Hodge decomposition of each substructure is orthogonal with respect to $h_i(u,v)=(\sqrt{-1})^{k_i}Q_i(u,{\overline{v}})$ and, on each summand of this Hodge decomposition, $h_i$ is definite.

\noindent For $i=1,2$, let $Q_i^{\pm}$ be the restriction of $Q_i$ to $V_i^{r_i,\pm}$. Let $K={\overline{\mathbb Q}}$ or ${\mathbb C}$ and identify ${\overline{\mathbb Q}}$ with a subfield of ${\mathbb C}$. We have a $K$-valued hermitian form $h_{i,K}^{\pm}(u,v)=(\sqrt{-1})^{r_i}Q_i^{\pm}(u,{\overline{v}})$ on $V^{r_i,\pm}_K$. By the discussion of the preceding paragraph,  the rational Hodge structure $V_i^{r_i,\pm}$ is a direct sum of rational Hodge substructures that are mutually orthogonal with respect to $Q_i^{\pm}$. In addition, the Hodge decomposition of each substructure is orthogonal with respect to $h_{i,{\mathbb C}}^{\pm}$, and $h_{i,{\mathbb C}}^{\pm}$ is definite on each summand of this decomposition. By the discussion of \S\ref{s:VHS}, the vector space $V^{r,s,\pm}_K$ associated to the rational Hodge structure $(V^{r,s,\pm},{\widetilde{\varphi}}^{(r,s,\pm)})$ carries a $K$-valued hermitian form $h_K^{\pm}$ such that
$$
h_{\mathbb C}^{\pm}(H^{i,r-i}(A_1)^{\pm}\otimes H^{j,s-j}(A_2)^{\pm},\; H^{i',r-i'}(A_1)^{\pm}\otimes H^{j',s-j'}(A_2)^{\pm})=0
$$
when $(i,j)\not=(i',j')$.

\noindent Recall that, using the fact that $(F^{d,d_1,d_2})^-$ is defined over ${\overline{\mathbb Q}}$, we have shown, for $i=1,2$, that $H^{d_i,0}(A_i)={\mathbb C}(\Omega_i\otimes_{\overline{\mathbb Q}}1_{\mathbb C})$ where $\Omega_i\in H^{d_i}(A_i,{\overline{\mathbb Q}})^-$. We now look at $(F^{d-1,d_1,d_2})^-$. We have
$$
(F^{d-1,d_1,d_2})^-=(F^{d,d_1,d_2})^-\oplus \left(\oplus_{i+j=1}H^{d_1-i, i}(A_1)^-\otimes_{\mathbb C}H^{d_2-j,j}(A_2)^-\right).
$$
$$
=(H^{d_1,0}(A_1)\otimes_{\mathbb C}H^{d_2,0}(A_2))
$$
$$
\oplus (H^{d_1-1,1}(A_1)^-\otimes_{\mathbb C}H^{d_2,0}(A_2))\oplus (H^{d_1,0}(A_1)\otimes_{\mathbb C}H^{d_2-1,1}(A_2)^-).
$$
By our discussion so far, the vector space $(F^{d-1,d_1,d_2})^-$ is of the form $V_{\overline{\mathbb Q}}\otimes_{\overline{\mathbb Q}}{\mathbb C}$ for $V_{\overline{\mathbb Q}}$ a subspace of $V^{d_1,d_2,-}_{\overline{\mathbb Q}}$, and $H^{d_1,0}(A_1)\otimes_{\mathbb C}H^{d_2,0}(A_2)$ is a vector subspace of $V^{d_1,d_2,-}_{\mathbb C}$ of the form $U_{1,{\mathbb C}}=U_{1,{\overline{\mathbb Q}}}\otimes_{\overline{\mathbb Q}}{\mathbb C}$, where $U_{1,{\overline{\mathbb Q}}}$ is a vector subspace of $V_{\overline{\mathbb Q}}$ of dimension 1 on which $h_{\overline{\mathbb Q}}^-$ is definite. Moreover, the direct sum
$$
U_2=(H^{d_1-1,1}(A_1)^-\otimes_{\mathbb C}H^{d_2,0}(A_2))\oplus (H^{d_1,0}(A_1)\otimes_{\mathbb C}H^{d_2-1,1}(A_2)^-)
$$
is orthogonal to
$$
H^{d_1,0}(A_1)\otimes_{\mathbb C}H^{d_2,0}(A_2)
$$
with respect to $h_{\mathbb C}^-$.
Therefore, by \S\ref{s:lemmas}, Lemma 9, it follows that $U_2$ has a basis in $U_2\cap \{v\otimes_{\overline{\mathbb Q}}1_{\mathbb C}: v\in V_{\overline{\mathbb Q}}\}$. If $H^{d_1-1,1}(A_1)^-$ and $H^{d_2-1,1}(A_2)^-$ are both trivial, then $U_2=\{0\}$. If only one, say $H^{d_1-1,1}(A_1)^-$, is trivial, we can use \S\ref{s:lemmas}, Lemma 6, to deduce that the other, now $H^{d_2-1,1}(A_2)^-$, has a basis in
$$
H^{d_2-1,1}(A_2)^-\cap \{v\otimes_{\overline{\mathbb Q}}1_{\mathbb C}: v\in H^{d_2}(A_2,{\overline{\mathbb Q}})^-\}.$$
If both are non-trivial, we may use \S\ref{s:lemmas}, Lemma 7, to deduce that
$$
(F^{d_i-1,d_i}(A_i))^-=H^{d_i,0}(A_i)\oplus H^{d_i-1,1}(A_i)^-
$$
has a basis in
$$
(H^{d_i,0}(A_i)\oplus H^{d_i-1,1}(A_i)^-)\cap \{v\otimes_{\overline{\mathbb Q}}1_{\mathbb C}: v\in H^{d_i}(A_i,{\overline{\mathbb Q}})^-\},\quad i=1,2.
$$
Let $V_{i,{\overline{\mathbb Q}}}$ be a ${\overline{\mathbb Q}}$-vector subspace of $V^{d_i,-}_{\overline{\mathbb Q}}$ such that $(F^{d_i-1,d_i}(A_i))^-=V_{i,{\overline{\mathbb Q}}}\otimes_{\overline{\mathbb Q}}{\mathbb C}$. Let $U_{i,{\overline{\mathbb Q}}}={\overline{\mathbb Q}}\Omega_i$, so that $U_{i,{\mathbb C}}=H^{d_i,0}(A_i)$. Then, applying Lemma 9 by using the fact that $H^{d_i,0}(A_i)$ and $H^{d_i-1,1}(A_i)^-$ are orthogonal with respect to $h_{i,{\mathbb C}}^-$, we deduce that $U_{i,2}=H^{d_i-1,1}(A_i)^-$ has a basis in $U_{i,2}\cap \{v\otimes_{\overline{\mathbb Q}}1_{\mathbb C}: v\in V_{i,{\overline{\mathbb Q}}}\}$.

\noindent We now proceed inductively, in a manner similar to the above passage from $(F^{d,d_1,d_2})^-$ to $(F^{d-1,d_1,d_2})^-$. Suppose, for $j=0,\ldots,k$, where $k<\min(d_1,d_2)$, and $i=1,2$, that  $H^{d_i-j,j}(A_i)^-$, has a basis in
$$
H^{d_i-j,j}(A_i)^-\cap \{v\otimes_{\overline{\mathbb Q}}1_{\mathbb C}: v\in H^{d_i}(A_i,{\overline{\mathbb Q}})^-\}.
$$
By the discussion at the beginning of the proof, the vector space $(F^{d-k-1,d_1,d_2})^-$ is of the form $V_{\overline{\mathbb Q}}\otimes_{\overline{\mathbb Q}}{\mathbb C}$ for $V_{\overline{\mathbb Q}}$ a subspace of $V^{d_1,d_2,-}_{\overline{\mathbb Q}}$ and the vector space
$$
(F^{d-k,d_1,d_2})^-=\oplus_{c'\le k}\oplus_{a+b=c'}(H^{d_1-a,a}(A_1)^-\otimes_{\mathbb C}H^{d_2-b,b}(A_2)^-)
$$
is of the form $U_{1,{\mathbb C}}'=U_{1,{\overline{\mathbb Q}}}'\otimes_{\overline{\mathbb Q}}{\mathbb C}$ where $U_{1,{\overline{\mathbb Q}}}'$ is a vector subspace of $V_{\overline{\mathbb Q}}$. By the induction hypothesis and the discussion of \S\ref{s:VHS}, the vector space $U_{1,{\overline{\mathbb Q}}}'$ is an orthogonal direct sum, with respect to $h_{\overline{\mathbb Q}}^-$, of ${\overline{\mathbb Q}}$-vector subspaces on which $h_{\overline{\mathbb Q}}^-$ is definite. We have
$$
(F^{d-k-1,d_1,d_2})^-=(F^{d-k,d_1,d_2})^-\oplus\left(\oplus_{i+j=k+1}H^{d_1-i,i}(A_1)^-\otimes_{\mathbb C}H^{d_2-j,j}(A_2)^-\right).
$$
$$
=U_{1,{\mathbb C}}'\oplus\left(\oplus_{i+j=k+1, 1\le i\le k}H^{d_1-i,i}(A_1)^-\otimes_{\mathbb C}H^{d_2-j,j}(A_2)^-\right)
$$
$$
\oplus\left(H^{d_1,0}(A_1)\otimes_{\mathbb C}H^{d_2-k-1, k+1}(A_2)^-\right)\oplus \left(H^{d_1-k-1,k+1}(A_1)^-\otimes_{\mathbb C}H^{d_2,0}(A_2)\right).
$$
Again, by the induction hypothesis and the discussion of \S\ref{s:VHS}, the vector space
$$
U_{1,{\mathbb C}}=U_{1,{\mathbb C}}'\oplus\left(\oplus_{i+j=k+1, 1\le i\le k}H^{d_1-i,i}(A_1)^-\otimes_{\mathbb C}H^{d_2-j,j}(A_2)^-\right)
$$
is of the form $U_{1,{\overline{\mathbb Q}}}\otimes_{\overline{\mathbb Q}}{\mathbb C}$ for a subspace $U_{1,{\overline{\mathbb Q}}}$ of $V_{\overline{\mathbb Q}}$ and
$U_{1,{\overline{\mathbb Q}}}$ is an orthogonal direct sum, with respect to $h_{\overline{\mathbb Q}}^-$, of ${\overline{\mathbb Q}}$-vector subspaces on which $h_{\overline{\mathbb Q}}^-$ is definite. Therefore, by \S\ref{s:lemmas}, Lemma 9,
it follows that
$$
U_2=\left(H^{d_1,0}(A_1)\otimes_{\mathbb C}H^{d_2-k-1, k+1}(A_2)^-\right)\oplus \left(H^{d_1-k-1,k+1}(A_1)^-\otimes_{\mathbb C}H^{d_2,0}(A_2)\right)
$$
has a basis in $U_2\cap \{v\otimes_{\overline{\mathbb Q}}1_{\mathbb C}: v\in V_{\overline{\mathbb Q}}\}$. If the vector spaces $H^{d_1-k-1,k+1}(A_1)^-$ and $H^{d_2-k-1,k+1}(A_2)^-$ are both trivial, then $U_2=\{0\}$. If only one, say $H^{d_1-k-1,k+1}(A_1)^-$, is trivial, we can use \S\ref{s:lemmas}, Lemma 6, to deduce that the other, now $H^{d_2-k-1,k+1}(A_2)^-$, has a basis in
$$
H^{d_2-k-1,k+1}(A_2)^-\cap \{v\otimes_{\overline{\mathbb Q}}1_{\mathbb C}: v\in H^{d_2}(A_2,{\overline{\mathbb Q}})^-\}.$$
If both are non-trivial, we may use \S\ref{s:lemmas}, Lemma 7, to deduce that
$$
H^{d_i,0}(A_i)\oplus H^{d_i-k-1,k+1}(A_i)^-
$$
has a basis in
$$
(H^{d_i,0}(A_i)\oplus H^{d_i-k-1,k+1}(A_i)^-)\cap \{v\otimes_{\overline{\mathbb Q}}1_{\mathbb C}: v\in H^{d_i}(A_i,{\overline{\mathbb Q}})^-\},\quad i=1,2.
$$
Let $V_{i,{\overline{\mathbb Q}}}$ be a ${\overline{\mathbb Q}}$-vector subspace of $V^{d_i,-}_{\overline{\mathbb Q}}$ such that
$$
H^{d_i,0}(A_i)\oplus H^{d_i-k-1,k+1}(A_i)^-=V_{i,{\overline{\mathbb Q}}}\otimes_{\overline{\mathbb Q}}{\mathbb C}.
$$
Let $U_{1,{\overline{\mathbb Q}}}={\overline{\mathbb Q}}\Omega_i$, so that $U_{1,{\mathbb C}}=H^{d_i,0}(A_i)$. Applying \S\ref{s:lemmas}, Lemma 9, by using the fact that $H^{d_i,0}(A_i)^-$ and $H^{d_i-k-1,k+1}(A_i)^-$ are orthogonal with respect to $h_{i,{\mathbb C}}^-$, we deduce that $U_2=H^{d_i-k-1,k+1}(A_i)^-$ has a basis in
$$
U_2\cap \{v\otimes_{\overline{\mathbb Q}}1_{\mathbb C}: v\in V_{\overline{\mathbb Q}}\}.
$$
If $\min(d_1,d_2)\le k< \max(d_1,d_2)$, then we may argue in the following manner. If $d_1=d_2$, then we are already done. If $d_2=\max(d_1,d_2)>\min(d_1,d_2)$, say, then we may use \S\ref{s:lemmas}, Lemma 6, applied to
$$
W_1=\left(H^{d_1,0}(A_1)\otimes_{\mathbb C}H^{d_2-k-1, k+1}(A_2)^-\right)\oplus \left(H^{d_1-k-1,k+1}(A_1)^-\otimes_{\mathbb C}H^{d_2,0}(A_2)\right)
$$
$$
=H^{d_1,0}(A_1)\otimes_{\mathbb C}H^{d_2-k-1,k+1}(A_2)^-,
$$
to deduce that  $U_2=H^{d_2-k-1,k+1}(A_2)^-$ has a basis in $U_2\cap \{v\otimes_{\overline{\mathbb Q}}1_{\mathbb C}: v\in V_{\overline{\mathbb Q}}\}$. Once we get to $k=\max(d_1,d_2)$, we stop. Notice that some of these steps may be redundant, since if $H^{p,q}(A_i)^-$ has a basis in  $\{v\otimes_{\overline{\mathbb Q}}1_{\mathbb C}: v\in V^{d_i,-}_{\overline{\mathbb Q}}\}$, then so does $H^{q,p}(A_i)={\overline{H}}^{p,q}(A_i)$.

\noindent To summarize, for $i=1,2$, we have shown that the Hodge filtration $F^{\ast,d_i}(A_i)^-$ of $H^{d_i}(A_i,{\mathbb Q}_{A_i})_{\mathbb C}^-$, and all the summands $H^{p,q}(A_i)^-$,  $p+q=d_i$, have bases in $H^{d_i}(A_i,{\overline{\mathbb Q}})$. Moreover, to do this, we have only used the assumption that the Hodge filtration $(F^{\ast, d_1,d_2})^-$ of the tensor product  of the rational Hodge structures $H^{d_1}(A_1,{\mathbb Q}_{A_1})^-\otimes H^{d_2}(A_2,{\mathbb Q}_{A_2})^-$ is defined over ${\overline{\mathbb Q}}$.

\noindent As $H^{d_i,0}(A_i)^+=\{0\}$, we cannot use these same arguments to show that, for $i=1,2$, the Hodge filtration $F^{\ast, d_i}(A_i)^+$ of $H^{d_i}(A_i,{\mathbb Q}_{A_i})_{\mathbb C}^+$ is defined over ${\overline{\mathbb Q}}$. Instead, we use the fact that, for any smooth connected projective variety $X$, we have $H^0(X,{\mathbb C})=H^{0,0}(X)=H^0(X,{\overline{\mathbb Q}})\otimes_{\overline{\mathbb Q}}{\mathbb C}$, and this vector space is of dimension 1. Now, by our assumptions, the Hodge filtration $(F^{\ast, d_1,0})^+$ of the Hodge structure $(V^{d_1,0,+},{\widetilde{\varphi}}^{(d_1,0,+)})$ of weight $d_1$ is defined over ${\overline{\mathbb Q}}$. Recall that
$$
V^{d_1,0,+}=H^{d_1}(A_1,{\mathbb Q})^+\otimes_{\mathbb Q} H^0(A_2,{\mathbb Q})^+=H^{d_1}(A_1,{\mathbb Q})^+\otimes_{\mathbb Q} H^0(A_2,{\mathbb Q}),
$$
since $H^0(A_2,{\mathbb Q})^+=H^0(A_2,{\mathbb Q})$. Therefore,
$$
(F^{d_1,d_1,0})^+= H^{d_1,0}(A_1)^+\otimes_{\mathbb C}H^{0,0}(A_2)=\{0\}.
$$
By \S\ref{s:lemmas}, Lemma 10, and the discussion of \S\ref{s:VHS}, all the summands $(V^{d_1,0,+}_{\mathbb C})^{p,q}$, $p+q=d_1$, in the Hodge decomposition of $(V^{d_1,0,+},{\widetilde{\varphi}}^{(d_1,0,+)})$
have a basis in $(V^{d_1,0,+}_{\mathbb C})^{p,q}\cap \{v\otimes_{\overline{\mathbb Q}}1_{\mathbb C}: v\in V^{d_1,0,+}_{\overline{\mathbb Q}}\}$. We have, for $i=1,\ldots, d_1$, that
$$
(V^{d_1,0,+}_{\mathbb C})^{d_1-i,i}= H^{d_1-i,i}(A_1)^+\otimes_{\mathbb C}H^{0,0}(A_2),
$$
so by \S\ref{s:lemmas}, Lemma 6, with $V_{2,{\overline{\mathbb Q}}}=H^0(A_2,{\overline{\mathbb Q}})$ and $U_1=H^{d_1-i,i}(A_1)^+$, we deduce that $H^{d_1-i,i}(A_1)^+$ has a basis in $H^{d_1}(A_1,{\overline{\mathbb Q}})^+$.
An analogous argument shows that $H^{d_2-j,j}(A_2)^+$ has a basis in $H^{d_2}(A_2,{\overline{\mathbb Q}})^+$, $j=1,\ldots,d_2$.

\noindent To summarize, with the hypotheses of the lemma, in particular using the fact that $(F^{\ast, d_1, d_2})^-$, $(F^{\ast, d_1,0})^+$ and $(F^{\ast, 0,d_2})^+$ are defined over ${\overline{\mathbb Q}}$, we have shown that, for $i=1,2$, the Hodge filtration of $H^{d_i}(A_i,{\mathbb Q}_{A_i})_{\mathbb C}$ is defined over ${\overline{\mathbb Q}}$.

\noindent More generally, we have that $(F^{k-p,p,0})^+$ is defined over ${\overline{\mathbb Q}}$ for $k=0,\ldots,d_1$, $p=0,\ldots k$, and
$$
(F^{k-p,p,0})^+=H^{k-p,p}(A_1)^+\otimes_{\mathbb C}H^{0,0}(A_2),
$$
so that by \S\ref{s:lemmas}, Lemma 6, it follows that $H^{k-p,p}(A_1)^+$ has a basis in $H^k(A_1,{\overline{\mathbb Q}})^+$. An analogous argument shows that $H^{\ell-q,q}(A_2)^+$ has a basis in $H^\ell(A_2,{\overline{\mathbb Q}})^+$, $\ell=0,\ldots,d_2$, $q=0,\ldots,\ell$.

\noindent We conclude by showing inductively that all the $H^{k-p,p}(A_1)^-$, $k=0,\ldots,d_1$, $p=0,\ldots k$, have bases in $H^k(A_1,{\overline{\mathbb Q}})^-$ and all the $H^{\ell-q,q}(A_2)^-$, $\ell=0,\ldots,d_2$, $q=0,\ldots,\ell$, have bases in  $H^\ell(A_2,{\overline{\mathbb Q}})^-$. We have already shown this for the $H^{d_1-i,i}(A_1)^-$, $i=0,\ldots, d_1$ and the $H^{d_2-j,j}(A_2)^-$, $j=0,\ldots, d_2$. Therefore, suppose that $\ell<d_2$. We know that $(F^{\ast, d_1,\ell})^-$ is defined over ${\overline{\mathbb Q}}$. By \S\ref{s:lemmas}, Lemma 10, and the discussion of \S\ref{s:VHS}, all the summands $(V^{d_1,\ell,-}_{\mathbb C})^{p,q}$, $p+q=d_1+\ell$, in the Hodge decomposition of $(V^{d_1,\ell,-},{\widetilde{\varphi}}^{(d_1,\ell,-)})$ have a basis in $(V^{d_1,\ell,-}_{\mathbb C})^{p,q}\cap \{v\otimes_{\overline{\mathbb Q}}1_{\mathbb C}: v\in V^{d_1,\ell,-}_{\overline{\mathbb Q}}\}$. Let $\ell$ be a fixed integer between $0$ and $d_2$. We have
$$
(V^{d_1,\ell,-}_{\mathbb C})^{(d_1+\ell,0)}=H^{d_1,0}(A_1)\otimes_{\mathbb C}H^{\ell, 0}(A_2)^-=\{0\},
$$
as $A_2$ is Calabi-Yau. We have
$$
(V^{d_1,\ell,-}_{\mathbb C})^{(d_1+\ell-1,1)}=(H^{d_1,0}(A_1)\otimes_{\mathbb C}H^{\ell-1,1}(A_2)^-)\oplus(H^{d_1-1,1}(A_1)^-\otimes_{\mathbb C}H^{\ell,0}(A_2)^-),
$$
$$
=H^{d_1,0}(A_1)\otimes_{\mathbb C}H^{\ell-1,1}(A_2)^-,
$$
again by the Calabi-Yau assumption on $A_2$. From \S\ref{s:lemmas}, Lemma 6, it follows that $H^{\ell-1,1}(A_2)^-$, which may be trivial, has a basis in $H^\ell(A_2,{\overline{\mathbb Q}})^-$. Assume that, for some integer $j$ between $0$ and $\ell-1$, we have shown that $H^{\ell-i,i}(A_2)^-$ has a basis in $H^{\ell}(A_2,{\overline{\mathbb Q}})^-$, for $i=0,\ldots, j$. We have,
$$
(V^{d_1,\ell,-}_{\mathbb C})^{(d_1+\ell-j-1,j+1)}
$$
$$
=(H^{d_1,0}(A_1)\otimes_{\mathbb C}H^{\ell-j-1,j+1}(A_2)^-)
$$
$$
\oplus(\oplus_{i=0}^j(H^{d_1-j-1+i,j+1-i}(A_1)^-\otimes_{\mathbb C}H^{\ell-i,i}(A_2)^-),
$$
where, in the last sum, we set $H^{d_1-j-1+i,j+1-i}(A_1)^-=\{0\}$, if $j-i+1>d_1$. By the induction hypothesis, and what we have shown already, we can again use \S\ref{s:lemmas}, Lemma 9, this time to deduce that $H^{d_1,0}(A_1)\otimes_{\mathbb C}H^{\ell-j-1,j+1}(A_2)^-$ has a basis in $\{v\otimes_{\overline{\mathbb Q}}1_{\mathbb C}: v\in V^{d_1,\ell,-}_{\overline{\mathbb Q}}\}$. Therefore, by \S\ref{s:lemmas}, Lemma 6, we deduce that $H^{\ell-j,j}(A_2)^-$ has a basis in $H^{\ell}(A_2,{\overline{\mathbb Q}})^-$. Therefore, by induction, this last statement is true for all $j=0,\ldots,\ell$. A similar argument for $A_1$ shows that $H^{k-i,i}(A_1)^-$ has a basis in $H^k(A_1,{\overline{\mathbb Q}})^-$, $i=0,\ldots,k$, $k=0,\ldots,d_1$.

\noindent Putting all these steps together, we have finally shown that, for $i=1,2$, the Hodge filtration of $H^{k_i}(A_i,{\mathbb Q}_{A_i})_{\mathbb C}=H^{k_i}(A_i,{\mathbb Q}_{A_i})_{\mathbb C}^+\oplus H^{k_i}(A_i,{\mathbb Q}_{A_i})_{\mathbb C}^-$ is defined over ${\overline{\mathbb Q}}$, for all $k_i\le\dim(A_i)$. This completes the first statement of the Proposition.

\noindent We now look at the Hodge filtration of the $H^{\ell_j}(R_j,{\mathbb Q}_{R_j})_{\mathbb C}$, $j=1,2$. Here, the argument is easier, as we do not have to worry about $\pm1$-eigenspaces. We let $(W_j^{r_j},{\widetilde{\varphi}}_j^{r_j})$ be the Hodge structure $H^{r_j}(R_j,{\mathbb Q}_{R_j})$, $j=1,2$, and $(W^k,{\widetilde{\varphi}}^k)$ the Hodge structure $H^k(R_1\times R_2,{\mathbb Q}_{R_1\times R_2})$. By \S\ref{s:lemmas}, we have
$$
W^k=\oplus_{r_1+r_2=k}  (W_1^{r_1}\otimes W_2^{r_2})
$$
where the rational Hodge structures on the vector spaces are related by
$$
{\widetilde{\varphi}}^k=\oplus_{r_1+r_2=k}  ({\widetilde{\varphi}}_1^{r_1}\otimes{\widetilde{\varphi}}_2^{r_2}).
$$
By the assumptions of the proposition, the Hodge filtration associated to $(W^k,{\widetilde{\varphi}}^{k})$ is defined over ${\overline{\mathbb Q}}$, so that by \S\ref{s:lemmas}, Lemma 4, the Hodge filtrations associated to each $(W_1^{r_1}\otimes W_2^{r_2}, {\widetilde{\varphi}}_1^{r_1}\otimes{\widetilde{\varphi}}_2^{r_2})$ are defined over ${\overline{\mathbb Q}}$. Again using \S\ref{s:lemmas}, Lemma 10, each summand in the Hodge decomposition of $(W_1^{r_1}\otimes_{\mathbb Q} W_2^{r_2})_{\mathbb C}$ has a basis in $(W_1^{r_1}\otimes_{\mathbb Q}W_2^{r_2})_{\overline{\mathbb Q}}$. This is in particular true for the case $r_2=0$, and we have,
$$
(W_1^{r_1}\otimes_{\mathbb Q} W_2^{0})_{\mathbb C}^{p,q}=\oplus_{i+i'=p, j+j'=q}(W_{1,{\mathbb C}}^{r_1})^{(i,j)}\otimes_{\mathbb C}(W_{2,{\mathbb C}}^{0})^{(i',j')}
$$
$$
=
(W_{1,{\mathbb C}}^{r_1})^{(p,q)}\otimes (W_{2,{\mathbb C}}^{0})^{(0,0)}.
$$
By \S\ref{s:lemmas}, it follows that $(W_{1,{\mathbb C}}^{r_1})^{(p,q)}$, has a basis in $W_{1,{\overline{\mathbb Q}}}^{p+q}$. We can apply a similar argument to $(W_{2,{\mathbb C}}^{r_2})^{(p,q)}$,
to show it has a basis in $W_{2,{\overline{\mathbb Q}}}^{p+q}$. This completes the proof of the proposition.

\noindent We end by remarking that the arguments of the proposition can no doubt be formulated in terms of period matrices and their Kronecker products. It is also likely that a much more efficient argument can be found.

\medskip

\noindent Consider two smooth projective families ${\mathcal A}_1\rightarrow S_1$, and ${\mathcal A}_2\rightarrow S_2$ of Calabi-Yau varieties defined over ${\overline{\mathbb Q}}$. Let $A_1(s_1)$ be the fiber of ${\mathcal A}_1$ at $s_1\in S_1$ and $A_2(s_2)$ the fiber of ${\mathcal A}_2$ at $s_2\in S_2$.  Suppose, for $i=1,2$, that $\dim(A_i(s_i))=d_i$ (independent of $s_i\in S_i$). Recall from \S\ref{s:VHS} that, by definition, the variety $A_i(s_i)$ is defined over ${\overline{\mathbb Q}}$ when $s_i\in S_i({\overline{\mathbb Q}})$. Suppose, for $i=1,2$, there is an involution ${\mathcal I}_i$ of ${\mathcal A}_i\rightarrow S_i$, leaving $S_i$ fixed, with ramification locus ${\mathcal R}_i\rightarrow S_i$ a smooth projective family defined over ${\overline{\mathbb Q}}$, whose fiber $R_i(s_i)$ at $s_i\in S_i$ consists of a union of smooth distinct hypersurfaces. For $s=(s_1,s_2)\in S=S_1\times S_2$, let $B(s)$ be the Calabi-Yau variety with involution constructed from the Calabi-Yau varieties with involution $A_1(s_1)$ and $A_2(s_2)$ using the method described before Proposition 2, and let ${\mathcal B}\rightarrow S$ be the smooth projective family whose fiber at $s\in S$ is $B(s)$. Let $R(s)$ be the ramification locus of the involution on $B(s)$ and let ${\mathcal R}\rightarrow S$ be the smooth projective family whose fiber at $s\in S$ is $R(s)$. We say that Expectation 2 holds for a smooth projective family ${\mathcal X}\rightarrow S$ defined over ${\overline{\mathbb Q}}$ when the set of fibers ${\mathcal X}_s$, $s\in S({\overline{\mathbb Q}})$, such that the Hodge filtration of $H^k({\mathcal X}_s,{\mathbb Q}_{{\mathcal X}_s})_{\mathbb C}$ is defined over ${\overline{\mathbb Q}}$ for all $k\le \dim({\mathcal X}_s)$, equals the set of $s\in S({\overline{\mathbb Q}})$ such that ${\mathcal X}_s$ has strong CM.

\noindent

\medskip

\noindent {\bf Theorem 1:} \emph{Assume that Expectation 2 holds for the families ${\mathcal A}_i\rightarrow S_i$ and ${\mathcal R}_i\rightarrow S_i$, $i=1,2$. Then, Expectation 2 holds for the family ${\mathcal B}\rightarrow S$. In other words, the set of $s=(s_1,s_2)\in S_1\times S_2$, with $s_i\in S_i({\overline{\mathbb Q}})$, $i=1,2$, such that the Hodge filtration of $H^k(B(s),{\mathbb Q}_{B(s)})_{\mathbb C}$ is defined over ${\overline{\mathbb Q}}$ for all $k\le \dim(B(s))$, equals the set of $s=(s_1,s_2)\in S_1\times S_2$, with $s_i\in S_i({\overline{\mathbb Q}})$, $i=1,2$, such that $B(s)$ has strong CM. Moreover, if for $s_i\in S_i({\overline{\mathbb Q}})$, $i=1,2$, the Hodge filtration of $H^k(B(s),{\mathbb Q}_{B(s)})_{\mathbb C}$ is defined over ${\overline{\mathbb Q}}$ for $k\le \dim(B(s))$, then $R(s)$ has strong CM.}

\medskip

\noindent {\bf Proof:} By the assumptions of the theorem, if $s_i\in S_i({\overline{\mathbb Q}})$, $i=1,2$, and if the Hodge filtrations of $H^k(A_i(s_i),{\mathbb Q}_{A_i(s_i)})_{\mathbb C}$, $k\le\dim(A_i)$, $H^k(R_i(s_i),{\mathbb Q}_{R(s_i)})_{\mathbb C}$, $k\le\dim(R_i)$, are defined over ${\overline{\mathbb Q}}$, then $A_i=A_i(s_i)$ and $R_i=R_i(s_i)$ have strong CM. Moreover, the varieties $B=B(s_1,s_2)$ and $R=R(s_1,s_2)$ are defined over ${\overline{\mathbb Q}}$. Suppose, for all $k\le\dim(B)$, that the Hodge filtration of $H^k(B,{\mathbb Q}_B)_{\mathbb C}$ is defined over ${\overline{\mathbb Q}}$. Then, by Proposition 2, the Hodge filtration of $H^k(A_i,{\mathbb Q}_{A_i})_{\mathbb C}$ is defined over ${\overline{\mathbb Q}}$ for all $k\le\dim(A_i)$, and the Hodge filtration of $H^k(R_i,{\mathbb Q}_{R_i})_{\mathbb C}$ is defined over ${\overline{\mathbb Q}}$ for all $k\le\dim(R_i)$, $i=1,2$. Therefore $R_i$ and $A_i$ have strong CM. By \cite{Roh}, Theorem 7.2.6 and Claim 7.2.7, it follows that $B$ and $R$ have strong CM, as required.

\noindent Iterating the general step up a Borcea-Voisin tower, we are able to construct infinitely many families of Calabi-Yau varieties, of arbitrary dimension, which satisfy Expectation 2. Each time we introduce a new family of Calabi-Yau varieties with involution into the construction, we must in general assume that the family of ramification loci of the involution also satisfies Expectation 2. For some examples in low dimension showing the necessity of the assumptions of Theorem 1, see \S\ref{s:exampleslow}.

\bigskip

\section{Viehweg-Zuo towers}\label{s:expect1}

\medskip

 \noindent In this section, we give examples of families of Calabi-Yau varieties for which the statement of Part 2 of the Problem of \S\ref{s:VHS} holds. They can serve as smooth projective families of Calabi-Yau varieties with involution ${\mathcal A}_1\rightarrow S_1$ and ${\mathcal A}_2\rightarrow S_2$ at the general inductive step of a Borcea-Voisin tower, as described in \S\ref{s:BVtower}. Our arguments entail fixing an algebraic point in the parameter space of the family and working with the projective variety given by the fiber above that point. It therefore suffices to show Expectation 2 of \S\ref{s:Schneider} for the set of such fibers.

 \noindent The varieties $X$ of this section are all curves, surfaces or Calabi-Yau 3-folds. As remarked in \S\ref{s:VHS},  when $X$ is a curve or a Calabi-Yau 3-fold, we have $H^n(X,{\mathbb Q})_{\rm prim}=H^n(X,{\mathbb Q})$, and when $X$ is a surface $H^2(X,{\mathbb Q})$ has CM if and only if $H^2(X,{\mathbb Q})_{\rm prim}$ has CM. We therefore work mainly with the full cohomology group $H^n(X,{\mathbb Q})$, $n=\dim(X)$, in this section and are still able to apply the results for primitive cohomology from \cite{VZ}.

\noindent In \cite{TCY}, we gave several examples of families of Calabi-Yau manifolds with involution for which Expectation 2 holds. One example, originating in work of Borcea, was a 2-step Borcea-Voisin tower of Calabi-Yau manifolds with elliptic curves as base. We then appealed to Th. Schneider's result \cite{Sch} that Expectation 2 is verified for elliptic curves. Another
 example was a quintic 3-fold arising as a 2-step iterated cyclic cover of a curve. We showed that if Expectation 2 is verified for the base curve, then it is verified for the iterated cyclic cover, so that we may apply the generalization of Schneider's Theorem in \cite{Co1} and \cite{SW}. Both examples fit into a family of Calabi-Yau manifolds with a dense set of CM fibers, and the arguments apply to any fiber over an algebraic point of the base of the family. We will revisit both examples in \S\ref{s:examples}. It turns out that Expectation 1 holds for our examples and that Expectation 2 is a consequence. Indeed, our first example in Theorem 2 of this section is a family of $K3$ surfaces, which have strong CM once they have CM in the sense of Definition 4 of \S\ref{s:VHS}. Our second example is a family of Calabi-Yau threefolds, and so again they have strong CM once they have CM.

 \noindent All our examples of suitable base manifolds for a Borcea-Voisin tower are taken from the families of Calabi-Yau manifolds with dense sets of CM fibers studied by Rohde in \cite{Roh}, which uses in part earlier work of Viehweg-Zuo \cite{VZ}. Some of the material from this section was derived from the discussion in \cite{HuS}. These families are the so-called Viehweg-Zuo towers, that is iterated cyclic covers of families of projective algebraic curves with dense sets of CM fibers. As pointed out by Rohde in \cite{Roh}, Chapter 7, the Borcea-Voisin towers and the Viehweg-Zuo towers are special cases of the same general construction described in that same reference. In all the examples of Viehweg-Zuo towers, it is shown in \cite{Bor}, \cite{Roh}, and \cite{VZ} that the Jacobian of the base curve has CM (we also say that the curve has CM) if and only if the Calabi-Yau variety given by its iterated cyclic cover has CM. We use this fact, and its proofs, to show that Expectation 1 and 2 hold for these examples.

 \noindent We now give the details. In \cite{Roh}, Rohde considers cyclic covers of ${\mathbb P}_1$ of the form
 $$
 y^m=x^{d_0}(x-a_1)^{d_1}(x-a_2)^{d_2}\ldots(x-a_n)^{d_n}(x-1)^{d_{n+1}}
 $$
 for integers $d_k$ with $1\le d_k\le m$, $k=0,\ldots,n+1$, and determines how many inequivalent such covers have sets of CM fibers which are dense in the Zariski topology of the parameter space: $a_i\in {\mathbb C}\setminus\{0,1\}$, $a_i\not= a_j$, $i\not= j$. The same question was solved independently by Moonen in \cite{Moo}. Moonen uses mixed characteristic methods and comments that a Hodge theoretic proof is desirable, which is exactly the approach used by Rohde. Some examples can also be found in the work of de Jong and Noot \cite{dJN}. When $m$ is prime, this is essentially work of Deligne-Mostow \cite{DM}, \cite{Mo}, as pointed out and exploited in the applications to transcendence of values of hypergeometric functions considered in \cite{CoWo1}, \cite{CoWo2}, \cite{DTT}, and \cite{TM}.

 \noindent The families of base curves with dense sets of CM fibers from the above list, and the number of their cyclic iterations, are chosen so that the resulting family of projective varieties are smooth hypersurfaces of degree $N+1$ in ${\mathbb P}_N$, and hence Calabi-Yau varieties. The list of suitable families of base curves is given in \cite{Roh}, \S7.3, p.152. For the smoothness of the families of Calabi-Yau varieties resulting from the iterated cyclic covers, we assume, with Rohde, that in addition the family ${\mathcal F}_0(m,n)$ of base curves has the form
 $$
 y^m=x(x-1)(x-a_1)(x-a_2)\ldots (x-a_n)
 $$
with $a_i\in{\mathbb C}$, $a_i\not=0,1$, $a_i\not= a_j$, $i\not=j$, where $m$ divides $n+3$. Together with requiring a dense set of CM fibers, this leaves the following possibilities,
$$
(m,n)=(2,1),\;(2,3),\;(3,3),\;(4,1),\;(5,2).
$$
The case $(2,1)$ gives a family of elliptic curves, which are Calabi-Yau varieties of dimension 1 with the involution $(x,y)\mapsto (x,-y)$, so there is no need for a cyclic iteration. In this example, Expectation 1 and 2 coincide and hold by Th. Schneider \cite{Sch}, as already discussed. For the case $(3,3)$, the construction of Viehweg-Zuo must be modified, and is treated in \cite{Roh}, Chapter 8. We do not treat this example in the present paper, and plan to come back to it in future work. We also do not treat the case $(2,3)$, for which Rohde remarks that it is not clear the Viehweg-Zuo construction can be made to work. The case $(4,1)$ yields a family of Calabi-Yau manifolds after one cyclic iteration of degree 4, as discussed in \cite{Roh}, \S7.4. Finally, the case $(5,2)$ is that considered in \cite{VZ}, where two cyclic iterations of degree 5 yield the family of quintic threefolds mentioned above. Although we treated this family in \cite{TCY}, here we give more detailed arguments. We therefore want to prove Expectation 1 and 2 for

\noindent (i) the family ${\mathcal F}_2(4,1)$ given by the zero locus in ${\mathbb P}_3$ of
$$
x_3^4+x_2^4+x_1(x_1-x_0)(x_1-a_1x_0)x_0,\qquad a_1\not=0,1,\infty.
$$
This is a family of cyclic covers of degree 4 of the family ${\mathcal F}_1(4,1)$ of curves in ${\mathbb P}_2$ given by the zero locus of
$$
x_2^4+x_1(x_1-x_0)(x_1-a_1x_0)x_0,\qquad a_1\not=0,1,\infty,
$$
which has CM for infinitely many $a_1\in{\overline{\mathbb Q}}$. The fiber of ${\mathcal F}_2(4,1)$ at $a_1$ carries both the involution $$I_2:[x_0:x_1:x_2:x_3]\mapsto [x_0:x_1:x_2:-x_3]$$ and $$I'_2:[x_0:x_1:x_2:x_3]\mapsto[x_0:x_1:x_3:x_2]$$ In both cases the divisor fixed by the involution is isomorphic to the curve in ${\mathcal F}_1(4,1)$ at $a_1$ (see \cite{Roh}, pp. 153-154).

\noindent (ii) the family ${\mathcal F}_3(5,2)$ given by the zero locus in ${\mathbb P}_4$ of
$$
x_4^5+x_3^5+x_2^5+x_1(x_1-x_0)(x_1-a_1x_0)(x_1-a_2)x_0,\qquad a_1,a_2\not=0,1,\infty,\;a_1\not=a_2.
$$
This is a family of cyclic covers of degree 5 of the family ${\mathcal F}_2(5,2)$ of surfaces in ${\mathbb P}_3$ given by the zero locus of
$$
x_3^5+x_2^5+x_1(x_1-x_0)(x_1-a_1x_0)(x_1-a_2x_0)x_0,\qquad a_1,a_2\not=0,1,\infty,\;a_1\not=a_2.
$$
This is in turn a cyclic cover of degree 5 of the family ${\mathcal F}_1(5,2)$ of curves in ${\mathbb P}_2$ given by the zero locus of
$$
x_2^5+x_1(x_1-x_0)(x_1-a_1x_0)(x_1-a_2x_0)x_0,,\qquad a_1,a_2\not=0,1,\infty,\;a_1\not=a_2,
$$
which has CM for a Zariski dense set of $(a_1, a_2)\in{\overline{\mathbb Q}}^2$. Each fiber of ${\mathcal F}_3(5,2)$ at $(a_1,a_2)$ carries the involution $$I_3:[x_0:x_1:x_2:x_3:x_4]\mapsto [x_0:x_1:x_2:x_4:x_3].$$ The divisor fixed by the involution is isomorphic to the surface given by the fiber of ${\mathcal F}_2(5,2)$ at $(a_1,a_2)$ (see \cite{Roh}, p.151).

\noindent We now consider the action of the $m$-th roots of unity on the fibers of the families ${\mathcal F}_d(m,n)$ above. In each case, the equation of the fibers, in projective coordinates, is of the form $x_j^m+f(x_0,\ldots,x_{j-1})=0$, with $f$ a homogeneous polynomial of degree $m$ without repeated roots. There is a natural action of $\zeta_m=\exp(2\pi i/m)$ on each fiber given by
$$[x_0:\ldots:x_{j-1}:x_j]\mapsto [x_0:\ldots:x_{j-1}:\zeta_mx_j],$$
which induces an action of $\zeta_m$ on the singular cohomology of the fibers with coefficients in ${\mathbb Q}(\zeta_m)$. This action commutes with the Hodge structure, and determines a decomposition of the singular cohomology $H^k(\;\cdot\;,{\mathbb Q}(\zeta_m))$, with coefficients in ${\mathbb Q}(\zeta_m)$, into eigenspaces $H^k(\;\cdot\;,{\mathbb Q}(\zeta_m))_i$ where $\zeta_m$ acts as $\zeta_m^i$, $i=1,\ldots,m-1$. As before, for a smooth complex projective variety $X$, we denote by $H^k(X,{\mathbb Q}_X)$ the usual Hodge structure on $H^k(X,{\mathbb Q})$ induced by the complex structure on $X$, $k\le\dim(X)$, and by $H^k(X,{\mathbb Q}_X)_i$ the induced Hodge structure on the ${\mathbb Q}(\zeta_m)$-eigenspace $H^k(\;\cdot\;,{\mathbb Q}(\zeta_m))_i$. Let $\Sigma^{(m)}$ denote the Fermat curve of degree $m$ given in homogeneous coordinates $[x_0:x_1:x_2]$ for ${\mathbb P}_2$ by  the equation
$$
x_2^m=x_1^m+x_0^m.
$$
 We fix a generator $g_m$ of $({\mathbb Z}/m{\mathbb Z})^\ast$ that acts as the automorphism of $\Sigma^{(m)}$ given by $g_m:[x_0:x_1:x_2]\mapsto [x_0:x_1:\zeta_m x_2].$ As above, let $H^1(\Sigma^{(m)},{\mathbb Q}_{\Sigma^{(m)}})_i$ be the eigenspace of $H^1(\Sigma^{(m)},{\mathbb Q}_{\Sigma^{(m)}})$ such that the induced action of $g_m$ is by $\zeta_m^i$. Remark that $H^a(X,{\mathbb Q}_X)_i\otimes H^b(\Sigma^{(m)},{\mathbb Q}_{\Sigma^{(m)}})_{m-i}$ has a ${\mathbb Q}$-structure.

\noindent By a rational Hodge structure $(W,\widetilde{\varphi})$ concentrated in bi-degree $(1,1)$, we mean a ${\mathbb Q}$-vector space $W$ together with a homomorphism of ${\mathbb R}$-algebraic groups $\widetilde{\varphi}:{\mathbb S}({\mathbb R})\rightarrow {\rm GL}(W)_{\mathbb R}$ such that for all $a,b\in{\mathbb R}$ with $a^2+b^2\not=0$ and all $w\in W$, we have
$$
\widetilde{\varphi}\left(\begin{pmatrix}{a&-b\cr b&a}\end{pmatrix}\right)(w)=(a^2+b^2)w.
$$
The smallest ${\mathbb Q}$-algebraic subgroup of ${\rm GL}(W)$ whose real points contain $\widetilde{\varphi}({\mathbb S}({\mathbb R}))$ is a torus with ${\mathbb Q}$-rational points isomorphic to the multiplicative group ${\mathbb Q}(\sqrt{-1})^\ast$, so that $(W,\widetilde{\varphi})$ has CM in this sense (see \S\ref{s:Schneider}). To $(W,\widetilde{\varphi})$ we associate the filtration
$$
\{0\}=F^2_{\mathbb C}\subset F^1_{\mathbb C}=F^0_{\mathbb C}=W_{\mathbb C},
$$
which is defined over ${\overline{\mathbb Q}}$ in the trivial sense that every ${\overline{\mathbb Q}}$-basis of $W_{\overline{\mathbb Q}}$ is automatically a ${\mathbb C}$-basis of $W_{\mathbb C}$.

\noindent Let $F_d^{(4)}(s)$ denote the fiber of ${\mathcal F}_d(4,1)$ at $s=a_1\in{\mathbb C}$, $a_1\not=0,1$, for $d=1,2$, and $F_d^{(5)}(s)$ denote the fiber of ${\mathcal F}_d(5,2)$ at $s=(a_1,a_2)\in{\mathbb C}^2$, $a_1,a_2\not=0,1$, $a_1\not=a_2$, for $d=1,2,3$.

\noindent We now look in more detail at the construction of the Viehweg-Zuo towers so as to make more evident the similarity to the Borcea-Voisin towers. We only treat the examples used in this section. For a more general and thorough presentation, see \cite{HuS}, \cite{Roh}, and \cite{VZ}. Fix $s\not=0,1$ and let $F^{(4)}_d=F^{(4)}_d(s)$, $d=1,2$. We write the equation for $F^{(4)}_1$ using homogeneous coordinates $[y_0:y_1:y_2]$ for ${\mathbb P}_2$ as
$$
y_2^4+y_1(y_1-y_0)(y_1-sy_0)y_0=0,
$$
and the equation for $\Sigma^{(4)}$ in homogeneous coordinates $[z_0:z_1:z_2]$ for ${\mathbb P}_2$ as
$$
z_2^4-z_1^4-z_0^4=0.
$$
The fixed points for the action of $\zeta_4$ on $F_1^{(4)}$ described above are those with $y_2=0$, namely the four points $[1:0:0]$, $[1:1:0]$, $[1:s:0]$, $[0:1:0]$. The fixed points for the action of $\zeta_4$ on $\Sigma^{(4)}$ are those with $z_2=0$, namely the four points $[1:1:0]$, $[1:\zeta_4:0]$, $[1:\zeta_4^2:0]$, $[1:\zeta_4^3:0]$. Consider the rational map, which is a morphism outside the set of 16 points given by $y_2=z_2=0$, given by
$$
\Phi^{(4)}:F_1^{(4)}\times \Sigma^{(4)}\dashrightarrow F_2^{(4)}\subset{\mathbb P}_2
$$
$$
([y_0:y_1:y_2],[z_0:z_1:z_2])\mapsto [x_0:x_1:x_2:x_3]=[z_2y_0:z_2y_1:y_2z_0:y_2z_1].
$$
It is easy to check directly that, where $\Phi^{(4)}$ is defined, its image indeed lies in $F_2^{(4)}$. Consider the action $\eta_2^{(4)}$ of $\zeta_4$ on $F_1^{(4)}\times \Sigma^{(4)}$ given by
$$
\eta_2^{(4)}:([y_0:y_1:y_2],[z_0:z_1:z_2])\mapsto([y_0:y_1:\zeta_4y_2],[z_0:z_1:\zeta_4 z_2]).
$$
The fixed points of $\eta_2^{(4)}$ are the 16 points on $F_1^{(4)}\times \Sigma^{(4)}$ with $y_2=z_2=0$. The points of $F_1^{(4)}\times \Sigma^{(4)}$ in the same orbit of $\eta_2^{(4)}$ have the same image in $F_2^{(4)}$, and the map $\Phi^{(4)}$, which is generically $4$ to $1$, induces a birational map
 $$
(F_1^{(4)}\times\Sigma^{(4)})/\langle\eta_2^{(4)}\rangle\dashleftarrow\dashrightarrow F_2^{(4)}.
 $$
 The smooth surface $F_2^{(4)}$ is a $K3$-surface, so in particular is Calabi-Yau.

 \noindent We can blow up $Y^{(4)}_2=F_1^{(4)}\times\Sigma^{(4)}$ at the 16 points fixed by $\eta_2^{(4)}$ to obtain a surface ${\widehat{Y}}^{(4)}_2$. On ${\widehat{Y}}^{(4)}_2$, we have an automorphism ${\widehat{\eta}}_2^{(4)}$, induced by $\eta_2^{(4)}$, whose ramification locus is the 16 exceptional divisors. The quotient ${\widehat{Y}}^{(4)}_2/{\widehat{\eta}}_2^{(4)}$ by the action of ${\widehat{\eta}}_2^{(4)}$ is a smooth surface. On the other hand, the quotient surface $(F_1^{(4)}\times\Sigma^{(4)})/\langle\eta_2^{(4)}\rangle$ has exactly 16 singular points coming from the 16 fixed points of $\eta_2^{(4)}$. We can resolve these singularities by blowing them up. By \cite{Roh}, Chapter 7, and \cite{VZ}, the resulting smooth surface is isomorphic to ${\widehat{Y}}^{(4)}_2/{\widehat{\eta}}_2^{(4)}$. By \cite{VZ}, Lemma 6.3, the surface $F_2^{(4)}$ is obtained from ${\widehat{Y}}^{(4)}_2/{\widehat{\eta}}_2^{(4)}$ by two successive blow-downs, introducing Hodge structures concentrated in bidegree $(1,1)$ and independent of $s$ which we ignore for the time being (see \cite{VZ}, Claim 7.1, Claim 7.2, and the comment below Corollary 8.4). Modulo these Hodge structures, it suffices to look at the Hodge structure
$$
H^2({\widehat{Y}}^{(4)}_2/{\widehat{\eta}}_2^{(4)},{\mathbb Q}_{{\widehat{Y}}^{(4)}_2/{\widehat{\eta}}_2^{(4)}})\simeq H^2({\widehat{Y}}_2^{(4)},{\mathbb Q}_{{\widehat{Y}}_2^{(4)}})^{{\widehat{\eta}}_2^{(4)}},
$$
where the superscript on the right hand side denotes the ${\widehat{\eta}}_2^{(4)}$-invariant part of $H^2({\widehat{Y}}_2^{(4)},{\mathbb Q}_{{\widehat{Y}}_2^{(4)}})$, which carries the rational Hodge substructure induced by $H^2({\widehat{Y}}_2^{(4)},{\mathbb Q}_{{\widehat{Y}}_2^{(4)}})$ (see \cite{Gro}). Using this together with \S\ref{s:lemmas}, Lemma 3, and the fact that ${\widehat{Y}}_2^{(4)}$ is obtained from $F_1^{(4)}\times \Sigma^{(4)}$ by blowing up 16 points fixed by $\eta_2^{(4)}$, we deduce that, up to a Hodge structure concentrated in bi-degree $(1,1)$,
$$
H^2({\widehat{Y}}_2^{(4)},{\mathbb Q}_{{\widehat{Y}}_2^{(4)}})^{{\widehat{\eta}}_2^{(4)}}\simeq H^2(F_1^{(4)}\times\Sigma^{(4)},{\mathbb Q}_{F_1^{(4)}\times\Sigma^{(4)}})^{\eta_2^{(4)}},
$$
where now the superscript denotes the $\eta_2^{(4)}$-invariant part. By \S\ref{s:lemmas}, Lemma 2, we have
$$
H^2(F_1^{(4)}\times\Sigma^{(4)},{\mathbb Q}_{F_1^{(4)}\times\Sigma^{(4)}})^{\eta_2^{(4)}}=\bigoplus_{r+s=2}\bigoplus_{i=0}^3H^r(F_1^{(4)},{\mathbb Q}_{F_1^{(4)}})_i\otimes H^s(\Sigma^{(4)},{\mathbb Q}_{\Sigma^{(4)}})_{4-i}.
$$
When $i\not=0$, we have
$$
H^0(F_1^{(4)},{\mathbb Q}_{F_1^{(4)}})_i=H^2(\Sigma^{(4)},{\mathbb Q}_{\Sigma^{(4)}})_{4-i}=\{0\}.
$$
Moreover,
$$
H^1(F_1^{(4)},{\mathbb Q}_{F_1^{(4)}})_0=H^1(\Sigma^{(4)},{\mathbb Q}_{\Sigma^{(4)}})_0=\{0\}.
$$
Finally, the $\eta_2^{(4)}$-invariant part of
$$
H^r(F_1^{(4)},{\mathbb Q}_{F_1^{(4)}})_0\otimes H^{2-r}(\Sigma^{(4)},{\mathbb Q}_{\Sigma^{(4)}})_0
$$
when $r=0,2$ is concentrated in bi-degree $(1,1)$ and is independent of $s$. We now put in the Hodge structures ${\mathbb W}$ and ${\mathbb W}'$ mentioned in the above arguments, concentrated in bi-degree $(1,1)$ and independent of $s$, to deduce that
$$
H^2(F_2^{(4)},{\mathbb Q}_{F_2^{(4)}})\oplus{\mathbb W}'=\bigoplus_{i=1}^3(H^1(F_1^{(4)},{\mathbb Q}_{F_1^{(4)}})_i\otimes H^1(\Sigma^{(4)},{\mathbb Q}_{\Sigma^{(4)}})_{4-i})\oplus{\mathbb W}.
$$

\noindent We now turn to $F^{(5)}_d=F^{(5)}_d(s)$, where $s=(a_1,a_2)$, $a_1\not=a_2$, $a_1,a_2\not=0$, $1$, $\infty$ is fixed, and $d=2,3$,
and look at the Hodge structure $H^3(F_3^{(5)},{\mathbb Q}_{F_3^{(5)}})$. We write the equation for $F_2^{(5)}$ using homogeneous coordinates $[y_0:y_1:y_2:y_3]$ for ${\mathbb P}_3$ as
$$
y_3^5+y_2^5+y_1(y_1-y_0)(y_1-a_1y_0)(y_1-a_2y_0)y_0=0.
$$
The fixed points for the action of $\zeta_5$ described above are those with $y_3=0$, whose locus is isomorphic to $F_1^{(5)}$, which is given by the equation,
$$
y_2^5+y_1(y_1-y_0)(y_1-a_1y_0)(y_1-a_2y_0)y_0=0.
$$
We write the equation for $\Sigma_5$ in ${\mathbb P}_2$ as before,
$$
z_2^5-z_1^5-z_0^5=0.
$$
The fixed points for the action of $\zeta_5$ on $\Sigma^{(5)}$ described above are the elements of $\{[1:\zeta_5^i:0]; i=0,\ldots,4\}$.
We have a rational map, which is a morphism outside the locus of points with $y_3=z_2=0$, given by
$$
\Phi^{(5)}: F_2^{(5)}\times\Sigma^{(5)}\dashrightarrow F_3^{(5)},
$$
where $\Phi^{(5)}$ maps
$$
([y_0:y_1:y_2:y_3],[z_0:z_1:z_2])
$$
to
$$
[x_0:x_1:x_2:x_3:x_4]=[z_2y_0:z_2y_1:z_2y_2:y_3z_0:y_3z_1].
$$
A simple computation shows that the image of this map indeed lies on $F_3^{(5)}$. Consider the action $\eta_3^{(5)}$ of $\zeta_5$ on $F_2^{(5)}\times \Sigma^{(5)}$, given by
$$
([y_0:y_1:y_2:y_3],[z_0:z_1:z_2])\mapsto ([y_0:y_1:y_2:\zeta_5y_3],[z_0:z_1:\zeta_5z_2]).
$$
The set of fixed points of $\eta^{(5)}_3$ is the union $Z$ of the 5 distinct curves
$$
F_1^{(5)}\times [1:\zeta_5^i:0],\qquad i=0,\ldots,4.
$$
The points of $F_2^{(5)}\times \Sigma^{(5)}$ in the same orbit of $\eta_3^{(5)}$ have the same image under $\Phi^{(5)}$, which is a generically $5$ to 1 map inducing a birational map
$$
(F_2^{(5)}\times\Sigma^{(5)})/\langle \eta^{(5)}_2\rangle\dashleftarrow\dashrightarrow F_3^{(5)}.
$$
The variety $F_3^{(5)}$ is a Calabi-Yau 3-fold and a resolution of $(F_2^{(5)}\times\Sigma^{(5)})/\langle \eta^{(5)}_2\rangle$. Denote by ${\widehat{Y}}_3^{(5)}$ the blow up of $F_2^{(5)}\times\Sigma^{(5)}$ along $Z$, and ${\widehat{\eta}}_3^{(5)}$ the automorphism of ${\widehat{Y}}_3^{(5)}$ induced by $\eta_3^{(5)}$. By \cite{Roh}, Chapter 7, and \cite{VZ}, the smooth quotient ${\widehat{Y}}_3^{(5)}/{\widehat{\eta}}_3^{(5)}$ is isomorphic to the resolution of $(F_2^{(5)}\times\Sigma^{(5)})/\langle \eta^{(5)}_3\rangle$ obtained by blowing up its singular locus. These are in turn isomorphic to the 3-fold obtained by two successive blow-ups of the Calabi-Yau threefold $F_3^{(5)}$, the first blow-up not affecting the Hodge structure and the second having center a curve isomorphic to $F_1^{(5)}$ by \cite{VZ}, 6.3(v), and Claim 7.1. By \S\ref{s:lemmas}, Lemma 3, we have an isomorphism of Hodge structures
$$
H^3({\widehat{Y}}_3^{(5)},{\mathbb Q}_{{\widehat{Y}}_3^{(5)}})\simeq H^3(F_2^{(5)}\times\Sigma^{(5)},{\mathbb Q}_{F_2^{(5)}\times\Sigma^{(5)}})\bigoplus
H^1(Z,{\mathbb Q}_Z)(-1),
$$
where $H^1(Z,{\mathbb Q}_Z)(-1)$ means $H^1(Z,{\mathbb Q}_Z)$ shifted by $(1,1)$ in bidegree. As $Z$ is fixed by $\eta_3^{(5)}$, we have
$$
H^3({\widehat{Y}}_3^{(5)},{\mathbb Q}_{{\widehat{Y}}_3^{(5)}})^{{\widehat{\eta}}_3^{(5)}}\simeq H^3(F_2^{(5)}\times\Sigma^{(5)},{\mathbb Q}_{F_2^{(5)}\times\Sigma^{(5)}})^{\eta_3^{(5)}}\bigoplus
H^1(Z,{\mathbb Q}_Z)(-1)
$$
$$
=\bigoplus_{r+s=3}\bigoplus_{i=0}^4(H^r(F_2^{(5)},{\mathbb Q}_{F_2^{(5)}})_i\otimes H^s(\Sigma^{(5)},{\mathbb Q}_{\Sigma^{(5)}})_{5-i})\bigoplus^5H^1(F_1^{(5)},{\mathbb Q}_{F_1^{(5)}})(-1),
$$
where $\oplus^5H^1(F_1^{(5)},{\mathbb Q}_{F_1^{(5)}})$ means 5 copies of $H^1(F_1^{(5)},{\mathbb Q}_{F_1^{(5)}})$. This term comes from the fact that $Z$ is the union of 5 distinct subvarieties isomorphic to $F_1^{(5)}$.  We have $H^2(\Sigma^{(5)},{\mathbb Q}_{\Sigma^{(5)}})_{5-i}=H^0(\Sigma^{(5)},{\mathbb Q}_{\Sigma^{(5)}})_{5-i}=\{0\}$, $i\not=0$, and $H^1(F_2^{(5)},{\mathbb Q}_{F_2^{(5)}})_0=H^3(F_2^{(5)},{\mathbb Q}_{F_2^{(5)}})_0=\{0\}$. Therefore the Hodge structure $H^3({\widehat{Y}}_3^{(5)},{\mathbb Q}_{{\widehat{Y}}_3^{(5)}})^{{\widehat{\eta}}_3^{(5)}}$ is isomorphic to
$$
\bigoplus_{i=1}^4(H^2(F_2^{(5)},{\mathbb Q}_{F_2^{(5)}})_i\otimes H^1(\Sigma^{(5)},{\mathbb Q}_{\Sigma^{(5)}})_{5-i})\bigoplus^{5}H^1(F_1^{(5)},{\mathbb Q}_{F_1^{(5)}})(-1).
$$
Allowing for the fact mentioned above that $F_3^{(5)}$ is a blow-down of ${\widehat{Y}}_3^{(5)}/{\widehat{\eta}}_3^{(5)}$, we have that the Hodge structure $H^3(F_3^{(5)},{\mathbb Q}_{F_3^{(5)}})$ is isomorphic to
$$
\bigoplus_{i=1}^4(H^2(F_2^{(5)},{\mathbb Q}_{F_2^{(5)}})_i\otimes H^1(\Sigma^{(5)},{\mathbb Q}_{\Sigma^{(5)}})_{5-i})\bigoplus^{4}H^1(F_1^{(5)},{\mathbb Q}_{F_1^{(5)}})(-1).
$$

\noindent The above discussion yields the following special cases of results of Viehweg-Zuo (\cite{VZ}, Proposition 7.4, and \S8) on general iterated cyclic covers with base a cyclic cover of ${\mathbb P}_1$.

\medskip

\noindent {\bf Proposition 3:} There are Hodge structures ${\mathbb W}$ and ${\mathbb W}'$, independent of $s$ and concentrated in bi-degree $(1,1)$, such that we have an isomorphism of rational Hodge structures,
$$
H^2(F_2^{(4)},{\mathbb Q}_{F_2^{(4)}})\oplus{\mathbb W}'\simeq \bigoplus_{i=1}^{3}( H^1(F_1^{(4)},{\mathbb Q}_{F_1^{(4)}})_i\otimes H^1(\Sigma^{(4)},{\mathbb Q}_{\Sigma^{(4)}})_{4-i})\oplus {\mathbb W}.
$$
The Hodge structure $H^3(F_3^{(5)},{\mathbb Q}_{F_3^{(5)}})$ is isomorphic to
$$\bigoplus_{i=1}^{4}(H^2(F_2^{(5)},{\mathbb Q}_{F_2^{(5)}})_i\otimes H^1(\Sigma^{(5)},{\mathbb Q}_{\Sigma^{(5)}})_{5-i})\oplus\bigoplus^{4}H^1(F_1^{(5)},{\mathbb Q}_{F_1^{(5)}})(-1).
$$
Let ${\mathcal X}\rightarrow S$ be a smooth projective family defined over ${\overline{\mathbb Q}}$, in the sense of \S\ref{s:BVtower}, whose fiber ${\mathcal X}_s$ at $s\in S$ has dimension $n$ (independent of $s$). We say that Expectation 1 holds for this family when the set of fibers ${\mathcal X}_s$, $s\in S({\overline{\mathbb Q}})$, such that the Hodge filtration of $H^n({\mathcal X}_s,{\mathbb Q}_{{\mathcal X}_s})_{\mathbb C}$ is defined over ${\overline{\mathbb Q}}$ equals the set of $s\in S({\overline{\mathbb Q}})$ such that $H^n({\mathcal X}_s,{\mathbb Q}_{{\mathcal X}_s})$ has CM.

\medskip

\noindent{\bf Theorem 2:}  Expectation 1 holds for the families ${\mathcal F}_2(4,1)$
 and ${\mathcal F}_3(5,2)$.

\medskip

\medskip

\noindent {\bf Proof:} Recall the notations and definitions from \S\ref{s:Schneider}. Let $X$ be a smooth projective variety of dimension $n$, and let $F^\ast_\Omega(X)$
be the Hodge filtration of $H^n_\Omega(X,{\mathbb C})$. Let $F^\ast(X)=\iota_1^{-1}(F^\ast_\Omega(X))$ be the induced filtration of $H^n(X,{\mathbb C})$,
and call it the Hodge filtration of $H^n(X,{\mathbb Q}_X)_{\mathbb C}$. We say that the Hodge filtration $F^\ast(X)$ of $H^n(X,{\mathbb Q}_X)_{\mathbb C}$
is defined over ${\overline{\mathbb Q}}$ if there is a ${\overline{\mathbb Q}}$-filtration of $H^n(X,{\overline{\mathbb Q}})$
inducing $F^\ast(X)$ by extension of scalars to ${\mathbb C}$. This is equivalent to the Hodge filtration $F_\Omega^\ast(X)$ of $H^n_\Omega(X,{\mathbb C})$ being
defined over ${\overline{\mathbb Q}}$, in the sense of \S\ref{s:Schneider}, Definition 2.

\noindent We first show that Expectation 1 holds for the family ${\mathcal F}_2(4,1)$, whose fiber at $s=a_1\in S_1={\mathbb C}\setminus\{0,1\}$ is $F_2^{(4)}=F_2^{(4)}(s)$. By \cite{Co1}, \cite{SW}, Expectation 1 holds for the family ${\mathcal F}_1(4,1)$ whose fiber at $s\in S_1$ is $F^{(4)}_1(s)$.
In other words, if $s$ has algebraic coordinates, and if the Hodge filtration $F^\ast(F^{(4)}_1(s))$
is defined over ${\overline{\mathbb Q}}$, then $F_1^{(4)}(s)$ has CM. Recall that, for a curve $X$, by this we mean that the rational Hodge structure $H^1(X,{\mathbb Q}_X)$ has CM.

\noindent Fix a parameter $s=\alpha$ with algebraic coordinates. Let $V=H^2(F_2^{(4)}(\alpha),{\mathbb Q})$, $W_1=H^1(F_1^{(4)}(\alpha),{\mathbb Q})$, $W_2=H^1(\Sigma^{(4)},{\mathbb Q})$. It is easy to see, for all $s$, that  $H^{2,0}(F_2^{(4)}(s))$ is the tensor product of the $\zeta_4^2$-eigenspace of $H^{1,0}(F_1^{(4)}(s))$, and the $\zeta_4^2$-eigenspace of $H^{1,0}(\Sigma^{(4)})$. Of course, both these are 1-dimensional, since $H^{2,0}(F_2^{(4)}(s))$ is also. By assumption, we have $
H^{2,0}(F_2^{(4)}(\alpha))={\mathbb C}(\Omega_2\otimes_{\overline{\mathbb Q}}1_{\mathbb C})$ where $\Omega_2$ is in $H^2(F_2^{(4)}(\alpha),{\overline{\mathbb Q}})$. As $\Sigma^{(4)}$ has complex multiplication, the $\zeta_4^{(2)}$-eigenspace of $H^{1,0}(\Sigma^{(4)})$ is of the form ${\mathbb C}(\Omega_1^{(4)}\otimes_{\overline{\mathbb Q}}1_{\mathbb C})$, for $\Omega_1\in H^1(\Sigma^{(4)},{\overline{\mathbb Q}})$.

Therefore, it is clear, for example by \S\ref{s:lemmas}, Lemma 6, that the $\zeta_4^2$-eigenspace of $H^{1,0}(F_1^{(4)}(\alpha))$ is of the form ${\mathbb C}(\Omega_1\otimes_{\overline{\mathbb Q}}1_{\mathbb C})$
for $\Omega_1$ in $H^1(F_1^{(4)}(\alpha),{\overline{\mathbb Q}})$. Therefore $F_1^{(4)}(\alpha)$ has CM and so by \cite{VZ}, it follows that $F_2^{(4)}(\alpha)$ has CM as required.

We now show that Expectation 1 holds for ${\mathcal F}_3(5,2)$. Fix a parameter $s=\alpha$ with algebraic coordinates. As by Proposition 3, the Hodge structure $H^1(F_1^{(5)}(\alpha),{\mathbb Q}_{F_1^{(5)}(\alpha)})(-1)$ is a direct summand of $H^3(F_3^{(5)}(\alpha),{\mathbb Q}_{F_3^{(5)}(\alpha)})$, it is clear that if the Hodge filtration of $H^3(F_3^{(5)}(\alpha),{\mathbb Q}_{F_3^{(5)}(\alpha)})_{\mathbb C}$ is defined over ${\overline{\mathbb Q}}$, then so is that of $H^1(F_1^{(5)}(\alpha),{\mathbb Q}_{F_1^{(5)}(\alpha)})_{\mathbb C}$, see also \S\ref{s:lemmas}, Lemma 4. Therefore $F_1^{(5)}(\alpha)$ has CM, and so $F_3^{(5)}(\alpha)$ has CM by \cite{VZ}. This completes the proof of Theorem 2.

\medskip

\noindent Finally, we briefly discuss the fixed loci of the involutions $I_2$ and $I_2'$ on $F_2^{(4)}$ and $I_3$ on $F_3^{(5)}$. In order for the families of this section to be part of a Borcea-Voisin tower construction, we have to check that these loci also satisfy Expectation 2. The fixed locus of both $I_2$ and $I_2'$ is isomorphic to $F_1^{(4)}$ for which Expectation 1 and 2 are equivalent and hold by \cite{Co1}, \cite{SW}. The fixed locus of $I_3$ is isomorphic to $F_2^{(5)}$ for which Expectation 2 also holds. For another approach via explicit computation of periods, see \S\ref{s:examples}.

\bigskip

\section{Examples of periods}\label{s:examples}

\bigskip

\noindent In this section, we reconsider some of the examples of the present paper from a more explicit viewpoint. We show that Examples 1, 3, and 4 below satisfy the following variant of Expectation 1 in which we consider just the $F^{n,n}:=H^{n,0}$-part of the Hodge filtration of a Calabi-Yau variety of dimension $n$. Notice that the $F^{n,n}$ part of $H^n(X,{\mathbb Q}_X)_{\mathbb C}$ and $H^n(X,{\mathbb Q}_X)_{{\rm prim},{\mathbb C}}$ are the same. For the iterated cyclic covers $F^{(4)}_2$, $F^{(5)}_3$ of \S\ref{s:expect1}, we recover in this way a more direct proof of Expectation 1. For a family of $K3$ surfaces obtained by applying the Borcea-Voisin construction of \S\ref{s:BVtower} to two base families of elliptic curves, we obtain a proof of Expectation 1.

\medskip

\noindent {\bf Expectation 3:} Let $X$ be a Calabi-Yau variety of dimension $n$, defined over ${\overline{\mathbb Q}}$. If the $F^{n,n}:=H^{n,0}(X)$-part of the Hodge filtration of $H^n(X,{\mathbb Q}_X)_{{\rm prim},\mathbb C}$ is defined over ${\overline{\mathbb Q}}$, then $X$ has CM. In other words, if $\Omega_n$ is a nowhere vanishing holomorphic $n$-form on $X$, defined over ${\overline{\mathbb Q}}$, whose normalized periods are all algebraic, then $X$ has CM.

\medskip

\noindent In the above statement, by the $F^{n,n}$-part being defined over ${\overline{\mathbb Q}}$ we mean that $F^{n,n}:=H^{n,0}=F^{n,n}_{\overline{\mathbb Q}}\otimes_{\overline{\mathbb Q}}{\mathbb C}$, for a $\overline{\mathbb Q}$-vector space $F^{n,n}_{\overline{\mathbb Q}}\subseteq H^n(X,{\overline{\mathbb Q}})$. By the normalized periods being algebraic we mean that the unnormalized periods $\int_\gamma\Omega_n$, as $\gamma$ ranges over the $n$-cycles on $X$, form a vector space of dimension 1 over ${\overline{\mathbb Q}}$.

\medskip

\noindent For a Calabi-Yau variety $X$ of dimension $n$, define the \emph{transcendental part} $T(X)$ of $H^n(X,{\mathbb Q}_X)_{\rm prim}$ to be the smallest rational Hodge substructure of $H^n(X,{\mathbb Q}_X)_{\rm prim}$ whose underlying ${\mathbb Q}$-vector space $T$ satisfies $H^{n,0}\subseteq T_{\mathbb C}$, where $T_{\mathbb C}=T\otimes_{\mathbb Q}{\mathbb C}$. If we work over ${\mathbb Z}$, instead of ${\mathbb Q}$, the corresponding notion is that of transcendental lattice. This notion has been extensively studied for $K3$-surfaces, in particular in the study of complex multiplication \cite{Zar}, and appears in the discussion of Calabi-Yau threefolds with CM in \cite{Bor}. Let $Q_T$ be the restriction to $T$ of the polarization $Q$ on $H^n(X,{\mathbb Q}_X)_{\rm prim}$. The resulting polarized Hodge structure $(T(X),Q_T)$ is irreducible, and the ${\mathbb Q}$-algebra of ${\mathbb Q}$-endomorphisms of $T$ that commute with the Mumford-Tate group is a number field. This can be easily seen by adapting its proof for $n=2,3$ in the above references. Consider the following variant of a generalized Schneider's theorem.

\medskip

\noindent {\bf Expectation 4:} Let $X$ be a Calabi-Yau variety of dimension $n$, defined over ${\overline{\mathbb Q}}$. If the $F^{n,n}$-part of the Hodge filtration of $H^n(X,{\mathbb Q}_X)_{{\rm prim},\mathbb C}$ is defined over ${\overline{\mathbb Q}}$, then $T(X)$ has CM.

\medskip

\noindent In Example 2 below, we show that if two families of Calabi-Yau varieties with involution satisfy Expectation 4, then the family of Calabi-Yau varieties in the next step of a Borcea-Voisin tower also satisfies Expectation 4. To show the strong CM property at each step of a Borcea-Voisin tower, we in general need the extra assumptions of \S\ref{s:BVtower}, including those on the ramification loci of the involutions. For some examples, see \S\ref{s:exampleslow}.

\medskip

\begin{enumerate}

\item We recover directly, in this example, the more general arguments of \S\ref{s:BVtower}. For more details, see \cite{BorK3CM}, \cite{Huy}, \cite{Sh2}, \cite{Zar}. A Calabi-Yau variety of dimension 1 is an elliptic curve. For $i=1,2$, let $\lambda_i$ be complex numbers not equal to $0,1,\infty$. Let $E_i=E_{\lambda_i}$ be the elliptic curve with affine equation
$$
y^2=x(x-1)(x-\lambda_i),
$$
and let $I_i$ be the involution $(x,y)\mapsto(x,-y)$ on $E_i$. It has fixed points the 4 torsion points of order 2 on $E_i$. If we apply the Borcea-Voisin construction of \S\ref{s:BVtower} with base varieties with involution $(E_1,I_1)$, $(E_2,I_2)$, the Calabi-Yau surface $K_{1,2}$, with involution $\sigma_{1,2}$, at the next step of the tower is the \emph{Kummer surface} associated to the abelian surface $E_1\times E_2$ (see, for example, \cite{Bad}, \S10.5). It is a $K3$ surface. The elliptic curves are complex tori and we can represent them as $E_1={\mathbb C}/({\mathbb Z}\omega_{11}+{\mathbb Z}\omega_{12})$ and $E_2={\mathbb C}/({\mathbb Z}\omega_{21}+{\mathbb Z}\omega_{22})$, where the $\omega_{1j}$, $\omega_{2k}$, $j,k=1,2$ are periods of the holomorphic 1-form $dx/y$ on $E_1$, $E_2$ respectively. Let $\gamma_{11}$, $\gamma_{12}$ be a basis of $H_1(E_{\lambda_1}, {\mathbb Q})$ corresponding to $\omega_{11}$, $\omega_{12}$, and $\gamma_{21}$, $\gamma_{22}$ be a basis of $H_1(E_{\lambda_2},{\mathbb Q})$ corresponding to $\omega_{21}$, $\omega_{22}$. That is $\omega_{1j}=\int_{\gamma_{1j}}dx/y$, $\omega_{2k}=\int_{\gamma_{2j}}dx/y$, $j,k=1,2$. Let ${\widehat Y}_{1,2}$ be the blow up of  $E_1\times E_2$ at the 16 fixed points of $I_{1,2}=I_1\times I_2$ and let ${\widehat{I}}_{1,2}$ be the involution on ${\widehat Y}_{1,2}$ induced by $I_{1,2}$. Let $\pi$ be the canonical projection
${\widehat Y}_{1,2}\rightarrow {\widehat Y}_{1,2}/{\widehat{I}}_{1,2}\simeq K_{1,2}$. The six distinct 2-cycles $\gamma_{ij}\times\gamma_{k\ell}$ on $E_1\times E_2$, where $\{ij\}\not=\{k\ell\}$, lift to ${\widehat Y}_{1,2}$, and then descend via $\pi$ to six 2-cycles $\gamma_{ijk\ell}$ on $K_{1,2}$. We also have sixteen 2-cycles $\delta_m$, $m=1,\ldots, 16$, on $K_{1,2}$ coming from the 16 exceptional divisors on the blow-up ${\widehat{Y}}_{1,2}$ of $E_1\times E_2$. Together the $\gamma_{ijk\ell}$ and $\delta_m$ form a basis of the 22-dimensional ${\mathbb Q}$-vector space $H_2(K_{1,2},{\mathbb Q})$. The holomorphic 2-form $dz_1\wedge dz_2$ on ${\mathbb C}^2$ induces a holomorphic 2-form on $E_1\times E_2$ and also a holomorphic 2-form $\Omega$ on $K_{1,2}$, which is the unique holomorphic 2-form on $K_{1,2}$ up to a constant factor. If $\lambda_1,\lambda_2\in{\overline{\mathbb Q}}$, this 2-form is defined over ${\overline{\mathbb Q}}$. As $\pi$ is a 2-1 map, we have the following \emph{unnormalized periods} of $\Omega$,
$$
\int_{\gamma_{1j2k}}\Omega=\frac12\int_{\gamma_{1j}\times\gamma_{2k}}dz_1\wedge dz_2=\frac12\omega_{1j}\omega_{2k}=\frac12\omega_{11}\omega_{21}\tau_{1j}\tau_{2k},
$$
which is also a period of a holomorphic 2-form on the abelian surface $E_1\times E_2$. Here, the normalized periods on the elliptic curves are given by $\tau_{11}=1$, $\tau_{12}=\omega_{12}/\omega_{11}$ on $E_1$, and $\tau_{21}=1$, $\tau_{22}=\omega_{22}/\omega_{21}$ on $E_2$. Observe that the numbers $1$, $\tau_{12}$, $\tau_{22}$, as well as the product $\tau_{12}\tau_{22}$, appear as normalized periods of the 2-form $\Omega$. The transcendence of the unnormalized periods $\omega_{1j}\omega_{2k}$, when $\omega_{1j}\not=\omega_{2k}$, is not known.

As in the paragraph preceding the statement of Expectation 4, the transcendental part $T(K_{1,2})$ of $H^2(K_{1,2},{\mathbb Q}_{K_{1,2}})_{\rm prim}$ is defined to be the smallest rational Hodge substructure whose underlying ${\mathbb Q}$-vector space $T$ satisfies $H^{2,0}\subseteq T_{\mathbb C}$. Since the Mumford-Tate group acts trivially on $H^{1,1}$, the Hodge structure on $H^2(K_{1,2},{\mathbb Q})$ has CM if and only if the Hodge structure $T(K_{1,2})$ has CM, which is determined by what happens on the $H^{2,0}$-part. The Hodge structure $T(K_{1,2})$ is irreducible and $H^{2,0}$ is a 1-dimensional subspace of $T_{\mathbb C}$. It can be shown that $H^2(K_{1,2},{\mathbb Q})=T\oplus{\rm Pic}(K_{1,2})_{\mathbb Q}$, where ${\rm Pic}(K_{1,2})_{\mathbb Q}= H^{1,1}\cap H^2(K_{1,2},{\mathbb Q})$ (see \cite{Huy}, \S2). The 16 exceptional curves on $K_{1,2}$ determine 16 distinct linearly independent classes in ${\rm Pic}(K_{1,2})_{\mathbb Q}$, so that $T$ has dimension at most 6. The $K3$ surface is called singular if the dimension of ${\rm Pic}(K_{1,2})_{\mathbb Q}$ is maximal. In this case $\dim_{\mathbb Q}({\rm Pic}(K_{1,2})_{\mathbb Q})=\dim_{\mathbb C}(H^{1,1})=20$ and $\dim_{\mathbb Q}T=2$, with $T_{\mathbb C}=H^{2,0}\oplus H^{0,2}$. For example $K_{1,2}$ is singular when $E_1$, $E_2$ are isogenous curves with CM.

\noindent If the $H^{2,0}$-part of the Hodge filtration of $H^2(K_{1,2},{\mathbb Q}_{K_{1,2}})$ is defined over ${\overline{\mathbb Q}}$, the normalized periods $\tau_{12}$, $\tau_{22}$, $\tau_{12}\tau_{22}$ of $\Omega$ must be algebraic numbers. In particular, the normalized periods $\tau_{12}$, $\tau_{22}$ of the 1-form $dx/y$ on both the elliptic curves $E_1$, $E_2$ must be algebraic numbers. If, in addition, $\lambda_1$, $\lambda_2$ are algebraic numbers, the curves $E_1$, $E_2$ both have CM by Schneider's theorem \cite{Sch}, so that $\tau_{12}$, $\tau_{22}$ are quadratic imaginary. It is known that $K_{1,2}$ then has CM since the Hodge structure $H^2(K_{1,2},{\mathbb Q}_{K_{1,2}})$ is the tensor product of the CM Hodge structures on $E_1$ and $E_2$ (see \cite{Bor}). In this case, for any cycle $\sigma$ in the submodule of $H_2(K_{1,2},{\mathbb Z})$ generated by the $\gamma_{ijk\ell}$, we have the \emph{unnormalized period}
$$
\int_\sigma\Omega=\alpha(\sigma)\omega_{11}\omega_{21},
$$
where $\alpha(\sigma)$ is an algebraic number in $K={\mathbb Q}(\tau_{12},\tau_{22})$, depending on $\sigma$, and corresponds to the normalized period of $\Omega$. As $K$ is the composite of two CM fields, it is also a CM field, that is a totally imaginary quadratic extension of a totally real field. If $E_1$ is isogenous to $E_2$, then $K$ is imaginary quadratic. In this case, the unnormalized periods are either zero or transcendental, being up to an algebraic factor just $\alpha(\sigma)$ times the square of the transcendental period of the holomorphic 1-form defined over ${\overline{\mathbb Q}}$ on either of the curves. If $E_1$ is not isogenous to $E_2$, the transcendence of this number is not known.

\noindent

\noindent Returning to the general situation, let $L$ be the ${\mathbb Q}$-algebra of ${\mathbb Q}$-linear endomorphisms of $T$ whose extension to $T_{\mathbb C}$ preserves the $T_{\mathbb C}\cap H^{p,q}$, $p+q=2$. Then $L$ is the ${\mathbb Q}$-algebra of ${\mathbb Q}$-linear endomorphisms of $T$ commuting with the Mumford-Tate group of $T$. As $L$ is faithfully represented by its action on the 1-dimensional space $H^{2,0}$, it is a number field. It can be shown that $L$ is either totally real or a CM field. For example, when $\tau_{12}$ and $\tau_{22}$ are quadratic imaginary, we saw above that $H^{2,0}\cap T_K\simeq K$, a CM field. Therefore $L\simeq K$, and the Mumford-Tate group is contained in the algebraic torus ${\mathbb T}_K$ defined over ${\mathbb Q}$ and satisfying ${\mathbb T}_K(A)=(A\otimes_{\mathbb Q}K)^\ast$ for any ring $A$. Therefore, the Mumford-Tate group is abelian. For more details see \cite{Huy}, \cite{Zar}.

\noindent As another example of a period of $\Omega$, consider the 2-cycle $\sigma$ on $K_{1,2}$ determined by $\gamma_{11}\times\gamma_{22}-\gamma_{12}\times\gamma_{21}$. We have the \emph{unnormalized period}
$$
\int_\sigma\Omega=\frac12\left(\omega_{11}\omega_{22}-\omega_{12}\omega_{21}\right)=\frac12\omega_{11}\omega_{21}\left(\tau_{22}-\tau_{12}\right),
$$
Thus, the difference $\tau_{22}-\tau_{12}$ appears as a normalized period. If $\lambda_1$, $\lambda_2$ are algebraic numbers, and the two elliptic curves do not have CM and are not isogenous, the transcendence of this number is not known (see also \cite{Sh2}).

\bigskip

\item We now consider the construction of \S\ref{s:BVtower}. Suppose that we have two families of Calabi-Yau varieties with involution, both defined over ${\overline{\mathbb Q}}$, and both satisfying Expectation 4. Then, the family of Calabi-Yau varieties constructed from them in the next step of the Borcea-Voisin tower also satisfies Expectation 4. The reasoning is similar to that of Example 1 of this section.

    \noindent We adopt the notation of \S\ref{s:BVtower}. For $i=1,2$, let $(A_i,I_i)$ be a Calabi-Yau variety of dimension $d_i$ with involution $I_i$ fixing a union $R_i$ of smooth distinct hypersurfaces. Assume that $A_i, I_i$ are defined over $\overline{\mathbb Q}$. Let ${\widehat{Y}}$ be the blow-up of $A_1\times A_2$ along $R_1\times R_2$ and let ${\widehat{I_{1,2}}}$ be the involution on ${\widehat{Y}}$ induced by $I_1\times I_2$. Let $(B,I)$ be the Calabi-Yau variety with involution defined over ${\overline{\mathbb Q}}$, where $B={\widehat{Y}}/{\widehat{I_{1,2}}}$.

    \noindent Recall the situation and notations of \S\ref{s:BVtower}. In particular, we have two smooth projective families ${\mathcal A}_1\rightarrow S_1$, ${\mathcal A}_2\rightarrow S_2$ of Calabi-Yau varieties defined over ${\overline{\mathbb Q}}$. Let $A_1(s_1)$ be the fiber of ${\mathcal A}_1$ at $s_1\in S_1$ and $A_2(s_2)$ the fiber of ${\mathcal A}_2$ at $s_2\in S_2$.  Suppose, for $i=1,2$, that $\dim(A_i(s_i))=d_i$ (independent of $s_i\in S_i$). Recall from \S\ref{s:VHS} that, by definition, the variety $A_i(s_i)$ is defined over ${\overline{\mathbb Q}}$ when $s_i\in S_i({\overline{\mathbb Q}})$. Suppose, for $i=1,2$, there is an involution ${\mathcal I}_i$ of ${\mathcal A}_i\rightarrow S_i$, leaving $S_i$ fixed, with ramification locus ${\mathcal R}_i\rightarrow S_i$ a smooth projective family defined over ${\overline{\mathbb Q}}$, whose fiber $R_i(s_i)$ at $s_i\in S_i$ consists of a union of smooth distinct hypersurfaces. For $s=(s_1,s_2)\in S=S_1\times S_2$, let $B(s)$ be the Calabi-Yau variety with involution constructed from the Calabi-Yau varieties with involution $A_1(s_1)$ and $A_2(s_2)$ using the method described before Proposition 2, \S\ref{s:BVtower}, and let ${\mathcal B}\rightarrow S$ be the smooth projective family whose fiber at $s\in S$ is $B(s)$. Let $R(s)$ be the ramification locus of the involution on $B(s)$ and let ${\mathcal R}\rightarrow S$ be the smooth projective family whose fiber at $s\in S$ is $R(s)$. We say that Expectation 4 holds for a smooth projective family ${\mathcal X}\rightarrow S$ defined over ${\overline{\mathbb Q}}$ when the set of fibers ${\mathcal X}_s$, $s\in S({\overline{\mathbb Q}})$, such that the $F^{n,n}$-part of the Hodge filtration of $H^n({\mathcal X}_s,{\mathbb Q}_{{\mathcal X}_s})_{{\rm prim},\mathbb C}$ is defined over ${\overline{\mathbb Q}}$, is the same as the set of $s\in S({\overline{\mathbb Q}})$ such that $T({\mathcal X}_s)$ has CM. Assume the families ${\mathcal A}_1\rightarrow S_1$, ${\mathcal A}_2\rightarrow S_2$ both satisfy Expectation 4. For $i=1,2$ and $s_i\in S_i$, let $\Omega_{d_i}(s_i)$ be a non-vanishing holomorphic $d_i$-form on $A_i(s_i)$. (If we are interested in the possible transcendence of the unnormalized periods of $\Omega_{d_i}(s_i)$, we should also assume it is defined over ${\overline{\mathbb Q}}$ whenever $s_i\in S_i({\overline{\mathbb Q}})$.) Therefore, it follows from our assumptions that, for $i=1,2$, if $s_i\in S_i({\overline{\mathbb Q}})$, and the unnormalized periods $\omega_{ij}(s_i)$ of $\Omega_{d_i}(s_i)$ are all algebraic multiples of each other, then the transcendental part $T(A_i(s_i))$ of $H^{d_i}(A_i(s_i),{\mathbb Q}_{A_i(s_i)})_{\rm prim}$ has CM. Here $j=1,\ldots,m_i$ ranges over a ${\mathbb Z}$-basis of $H_{d_i}(A_i(s_i),{\mathbb Z})$. We can suppose that $\omega_{11}\not=0$ and $\omega_{21}\not=0$, and consider the normalized periods $\tau_{ij}=\tau_{ij}(s_i)=\omega_{ij}/\omega_{i1}$, $j=1,\ldots m_i$, $i=1,2$. Therefore, if the normalized periods of $\Omega_{d_i}(s_i)$ are all algebraic, then $T(A_i(s_i))$ has CM.

\noindent For $i=1,2$, fix $s_i\in S_i({\overline{\mathbb Q}})$, and let $A_i=A_i(s_i)$ and $\Omega_i=\Omega_{d_i}(s_i)$. Let $B=B(s_1,s_2)$ be the Calabi-Yau variety of dimension $d=d_1+d_2$, with involution, constructed from $(A_1(s_1), I_1(s_1)$) and $(A_2(s_2),I_2(s_2))$, at the next step of the Borcea-Voisin tower as in \S\ref{s:BVtower}. Suppose that the filtration $F^{d,d}\subseteq H^d(B,{\mathbb C})$ is defined over ${\overline{\mathbb Q}}$, so that if $\Omega_d$ is a nowhere vanishing holomorphic $d$-form on $B$, then its normalized periods are algebraic. From the argument at the beginning of Proposition 2, \S\ref{s:BVtower}, we have
    $$
    F^{d,d}(B)=H^{d,0}(B)= H^{d,0}({\widehat Y})^{{\widehat{I}}_{1,2}}.$$
    Moreover, from Lemma 2 and Lemma 3, \S7, of the present paper and Lemma 7.2.4 of \cite{Roh},
    $$
    H^{d,0}({\widehat{Y}})=H^{d,0}(A_1\times A_2)=H^{d_1,0}(A_1)\otimes H^{d_2,0}(A_2)=H^{d,0}({\widehat{Y}})^{{\widehat{I}}_{1,2}}.
    $$
    Therefore, if the filtration $F^{d,d}(B)\subseteq H^d(B,{\mathbb C})$ is defined over ${\overline{\mathbb Q}}$, then so is the filtration $F^{d,d}({\widehat{Y}})\subseteq H^d({\widehat{Y}},{\mathbb C})$. We can think of $\Omega_1\Omega_2$ as a nonvanishing holomorphic $d$-form on $A_1\times A_2$ which first lifts to ${\widehat{Y}}$ and then descends via the canonical projection $\pi:{\widehat{Y}}\rightarrow B$ to a nowhere vanishing holomorphic $d$-form $\Omega_d$ on $B$. If $\Omega_1$ and $\Omega_2$ are defined over ${\overline{\mathbb Q}}$, then so is $\Omega_d$. As noted already above, this is not important for the transcendence properties of the normalized periods, but is important for the transcendence properties of the unnormalized periods. From \cite{Roh}, Lemma 7.2.4, for $i=1,2$, the form $\Omega_i$ changes sign under the induced action of $I_i$. Therefore $\Omega_d$ is invariant under $I_1\times I_2$, as required. Let $\gamma_{ij}\in H_{d_i}(A_i,{\mathbb Z})$ generate $H_{d_i}(A_i,{\mathbb Q})$, and let $\omega_{ij}=\int_{\gamma_{ij}}\Omega_{d_i}$, for $j=1,\ldots m_i$, $i=1,2$. The cycles given by $\gamma_{1j}\times\gamma_{2k}$ on $A_1\times A_2$ lift to ${\widehat{Y}}$ and then descend via $\pi$ to give $m_1m_2$ elements $\gamma_{1j2k}$ of the ${\mathbb Z}$-module of $d=d_1+d_2$-cycles on $B$.  As $\pi$ is a 2-1 map, we have the following unnormalized periods of $\Omega_d$,
$$
\int_{\gamma_{1j2k}}\Omega_d=\,\frac12\omega_{1j}\omega_{2k}=\,\frac12\omega_{11}\omega_{21}\tau_{1j}\tau_{2k},\quad j=1,\ldots,m_1,\;\; k=1,\ldots, m_2,
$$
which are also unnormalized periods of $\Omega_1\Omega_2$ on $A_1\times A_2$.
Therefore, the products $\tau_{1j}\tau_{2k}$ appear as normalized periods of $\Omega_d$. As $\tau_{11}=\tau_{21}=1$, these products include all the normalized periods $\tau_{1j}$, $\tau_{2k}$ of $\Omega_1$ and $\Omega_2$. Therefore, if the normalized periods of $\Omega_d$ are all algebraic, then the normalized periods of both $\Omega_1$ and $\Omega_2$ are all algebraic. As we have assumed, for $i=1,2$, that Expectation 4 holds for the family ${\mathcal A}_i\rightarrow S_i$, it follows that $T(A_i)=T(A_i(s_i))$ has CM. By \cite{VZ}, Lemma 8.1, it follows that the tensor product of polarized Hodge structures $T(A_1)\otimes T(A_2)$ is a polarized Hodge structure with CM. For $i=1,2$, let $T_i$ be the underlying ${\mathbb Q}$-vector space of $T(A_i)$. Then
$$
H^{d,0}(B)=H^{d_1,0}(A_1)\otimes_{\mathbb C}H^{d_2,0}(A_2)\subseteq T_{1,{\mathbb C}}\otimes_{\mathbb C}T_{2,{\mathbb C}},
$$
so that by the definition of $T(B)=T(B(s_1,s_2))$, it is a polarized Hodge substructure of $T(A_1)\otimes T(A_2)$. By \cite{DeSS}, rational polarized Hodge structures are semi-simple, so that $T(B)$ appears as a direct summand in the decomposition of $T(A_1)\otimes T(A_2)$ into a direct sum of irreducible polarized Hodge structures. Indeed, we have a direct sum decomposition $T(A_1)\otimes T(A_2)\simeq T(B)\oplus T(B)^\perp$, where $T(B)^\perp$ is the orthogonal complement of $T(B)$ with respect to the polarization on $T(A_1)\otimes T(A_2)$. By Lemma 1, \S\ref{s:lemmas}, it follows that there is a group surjection from the Mumford-Tate group of $T(A_1)\otimes T(A_2)$ to the Mumford-Tate group of $T(B)$. As $T(A_1)\otimes T(A_2)$ has CM, its Mumford-Tate group is abelian. Therefore, the Mumford-Tate group of $T(B)$ is abelian, and $T(B)=T(B(s_1,s_2))$ has CM, as required.

\bigskip

\item Consider the $K3$-surface $F_2^{(4)}=F_2^{(4)}(s)$ and the map from $F_1^{(4)}\times\Sigma^{(4)}$, $F_1^{(4)}=F_1^{(4)}(s)$, to $F_2^{(4)}$ as in the discussion preceding Proposition 3, \S\ref{s:expect1}. We have the affine model of $F_2^{(4)}$ on $x_0\not=0$:
     $$
     X_3^4+X_2^4+X_1(X_1-1)(X_1-s)=0,
     $$
     of $F_1^{(4)}$ on $y_0\not=0$:
     $$
     Y_2^4+Y_1(Y_1-1)(Y_1-s)=0,
     $$
      and of $\Sigma^{(4)}$ on $z_2\not=0$:
      $$
      Z_0^4+Z_1^4=1.
      $$
      On $z_2y_0\not=0$, the map from $F_1^{(4)}\times\Sigma^{(4)}$ to $F_2^{(4)}$ is given by
      $$X_1=Y_1,\quad X_2=Y_2Z_0,\quad X_3=Y_2Z_1.$$
      On $y_0\not=0$, the action of $\zeta_4$ on $F_1^{(4)}$ is $(Y_1, Y_2)\mapsto (Y_1,\zeta_4Y_2)$, and on $z_2\not=0$, the action of $\zeta_4$ on $\Sigma^{(4)}$ is $(Z_0,Z_1)\mapsto(\zeta_4^{-1}Z_0,\zeta_4^{-1}Z_1)$. A nowhere vanishing holomorphic $2$-form on $F_2^{(4)}$ is given by
      $$
     \Omega_2(s)= X_3^{-3}dX_1\wedge dX_2.
      $$
      This $2$-form is defined over ${\overline{\mathbb Q}}$ when $s$ is an algebraic number.
      A direct computation shows that $\Omega_2(s)$ corresponds to
      $$
      \left(Y_2^{-2}dY_1\right)\left(Z_1^{-3}dZ_0\right)
      $$
      in the tensor product of the $\zeta_4^2$-eigenspaces $H^1(F_1^{(4)},{\mathbb C})_2\otimes H^1(\Sigma^{(4)},{\mathbb C})_2$ (see Proposition 3, \S\ref{s:expect1}). The unnormalized period of the $1$-form $Z_1^{-3}dZ_0$ over any $1$-cycle $\gamma$ on $\Sigma^{(4)}$ is of the form
      $$
      \beta(\gamma)\frac{\Gamma(\frac14)\Gamma(\frac14)}{\Gamma(\frac12)},
      $$
      where $\beta(\gamma)$ is an algebraic number, depending on $\gamma$, which corresponds to a normalized period (see \cite{TCY}, Appendix).  The Jacobian of the genus 3 curve $\Sigma^{(4)}$ is isogenous to the third power $E_0^3$ of an elliptic curve with CM by ${\mathbb Q}(\zeta_4)$. As all elliptic curves with CM by ${\mathbb Q}(\zeta_4)$ are isogenous, we can assume that $E_0$ has equation
      $$
y^3=4x^2-4x.
      $$
      The unnormalized periods of $dx/y$ on $E_0$ equal $\Gamma(1/4)^2/\Gamma(1/2)$, up to multiplication by an algebraic number.

      \noindent The Jacobian of $F_1^{(4)}(s)$ is isogenous to $E(s)\times T(s)$, where $E(s)$ is the elliptic curve
      $$
      w^2=u(u-1)(u-s), \qquad s\not=0,1,\infty
      $$
      and $T(s)$ is an abelian variety of dimension 2 whose endomorphism algebra contains ${\mathbb Q}(\zeta_4)$ (see \cite{CoWo1}, \cite {CoWu}, \cite{Wo}). The holomorphic $1$-form $Y_2^{-2}dY_1$ corresponds to the $1$-form $du/w$ on $E(s)$ whose unnormalized periods are
      $$
      \omega_2(s)=\int_1^\infty u^{-1/2}(u-1)^{-1/2}(u-s)^{-1/2}du=\pi F(1/2,1/2,1;s),
      $$
      $$
     \omega_1(s)=\int_{-\infty}^0 u^{-1/2}(u-1)^{-1/2}(u-s)^{-1/2}du=\zeta_4\pi F(1/2,1/2,1;1-s),$$
where $F(a,b,c;x)$, $c\not=0,-1,-2,\ldots$, is the multi-valued classical Gauss hypergeometric function given by the analytic continuation of the power series
$$
F(a,b,c;x)=\sum_{n=0}^\infty\frac{(a)_n(b)_n}{(c)_n}\frac{x^n}{n!},\qquad |x|<1,
$$
where $(w)_n=w(w+1)\ldots(w+n-1)$, $w\in{\mathbb C}$. The monodromy group of $F(1/2,1/2,1;s)$ is the principal congruence subgroup $\Gamma[2]$. When $s$ is an algebraic number $\omega_1(s)$, $\omega_2(s)$ are transcendental numbers \cite{Sch2}. The corresponding normalized period
$$
\tau(s)=\omega_2(s)/\omega_1(s),
$$
of $E(s)$ is the Schwartz triangle map for the hyperbolic triangle with angles $0$, $0$, $0$ at the vertices (\cite{Wo}). By Schneider's theorem \cite{Sch}, if $s$ is algebraic, then $\tau(s)$ is an algebraic number if and only if $E(s)$ has CM.

\noindent Overall, the unnormalized periods of $\Omega_2$ are, up to a non-zero algebraic factor,
$$
P_1(s)=\,\frac{\Gamma(\frac14)\Gamma(\frac14)}{\Gamma(\frac12)}\,\omega_1(s),
$$
and
$$
P_2(s)=\,\frac{\Gamma(\frac14)\Gamma(\frac14)}{\Gamma(\frac12)}\,\omega_2(s).
$$
They satisfy the Picard-Fuchs equation given for the classical Gauss hypergeometric function $F(1/2,1/2,1,s)$ appearing above, and given by
$$
s(1-s)\frac{d^2F}{ds^2}+(1-2s)\frac{dF}{ds}-\frac14F=0.
$$
When $E(s)$ is isogenous to $E_0$, we have $s\in {\overline{\mathbb Q}}$. The periods $P_1(s)$ and $P_2(s)$ are therefore transcendental, being essentially the square of $(\Gamma(1/4))^2/\Gamma(1/2)$. When $E(s)$ and $E_0$ are not isogenous, it is not known whether $P_1(s)$, $P_2(s)$ is transcendental when $s\in {\overline{\mathbb Q}}$. Notice that, by \cite{TCY}, Appendix, the period $\Gamma(1/4)^4/\Gamma(1/2)^2$ occurs as the period of the $(2,0)$-form on the Fermat quartic $x^4+y^4+z^4=1$. This period occurs here when a fiber of the family $F_2^{(4)}(s)$ is isomorphic to the Fermat quartic, in which case $E(s)$ is isogenous to $E_0$.

\noindent On the other hand, we see overall that, up to multiplication by a non-zero algebraic number, the ratio $\tau(s)=\omega_2(s)/\omega_1(s)$ is a normalized period of $\Omega_2$. If $s$ is algebraic, and if the Hodge filtration of $H^2(F_2^{(4)},{\mathbb Q}_{F_2^{(4)}})_{\mathbb C}$ is defined over ${\overline{\mathbb Q}}$, then $\tau(s)$ must be algebraic, so that $E(s)$ has CM.

\noindent We now look at the other factor of the Jacobian of $F_1^{(4)}=F_1^{(4)}(s)$. As $T(s)$ is of dimension 2, the vector space of $(1,0)$-forms has dimension 2, and by \cite{CW}, consists soley of the $\zeta_4$-eigenspace for the induced action of $\zeta_4$ on the holomorphic $1$-forms of $T(s)$. A basis of $H^{1,0}$, in terms of the affine coordinates $(Y_1,Y_2)$ on $F_1^{(4)}$, is given by
$$
Y_2^{-3}(Y_1-s)dY_1,\quad (Y_2^{-3})Y_1dY_1,
$$
which have unnormalized periods (up to a non-zero algebraic factor)
$$
\omega_{11}(s)=\int_1^\infty Y_2^{-3}(Y_1-s)dY_1=\;\frac{\Gamma(\frac14)\Gamma(\frac14)}{\Gamma(\frac12)}\;F\left(\frac14,\frac34,\frac12;s\right),
$$
$$
\omega_{12}(s)=\;\int_0^s Y_2^{-3}(Y_1-s)dY_1=\; \frac{\Gamma(\frac54)\Gamma(\frac14)}{\Gamma(\frac32)}\;s^{1/2}F\left(\frac54,\frac34,\frac32;s\right).
$$
$$
=\;2\frac{\Gamma(\frac14)\Gamma(\frac14)}{\Gamma(\frac12)}\;s^{1/2}F\left(\frac54,\frac34,\frac32;s\right)
$$
$$
\omega_{21}(s)=\;\int_1^\infty(Y_2^{-3})Y_1dY_1=\;\frac{\Gamma(\frac14)}{\Gamma(\frac14)\Gamma(\frac12)}\;F\left(\frac14,\frac34,\frac12;s\right),
$$
$$
\omega_{22}(s)=\;\int_0^s (Y_2^{-3})Y_1dY_1=\;\frac{\Gamma(\frac54)\Gamma(\frac14)}{\Gamma(\frac32)}\;s^{1/2} F\left(\frac54,\frac34,\frac32;s\right)
$$
$$
=\;2 \frac{\Gamma(\frac14)\Gamma(\frac14)}{\Gamma(\frac12)}\;s^{1/2} F\left(\frac54,\frac34,\frac32;s\right).
$$
Here, we have used the functional equation for the $\Gamma$-function given by $\Gamma(z+1)=z\Gamma(z)$, $z\not=0,-1,-2,\ldots$. The monodromy group of all the above (contiguous) hypergeometric functions is the spherical triangle group with signature $(2,2,2)$, which is the dihedral group with $4$ elements. As this group is finite, the hypergeometric functions are algebraic over ${\overline{\mathbb Q}}(s)$, and therefore take algebraic values when $s\in{\overline{\mathbb Q}}$. Indeed, it is well-known that 
$$
F\left(\frac14,\frac34,\frac12;s\right)=\frac{(1-\sqrt s)^{-1/2}+(1+\sqrt s)^{-1/2}}{2}
$$
and
$$
F\left(\frac54,\frac34,\frac32;s\right)=\frac{(1-\sqrt s)^{-1/2}-(1+\sqrt s)^{-1/2}}{s^{\frac12}}\;\;.
$$
The corresponding normalized period is given by the following Schwarz triangle map for the spherical triangle with vertices $(\pi/2,\pi/2,\pi/2)$,
$$
{\mathcal T}(s)=\frac{\omega_{12}(s)}{\omega_{11}(s)}=\frac{\omega_{22}(s)}{\omega_{21}(s)}=s^{1/2}\,\frac{F\left(\frac54,\frac34,\frac32;s\right)}{F\left(\frac14,\frac34,\frac12;s\right)}
$$
$$
=\,2\,\frac{(1-\sqrt s)^{-1/2}-(1+\sqrt s)^{-1/2}}{(1-\sqrt s)^{-1/2}+(1+\sqrt s)^{-1/2}}\;.
$$
The unnormalized periods $\omega_{ij}$, $i,j=1,2$ are all transcendental when $s\in{\overline{\mathbb Q}}$, as follows from the transcendence of  $B(\frac14,\frac14):=\Gamma(\frac14)\Gamma(\frac14)/\Gamma(\frac12)$ which is, as already remarked, a period of a differential form of the first kind defined over ${\overline{\mathbb Q}}$ on the Fermat curve $\Sigma^{(4)}$, up to multiplication by a non-zero algebraic number. The normalized period ${\mathcal T}(s)$ is algebraic for all $s\in{\overline{\mathbb Q}}$. By \cite{Co1},\cite{SW}, it follows that the abelian surface $T(s)$ has CM, but we can also argue as follows. If ${\mathcal P}(s)$ is the ${\overline{\mathbb Q}}$-vector space generated by the numbers $\omega_{ij}(s)$, $i,j=1,2$, then $\dim_{\overline{\mathbb Q}}({\mathcal P})=1$ for all $s\in{\overline{\mathbb Q}}$. By \cite{Wu3}, if $T(s)$ is simple and $s\in{\overline{\mathbb Q}}$, then
$$
\dim_{\overline{\mathbb Q}}({\mathcal P})[L(s):{\mathbb Q}]=2(\dim_{\mathbb C}(T(s)))^2,
$$
where $L(s)={\rm End}_0(T(s))$ is the endomorphism algebra of $T(s)$. But this cannot be true as $[L(s):{\mathbb Q}]$ is at most $2\dim(T(s))=4$. Therefore $T(s)$ is not simple, and must be isogenous to a product of two elliptic curves $E_1\times E_2$. If $E_1$ and $E_2$ are not isogenous, this again contradicts $\dim_{\overline{\mathbb Q}}({\mathcal P})=1$. Thus $E_1$ is isogenous to $E_2$ and both elliptic curves have endomorphism algebra ${\mathbb Q}(\zeta_4)$. Since any two elliptic curves with endomorphism algebra ${\mathbb Q}(\zeta_4)$ are isogenous, we can suppose $T(s)$ is isogenous to $E_0^2$. In fact $T(s)$ is isogenous to $E_0^2$ for all $s$, algebraic or not, by \cite{Shi}, Proposition 14.

\noindent Therefore, for $s$ algebraic, if the normalized periods of $\Omega_2(s)$ on $F_2^{(4)}(s)$ are algebraic, then the curve $F^{(4)}_1(s)$ has CM. This happens if and only if the normalized period $\tau(s)$ on $E(s)$ is algebraic. By \cite{VZ}, if $F^{(4)}_1(s)$ has CM, then $F^{(4)}_2(s)$ has CM. It follows that Expectation 3 holds for the family ${\mathcal F}_2(4,1)$ whose fiber at $s\in {\mathbb C}\setminus\{0,1\}$ is $F_2^{(4)}(s)$.

\item References for this example are \cite{CW}, \cite{CoWo2}, \cite{DM}, \cite{DTT}, \cite{KR}, \cite{TM}. Consider the quintic 3-fold $F_3^{(5)}=F_3^{(5)}(s)=F_3^{(5)}(a_1,a_2)$ and the map from $F_2^{(5)}\times\Sigma^{(5)}$ to $F_3^{(5)}$, where $F_2^{(5)}=F_2^{(5)}(s)$. In affine coordinates, the equation of $F_3^{(5)}$ is:
     $$
     X_4^5+X_3^5+X_2^5+X_1(X_1-1)(X_1-a_1)(X_1-a_2)=0,
     $$
     of $F_2^{(5)}$ is:
     $$
    Y_3^5+Y_2^5+Y_1(Y_1-1)(Y_1-a_1)(Y_1-a_2)=0,
     $$
      and of $\Sigma^{(5)}$ is:
      $$
      Z_0^5+Z_1^5=1.
      $$
      The map from $F_2^{(5)}\times\Sigma^{(5)}$ to $F_3^{(5)}$ is given by
      $$X_1=Y_1,\quad X_2=Y_2,\quad X_3=Y_3Z_0,\quad X_4=Y_3Z_1.$$
      The action of $\zeta_5$ on $F_2^{(5)}$ is $(Y_1,Y_2,Y_3)\mapsto (Y_1,Y_2,\zeta_5Y_3)$, and on $\Sigma^{(5)}$ is $(Z_0,Z_1)\mapsto(\zeta_5^{-1}Z_0,\zeta_5^{-1}Z_1)$. A nowhere vanishing holomorphic $3$-form on $F_3^{(5)}$ is given by
      $$
    \Omega_3=\Omega_3(a_1,a_2)= X_4^{-4}dX_1\wedge dX_2\wedge dX_3.
      $$
      This $3$-form is defined over ${\overline{\mathbb Q}}$ when $a_1,a_2$ are algebraic numbers.
      A direct computation shows that $\Omega_3(a_1,a_2)$ corresponds to
      $$
      \left(Y_3^{-3}dY_1\wedge dY_2\right)\left(Z_1^{-4}dZ_0\right)
      $$
      in the tensor product $H^1(F_1^{(5)},{\mathbb C})_2\otimes H^1(\Sigma^{(5)},{\mathbb C})_3$ (see Proposition 3, \S\ref{s:expect1}). Now consider the map from $F_1^{(5)}\times\Sigma^{(5)}$, $F_1^{(5)}=F_1^{(5)}(a_1,a_2)$, to $F_2^{(5)}$. In a similar manner to the previous example, let $(W_1,W_2)$ be affine coordinates on $F_1^{(5)}$ with $\zeta_5$ acting by $(W_1,W_2)\rightarrow (W_1,\zeta_5 W_2)$. Then the map is given by $Y_1=W_1$, $Y_2=W_2Z_0$, $Y_3=W_2Z_1$. We can then write $\Omega_3$ as
      $$
      \Omega_3=\left(W_2^{-2}dW_1\right)\left(Z_1^{-3}dZ_0\right)\left(Z_1^{-4}dZ_0\right)
      $$
      corresponding to an element of
      $$
      H^1(F_1^{(5)},{\mathbb C})_3\otimes H^1(\Sigma^{(5)},{\mathbb C})_2\otimes H^1(\Sigma^{(5)},{\mathbb C})_3.
      $$

      The Jacobian variety $J(a_1,a_2)$ of $F_1^{(5)}(a_1,a_2)$ has dimension 6. The Jacobian variety of $\Sigma^{(5)}$ also has dimension 6 and is isogenous to the third power $A_0^3$ of a simple abelian surface $A_0$ with CM by ${\mathbb Q}(\zeta_5)$. The periods of the differential forms of the first kind on $A_0$ generate over ${\overline{\mathbb Q}}$ a vector space ${\mathcal P}$ of dimension 2 with basis $B(1/5,1/5)$ and $B(2/5,2/5)$, where $B(a,b)=\Gamma(a)\Gamma(b)/\Gamma(a+b)$. Note that the action of $\zeta_5$ on $\Sigma^{(5)}$ does not induce the action $(u,w)\mapsto (u,\zeta_5^{-1}w)$ on the image curve $w^5=u^{5-r}(1-u)^{5-s}$, $1\le r,s$, $r+s\le 4$, considered in the references given at the beginning of this example. Here $u=Z_0^5, w=Z_0^{5-r}Z_1^{5-s}$, and the induced action of $\zeta_5$ on the holomorphic differential form $du/w$ is multiplication by $\zeta_5^{-r-s}$. We have,
       $$
       \int_0^1\frac{du}w=\,\int_0^1u^{\frac{r}{5}-1}(1-u)^{\frac{s}{5}-1}du=\,B\left(\frac{r}{5},\frac{s}{5}\right)=
       \,\frac{\Gamma\left(\frac{r}{5}\right)\Gamma\left(\frac{s}{5}\right)}{\Gamma\left(\frac{r+s}{5}\right)}.
       $$
       The holomorphic $1$-forms defined over ${\overline{\mathbb Q}}$ on $\Sigma^{(5)}$ which are in the $\zeta_5^3$-eigenspace $H^1(\Sigma^{(5)},{\mathbb C})_3$ have periods in the $1$-dimensional subspace of ${\mathcal P}$ generated by $B(1/5,1/5)$. The holomorphic $1$-forms defined over ${\overline{\mathbb Q}}$ on $\Sigma^{(5)}$ which are in the $\zeta_5^2$-eigenspace have periods in the $1$-dimensional subspace of ${\mathcal P}$ generated by $B(2/5,2/5)$. Here, we have used the fact that $\Gamma(z)\Gamma(1-z)=\pi/\sin(\pi z)$, so that
       $$\frac{\Gamma\left(\frac{1}{5}\right)}{\Gamma\left(\frac{3}{5}\right)}=
       \alpha \frac{\pi}{\Gamma\left(\frac{4}{5}\right)}
       \frac{\Gamma\left(\frac{2}{5}\right)}{\pi}=\,\alpha\frac{\Gamma\left(\frac{2}{5}\right)}{\Gamma\left(\frac{4}{5}\right)},
       $$
       for a non-zero algebraic number $\alpha$, so that $B(1/5,2/5)=\alpha B(2/5,2/5)$.

       \noindent For an integer $n$ with $1\le n\le4$, let $r_n$ be the dimension of the subspace of holomorphic $1$-forms on $F_1^{(5)}$ in the $\zeta_5^n$-eigenspace $(H^1(F_1^{(5)},{\mathbb C}))_n$. Then $r_1=3$, $r_2=2$, $r_3=1$, $r_4=0$. (Our $r_n$ is the $r_{-n}$ of the references at the beginning of this example). Notice that $r_i+r_{5-i}=3$. By \cite{CoWo2}, Th\'eor\`eme 3, the periods of the $\zeta_5^3$-eigenform $W^{-2}_2dW_1$ are given, up to multiplication by non-zero algebraic numbers, by
       $$
       P_1(a_1,a_2)=\,\int_0^1W_1^{-2/5}(W_1-1)^{-2/5}(W_1-a_1)^{-2/5}(W_1-a_2)^{-2/5}dW_1$$
       $$
       P_2(a_1,a_2)=\,\int_1^{a_1}W_1^{-2/5}(W_1-1)^{-2/5}(W_1-a_1)^{-2/5}(W_1-a_2)^{-2/5}dW_1
       $$
       $$
       P_3(a_1,a_2)=\,\int_0^{a_2}W_1^{-2/5}(W_1-1)^{-2/5}(W_1-a_1)^{-2/5}(W_1-a_2)^{-2/5}dW_1,
       $$
     and, in the neighborhood of $a_1=1, a_2=0$, can be expressed explicitly in terms of Appell hypergeometric series. The periods are all solutions of the system of Appell partial differential equations satisfied by
       $$
       F_1(a,b,b',c;a_1,a_2)=\sum_{m,n=0}^\infty\frac{(a)_{m+n}(b)_m(b')_n}{(c)_{m+n}}\frac{a_1^m}{m!}\frac{a_2^n}{n!}, \qquad |a_1|,|a_2|<1
       $$
       where $a=3/5$, $b=2/5$, $b'=2/5$, $c=6/5$. This system of partial differential equations is given by
       $$
       a_1(1-a_1)\frac{\partial^2F}{\partial a_1^2}+a_2(1-a_1)\frac{\partial^2F}{\partial a_1\partial a_2}+\left(\frac65-2a_1\right)\frac{\partial F}{\partial a_1}-\frac25a_2\frac{\partial F}{\partial a_2}-\frac6{25}F=0,
       $$
       $$
       a_2(1-a_2)\frac{\partial^2F}{\partial a_2^2}+a_1(1-a_2)\frac{\partial^2F}{\partial a_1\partial a_2}+\left(\frac65-2a_2\right)\frac{\partial F}{\partial a_2}-\frac25a_1\frac{\partial F}{\partial a_1}-\frac6{25}F=0.
       $$
       The solutions are multi-valued functions of $a_1$, $a_2$, spanning a complex vector space of dimension 3, with monodromy group an infinite arithmetic group acting discontinuously on the complex 2-ball, and corresponding to the signature $(\mu_i)_{i=0}^4=(2/5,2/5,2/5,2/5,2/5)$ in the list of Deligne-Mostow \cite{DM}, see also \cite{CoWo2}. As a function of $(a_1,a_2)$, the period $P_1(a_1,a_2)$ is holomorphic in a neighborhood of $(1,0)$ and $P_1(1,0)=B(1/5,1/5)$. Neither of the periods $P_2(a_1,a_2)$, $P_3(a_1,a_2)$ is holomorphic in a neighborhood of $(1,0)$, but their values at this point are well-defined and both equal zero. For more information on Appell hypergeometric series, see \cite{AK}, and \cite{Yos}, Chapter 6. We see therefore that, up to algebraic factors, the unnormalized periods of $\Omega_3(a_1,a_2)$ are given by
       $$
       B(1/5,1/5)B(2/5,2/5)P_i(a_1,a_2),\qquad i=1,2,3
       $$
and nothing is known in general about their transcendence when $a_1,a_2$ are algebraic numbers. The corresponding normalized periods are, if we divide by the period holomorphic in a neighborhood of $(1,0)$,
$$
P_2(a_1,a_2)/P_1(a_1,a_2),\quad P_3(a_1,a_2)/P_1(a_1,a_2),
$$
and if these take algebraic values, then $F_1^{(5)}(a_1,a_2)$ has CM by \cite{Co1},\cite{SW} (see also \cite{SSW}, \S1.5). Indeed, these normalized periods correspond to points in the Shimura period domain for the analytic family associated to the signature $(2/5,2/5,2/5,2/5,2/5)$, as in \cite{CoWo2}, which in this case is the complex 2-ball. By \cite{VZ}, it then follows that $F_3^{(5)}$ has CM. Note that, to deduce the generalization of Schneider's Theorem, it suffices to consider the periods of $\Omega_3$, so that Expectation 3 holds for the family ${\mathcal F}_3(5,2)$ whose fiber at $(a_1,a_2)\in {\mathbb C}^2$, $a_1\not=a_2$, $a_1,a_2\not=0,1$ is $F_3^{(5)}(a_1,a_2)$.

\noindent Notice that, if $(a_1,a_2)$ is such that $F_1^{(5)}$ is isomorphic to $\Sigma^{(5)}$, then the periods $P_i(a_1,a_2)$, $i=1,2,3$ are all of the form $B(1/5,1/5)$ times an algebraic number, since we are working with a period on $F_1^{(5)}$ in the $\zeta_5^3$-eigenspace. The unnormalized periods of $\Omega_3$ are then, up to an algebraic factor, all equal to
$$
B(1/5,1/5)B(2/5,2/5)B(1/5,1/5)=\left(\Gamma(1/5)\right)^4\left(\Gamma(4/5)\right)^{-1}.
$$
This is the unnormalized period on the Fermat quintic three-fold by \cite{TCY}, Appendix, and as noted there its transcendence is unproven. It is well known that there is a rational map from $\Sigma^{(5)}\times \Sigma^{(5)}\times \Sigma^{(5)}$ to this Fermat quintic, which is a special case of the discussion preceding Proposition 3, see for example \cite{HuS}.

\end{enumerate}

\bigskip

\section{Hodge theory and linear algebra lemmas}\label{s:lemmas}

In this section we collect, for the convenience of the reader, the main lemmas that we use in this paper. We begin by stating some well-known lemmas from Hodge theory, which, for the most part,  predate the references we give for them.

\bigskip

\noindent{\bf Lemma 1:} \cite{Bor},\cite{VZ} (i) Let $(V_1,\varphi_1)$ and $(V_2,\varphi_2)$ be two polarized rational Hodge structures of weight $n$ and $\varphi_1\oplus\varphi_2$ the induced Hodge structure on $V_1\oplus V_2$. Then,
$$
M_{\varphi_1\oplus\varphi_2}\subset M_{\varphi_1}\times M_{\varphi_2}\subset {\rm{GL}}(V_1)\times{\rm{GL}}(V_2)\subset {\rm{GL}}(V_1\oplus V_2),
$$
and the projections
$$
M_{\varphi_1\oplus\varphi_2}\rightarrow M_{\varphi_1},\qquad M_{\varphi_1\oplus\varphi_2}\rightarrow M_{\varphi_2}
$$
are surjective.

\medskip

(ii) The Mumford-Tate group does not change under Tate twists (that is, shifts in bi-degree $(1,1)$).

\medskip

(iii) The Mumford-Tate group of a Hodge structure concentrated in bi-degree $(p,p)$, $p\in{\mathbb Z}$, is trivial.

\medskip

(iv) Let $\varphi_1\otimes\varphi_2$ be the induced Hodge structure on $V_1\otimes V_2$. Then $\varphi_1\otimes\varphi_2$ has CM if and only if both $\varphi_1$ and $\varphi_2$ have CM.

\bigskip

\noindent {\bf Lemma 2:} \cite{Cat} Let $X_1$ and $X_2$ be compact K\"ahler manifolds. Then, for any integers $k,r,s\ge0$, we have
$$
H^k(X_1\times X_2,{\mathbb Q})=\oplus_{i+j=k}H^i(X_1,{\mathbb Q})\otimes H^j(X_2,{\mathbb Q})
$$
and
$$
H^{r,s}(X_1\times X_2)=\oplus_{p+p'=r,\, q+q'=s}H^{p,q}(X_1)\otimes H^{p',q'}(X_2).
$$

\bigskip

\noindent {\bf Lemma 3:} \cite{Roh}, \cite{Voi}. Let $X$ be an algebraic manifold of dimension $n$ and let $\widehat{X}$ be the blow-up of $X$ along a submanifold $Z$ of codimension 2 in $X$. Then, for all $k$, we have an isomorphism of Hodge structures
$$
H^k(X,{\mathbb Q}_X)\oplus H^{k-2}(Z,{\mathbb Q}_Z)(-1)\simeq H^k(\widehat{X},{\mathbb Q}_{\widehat{X}}),
$$
where $H^{k-2}(Z,{\mathbb Q}_Z)(-1)$ is $H^{k-2}(Z,{\mathbb Q}_Z)$ shifted by $(1,1)$ in bidegree. Thus, the Mumford-Tate group of $H^k(\widehat{X},{\mathbb Q}_{\widehat{X}})$ is abelian if and only if the Mumford-Tate groups of both $H^k(X,{\mathbb Q}_X)$ and $H^{k-2}(Z,{\mathbb Q}_Z)$ are abelian. What's more, if $X$ is a smooth surface and $Z$ is a point of $X$, then the Mumford--Tate groups of $H^2(X,{\mathbb Q}_X)$ and $H^2(\widehat{X},{\mathbb Q}_{\widehat{X}})$ are isomorphic.

\medskip

\noindent{\bf Lemma 4:} Let $(V,{\widetilde{\varphi}})$ be a rational Hodge structure of weight $n$ and $F^{\ast,n}_{\mathbb C}$ the associated Hodge filtration of $V_{\mathbb C}$. Then we say that $F^{\ast,n}_{\mathbb C}$ is defined over ${\overline{\mathbb Q}}$ if it is induced, by extension of scalars to ${\mathbb C}$, from a ${\overline{\mathbb Q}}$-filtration of $V_{\overline{\mathbb Q}}$. Let $(V_1,{\widetilde{\varphi}}_1)$ and $(V_2,{\widetilde{\varphi}}_2)$ be two rational Hodge structures of weight $n$ and ${\widetilde{\varphi}}={\widetilde{\varphi}}_1\oplus{\widetilde{\varphi}}_2$ the induced Hodge structure on $V=V_1\oplus V_2$. If the Hodge filtration $F^{\ast,n}$ of $V_{\mathbb C}$ associated to $(V,{\widetilde{\varphi}})$ is defined over ${\overline{\mathbb Q}}$, then, for $i=1,2$, so is the Hodge filtration $F_i^{\ast,n}$ of $V_{i,{\mathbb C}}$ associated to ${\widetilde{\varphi}}_i$.

\medskip

\noindent {\bf Proof:} Let $F^{\ast,n}_{\mathbb C}$ be the Hodge filtration associated to $(V,{\widetilde{\varphi}})$, and for $i=1,2$, let $F^{\ast,n}_{i,{\mathbb C}}$ be the Hodge filtration associated to $(V_i,{\widetilde{\varphi}}_i)$. Then, as ${\widetilde{\varphi}}={\widetilde{\varphi}}_1\oplus{\widetilde{\varphi}}_2$, we have, for all $p,q$ with $p+q=n$,
$$
V_{\mathbb C}^{p,q}=V_{1,{\mathbb C}}^{p,q}\oplus V_{2,{\mathbb C}}^{p,q}.
$$
Therefore, for $p=0,\ldots,n$,
$$
F^{p,n}=F_1^{p,n}\oplus F_2^{p,n}.
$$
By assumption, there is a ${\overline{\mathbb Q}}$-vector subspace $F^{p,n}_{\overline{\mathbb Q}}$ of $V_{\overline{\mathbb Q}}$ such that we have $F^{p,n}=F^{p,n}_{\overline{\mathbb Q}}\otimes_{\overline{\mathbb Q}}{\mathbb C}$. Let $e_1,\ldots, e_f$ be a basis of $F^{p,n}_{\overline{\mathbb Q}}$, where $f=\dim_{\mathbb C}(F^{p,n})$. As
$$
V_{\overline{\mathbb Q}}=V_{1,{\overline{\mathbb Q}}}\oplus V_{2,{\overline{\mathbb Q}}},
$$
we have $e_j=E_{1,j}\oplus E_{2,j}$, for some $E_{1,j}\in V_{1,{\overline{\mathbb Q}}}$ and $E_{2,j}\in V_{2,{\overline{\mathbb Q}}}$, $j=1,\ldots,f$, which are unique. Therefore
$$
e_j\otimes_{\overline{\mathbb Q}}1_{\mathbb C}=(E_{1,j}\otimes_{\overline{\mathbb Q}}1_{\mathbb C})\oplus (E_{2,j}\otimes_{\overline{\mathbb Q}}1_{\mathbb C})
$$
is the unique decomposition of $e_j\otimes_{\overline{\mathbb Q}}1_{\mathbb C}\in F^{p,n}$ as a direct sum of an element $E_{1,j}\otimes_{\overline{\mathbb Q}}1_{\mathbb C}\in F_1^{p,n}$ and an element $E_{2,j}\otimes_{\overline{\mathbb Q}}1_{\mathbb C}\in F_2^{p,n}$.
Let $\omega_1\in F_1^{p,n}$. Then $\omega_1\oplus 0$ is in $F^{p,n}$, so that there are complex numbers $\lambda_j$, $j=1,\ldots,f$, such that
$$
\omega_1\oplus 0=\sum_{j=1}^f\lambda_j(e_j\otimes_{\overline{\mathbb Q}}1_{\mathbb C})=\left(\sum_{j=1}^f\lambda_j(E_{1,j}\otimes_{\overline{\mathbb Q}}1_{\mathbb C})\right)\oplus\left(\sum_{j=1}^f\lambda_j(E_{2,j}\otimes_{\overline{\mathbb Q}}1_{\mathbb C})\right).
$$
Therefore $\sum_{j=1}^f\lambda_j(E_{2,j}\otimes_{\overline{\mathbb Q}}1_{\mathbb C})=0$ and $\omega_1=\sum_{j=1}^f\lambda_j(E_{1,j}\otimes_{\overline{\mathbb Q}}1_{\mathbb C})$. Therefore $F^{p,n}_1$ is spanned by the $(E_{1,j}\otimes_{\overline{\mathbb Q}}1_{\mathbb C})\in F^{p,n}_1\cap \{v\otimes_{\overline{\mathbb Q}}1_{\mathbb C}:v\in V_{1,{\overline{\mathbb Q}}}\}$, $j=1,\ldots,f$. Choosing a maximal ${\overline{\mathbb Q}}$-linearly independent subset of this spanning set gives a basis of $F^{p,n}_1$ in $F^{p,n}_1\cap \{v\otimes_{\overline{\mathbb Q}}1_{\mathbb C}:v\in V_{1,{\overline{\mathbb Q}}}\}$. As this works for any $p$, it follows that the Hodge filtration $F^{\ast,n}_1$ is defined over ${\overline{\mathbb Q}}$. The analogous reasoning shows that $F^{\ast,n}_2$ is defined over ${\overline{\mathbb Q}}$.

\bigskip

\noindent The remaining lemmas of this section are, for the most part, straightforward applications of elementary multilinear algebra. They are nonetheless crucial for the proofs of this paper, and we include them for the reader's convenience. For an ``online'' reference which thoroughly treats the fundamentals of tensor products and bilinear forms, see \cite{CoK}.

\bigskip

\noindent{\bf Lemma 5:} Let $V_{1,{\overline{\mathbb Q}}}$, $V_{2,{\overline{\mathbb Q}}}$ be two ${\overline{\mathbb Q}}$-vector spaces and let $W_{\overline{\mathbb Q}}=V_{1,{\overline{\mathbb Q}}}\otimes_{\overline{\mathbb Q}}V_{2,{\overline{\mathbb Q}}}$. Suppose $\Omega\in W_{\overline{\mathbb Q}}$ is such that $\Omega\otimes_{\overline{\mathbb Q}}1_{\mathbb C}=\Omega_1\otimes_{\mathbb C}\Omega_2$ with $\Omega_i\in V_{i,{\mathbb C}}$, $\Omega_i\not=0$, $i=1,2$. Then, for $i=1,2$, there is a non-zero complex scalar $c_i$ such that $c_i\Omega_i\in \{v_i\otimes_{\overline{\mathbb Q}}1_{\mathbb C}: v_i\in V_{i,{\overline{\mathbb Q}}}\}$.

\medskip

{\noindent Proof:} Let $g_i=\dim(V_{i,{\overline{\mathbb Q}}})$, $i=1,2$. Let $e_1,\ldots, e_{g_1}$ be a basis of $V_{1,{\overline{\mathbb Q}}}$ and $f_1,\ldots,f_{g_2}$ be a basis of $V_{2,{\overline{\mathbb Q}}}$. As $\Omega_1\not=0$, $\Omega_2\not=0$, by reordering the bases if necessary, we can assume that
$$
\Omega_1=\sum_{i=1}^{g_1}\alpha_i(e_i\otimes_{\overline{\mathbb Q}}1_{\mathbb C}),\quad \Omega_2=\sum_{j=1}^{g_2}\beta_j(f_j\otimes_{\overline{\mathbb Q}}1_{\mathbb C}),
$$
where $\alpha_i, \beta_j\in{\mathbb C}$ and $\alpha_1\not=0$, $\beta_1\not=0$. As  $\Omega\in W_{\overline{\mathbb Q}}$, it has an expansion in terms of the basis $e_i\otimes_{\overline{\mathbb Q}}f_j$, $i=1,\ldots,g_1$, $j=1,\ldots,g_2$, of $W_{\overline{\mathbb Q}}$ of the form
$$
\Omega=\sum_{i,j}a_{ij}e_i\otimes_{\overline{\mathbb Q}}f_j,\qquad a_{ij}\in{\overline{\mathbb Q}}.
$$
We also have
$$
\Omega\otimes_{\overline{\mathbb Q}}1_{\mathbb C}=\Omega_1\otimes_{\mathbb C}\Omega_2
$$
$$
=\sum_{ij}\alpha_i\beta_j((e_i\otimes_{\overline{\mathbb Q}}f_j)\otimes_{\overline{\mathbb Q}}1_{\mathbb C}).
$$
Comparing coefficients of the $e_i\otimes_{\overline{\mathbb Q}}f_j$ in the expression for $\Omega$ with coefficients of the $(e_i\otimes_{\overline{\mathbb Q}}f_j)\otimes_{\overline{\mathbb Q}}1_{\mathbb C}$ in the expansion for $\Omega\otimes_{\overline{\mathbb Q}}1_{\mathbb C}$ we have $\alpha_i\beta_j=a_{ij}\in{\overline{\mathbb Q}}$. In particular $\alpha_1\beta_1=a_{11}\not=0$. For any $j=1,\ldots,g_2$ we have $\alpha_1\beta_j=a_{1j}$. Therefore $(\alpha_1\beta_j)/(\alpha_1\beta_1)=\beta_j/\beta_1=a_{1j}/a_{11}$, so that $\beta_j=\beta_1(a_{1j}/a_{11})$ and, as $a_{1j}/a_{11}$ is in ${\overline{\mathbb Q}}$, we can take $c_2=\beta_1^{-1}$ and deduce that
$$
c_2\Omega_2\in \{v_2\otimes_{\overline{\mathbb Q}}1_{\mathbb C}: v_2\in V_{2,{\overline{\mathbb Q}}}\}.
$$
A similar argument for $\Omega_1$ concludes the proof of the lemma.

\bigskip

\noindent{\bf Lemma 6:} Let $V_{1,{\overline{\mathbb Q}}}$, $V_{2,{\overline{\mathbb Q}}}$ be two ${\overline{\mathbb Q}}$-vector spaces and let $W_{\overline{\mathbb Q}}=V_{1,{\overline{\mathbb Q}}}\otimes_{\overline{\mathbb Q}}V_{2,{\overline{\mathbb Q}}}$. Let $U_1$ be a ${\mathbb C}$-vector subspace of $V_{1,\mathbb C}=V_{1,{\overline{\mathbb Q}}}\otimes_{\overline{\mathbb Q}}{\mathbb C}$ and $\Omega_2\in V_{2,{\overline{\mathbb Q}}}$, with $\Omega_2\not=0$. If $W_1=U_1\otimes_{\mathbb C}{\mathbb C}(\Omega_2\otimes_{\overline{\mathbb Q}}1_{\mathbb C})$ has a basis in $W_1\cap\{w\otimes_{\overline{\mathbb Q}}1_{\mathbb C}: w\in W_{\overline{\mathbb Q}}\}$, then $U_1$ has a basis in $U_1\cap\{u\otimes_{\overline{\mathbb Q}}1_{\mathbb C}: u\in V_{1,{\overline{\mathbb Q}}}\}$.

\medskip

\noindent {\bf Proof:}  Let $g_i=\dim_{\overline{\mathbb Q}}(V_{i,{\overline{\mathbb Q}}})$, $i=1,2$. Let $\{e_i\}_{i=1}^{g_1}$ be a basis of $V_{1,{\overline{\mathbb Q}}}$ and $\{f_j\}_{j=1}^{g_2}$ be a basis of $V_{2,{\overline{\mathbb Q}}}$. Then the $e_i\otimes_{\overline{\mathbb Q}}f_j$, $i=1,\ldots,g_1$, $j=1,\ldots,g_2$, form a basis ${\mathcal B}$ of $W_{\overline{\mathbb Q}}$. We can assume that
$$
\Omega_2=f_1+b_2f_2+\ldots+b_{g_2}f_{g_2},\qquad b_j\in{\overline{\mathbb Q}},\;j=2,\ldots,g_2.
$$
Let $w_1\in W_1\cap\{w\otimes_{\overline{\mathbb Q}}1_{\mathbb C}: w\in W_{\overline{\mathbb Q}}\}$. As all tensors in $W_1=U_1\otimes_{\mathbb C}{\mathbb C}(\Omega_2\otimes_{\overline{\mathbb Q}}1_{\mathbb C})$ are elementary, there is a $u_1\in U_1$ such that $w_1=u_1\otimes_{\mathbb C}(\Omega_2\otimes_{\overline{\mathbb Q}}1_{\mathbb C})$. We have
$$
u_1=\alpha_1(e_1\otimes_{\overline{\mathbb Q}}1_{\mathbb C})+\alpha_2(e_2\otimes_{\overline{\mathbb Q}}1_{\mathbb C})+\ldots+\alpha_{g_1}(e_{g_1}\otimes_{\overline{\mathbb Q}}1_{\mathbb C}),\qquad \alpha_i\in{\mathbb C}.
$$
For $i=1,\ldots,g_1$, the coefficient of $(e_i\otimes_{\overline{\mathbb Q}}f_1)\otimes_{\overline{\mathbb Q}}1_{\mathbb C}$ in the expansion of $w_1$ in terms of the $b\otimes_{\overline{\mathbb Q}}1_{\mathbb C}$, $b\in{\mathcal B}$, equals $\alpha_i$. As $w_1\in \{w\otimes_{\overline{\mathbb Q}}1_{\mathbb C}: w\in W_{\overline{\mathbb Q}}\}$, we have $\alpha_i\in {\overline{\mathbb Q}}$, $i=1,\ldots,g_1$. Therefore $u_1\in U_1\cap\{u\otimes_{\overline{\mathbb Q}}1_{\mathbb C}: u\in V_{1,{\overline{\mathbb Q}}}\}$. Now let $w_k$, $k=1,\ldots,h$, $h=\dim_{\mathbb C}W_1$, be a basis of $W_1$ in $W_1\cap\{w\otimes_{\overline{\mathbb Q}}1_{\mathbb C}: w\in W_{\overline{\mathbb Q}}\}$. By the above discussion $w_k=(u_k\otimes_{\overline{\mathbb Q}}\Omega_2)\otimes_{\overline{\mathbb Q}}1_{\mathbb C}$, for some $u_k\in V_{1,{\overline{\mathbb Q}}}$ with $u_k\otimes_{\overline{\mathbb Q}}1_{\mathbb C}\in U_1$. Let $u\in U_1$. Then $u\otimes_{\mathbb C}(\Omega_2\otimes_{\overline{\mathbb Q}}1_{\mathbb C})\in W_1$ and we have an expansion of the form
$$
u\otimes_{\mathbb C}(\Omega_2\otimes_{\overline{\mathbb Q}}1_{\mathbb C})=\sum_{k=1}^hc_kw_k=\sum_{k=1}^hc_k((u_k\otimes_{\overline{\mathbb Q}}\Omega_2)\otimes_{\overline{\mathbb Q}}1_{\mathbb C}),
$$
$$
=\left(\sum_{k=1}^h  c_k(u_k\otimes_{\overline{\mathbb Q}}1_{\mathbb C})\right)\otimes_{\mathbb C} (\Omega_2\otimes_{\overline{\mathbb Q}}1_{\mathbb C}), \qquad c_k\in{\mathbb C}.
$$
Therefore
$$
u=\sum_{k=1}^h  c_k(u_k\otimes_{\overline{\mathbb Q}}1_{\mathbb C}),\qquad c_k\in{\mathbb C}.
$$
As $u$ is an arbitrary element of $U_1$, we see that $U_1$ has the spanning set $\{u_k\otimes_{\overline{\mathbb Q}}1_{\mathbb C}\}_{k=1}^h\subset U_1$. Therefore $U_1$ has a basis in $U_1\cap\{u\otimes_{\overline{\mathbb Q}}1_{\mathbb C}: u\in V_{1,{\overline{\mathbb Q}}}\}$, as required (in fact $h=\dim_{\mathbb C}(U_1)=\dim_{\mathbb C}(W_1)$, so this spanning set is a basis).

\medskip

\noindent{\bf Lemma 7:} Let $V_{1,{\overline{\mathbb Q}}}$ and $V_{2,{\overline{\mathbb Q}}}$ be two ${\overline{\mathbb Q}}$-vector spaces and $V_{\overline{\mathbb Q}}=V_{1,{\overline{\mathbb Q}}}\otimes_{\overline{\mathbb Q}}V_{2,{\overline{\mathbb Q}}}$. Let $W_1\not=\{0\}$ be a ${\mathbb C}$-vector subspace of $V_{1,{\mathbb C}}=V_{1,{\overline{\mathbb Q}}}\otimes_{\overline{\mathbb Q}}{\mathbb C}$ and $W_2\not=\{0\}$ a ${\mathbb C}$-vector subspace of $V_{2,{\mathbb C}}=V_{2,{\overline{\mathbb Q}}}\otimes_{\overline{\mathbb Q}}{\mathbb C}$. Let $\Omega_1\in V_{1,{\overline{\mathbb Q}}}$ with $\Omega_1\otimes_{\overline{\mathbb Q}}1_{\mathbb C}\notin W_1$, and $\Omega_2\in V_{2,{\overline{\mathbb Q}}}$ with $\Omega_2\otimes_{\overline{\mathbb Q}}1_{\mathbb C}\notin W_2$. Suppose that the ${\mathbb C}$-vector subspace $W_{\mathbb C}$ of $V_{\mathbb C}=V_{\overline{\mathbb Q}}\otimes_{\overline{\mathbb Q}}{\mathbb C}$ given by the internal direct sum
$$
W_{\mathbb C}=(W_1\otimes_{\mathbb C}{\mathbb C}(\Omega_2\otimes_{\overline{\mathbb Q}}1_{\mathbb C}))\oplus({\mathbb C}(\Omega_1\otimes_{\overline{\mathbb Q}}1_{\mathbb C})\otimes_{\mathbb C}W_2)
$$
has a basis in $W_{\mathbb C}\cap \{v\otimes_{\overline{\mathbb Q}}1_{\mathbb C}:v\in V_{\overline{\mathbb Q}}\}$. Then $W_1'={\mathbb C}(\Omega_1\otimes_{\overline{\mathbb Q}}1_{\mathbb C})\oplus W_1$ has a basis in $W_1'\cap\{v_1\otimes_{\overline{\mathbb Q}}1_{\mathbb C}:v_1\in V_{1,{\overline{\mathbb Q}}}\}$ and $W_2'={\mathbb C}(\Omega_2\otimes_{\overline{\mathbb Q}}1_{\mathbb C})\oplus W_2$ has a basis in $W_2'\cap\{v_2\otimes_{\overline{\mathbb Q}}1_{\mathbb C}:v_2\in V_{2,{\overline{\mathbb Q}}}\}$.

\medskip

\noindent{\bf Proof:} For $i=1,2$, let $g_i=\dim_{\overline{\mathbb Q}}(V_{i,{\overline{\mathbb Q}}})$, and $h_i=\dim_{\mathbb C}(W_i)$. Let $\{e_i\}_{i=1}^{g_1}$ be a ${\overline{\mathbb Q}}$-basis of $V_{1,{\overline{\mathbb Q}}}$ and $\{f_j\}_{j=1}^{g_2}$ be a ${\overline{\mathbb Q}}$-basis of $V_{2,{\overline{\mathbb Q}}}$. As $\Omega_1$, $\Omega_2$ are not zero, by reordering these bases and rescaling them, if necessary, we can assume that $\Omega_1$ is of the form
$$
\Omega_1=e_1+\sum_{i=2}^{g_1}a_ie_i,\qquad a_i\in{\overline{\mathbb Q}}
$$
and that $\Omega_2$ is of the form
$$
\Omega_2=f_1+\sum_{j=2}^{g_2}b_jf_j,\qquad b_j\in{\overline{\mathbb Q}}.
$$
Let $w\in V_{\overline{\mathbb Q}}$ be such that $w\otimes_{\overline{\mathbb Q}}1_{\mathbb C}\in W_{\mathbb C}$. As every tensor in $W_1\otimes_{\mathbb C}{\mathbb C}(\Omega_2\otimes_{\overline{\mathbb Q}}1_{\mathbb C})$ and ${\mathbb C}(\Omega_1\otimes_{\overline{\mathbb Q}}1_{\mathbb C})\otimes_{\mathbb C}W_2$ is elementary, there are elements $\omega\in W_1$ and $\eta\in W_2$ such that
$$
w\otimes_{\overline{\mathbb Q}}1_{\mathbb C}=\omega\otimes_{\mathbb C}(\Omega_2\otimes_{\overline{\mathbb Q}}1_{\mathbb C})+(\Omega_1\otimes_{\overline{\mathbb Q}}1_{\mathbb C})\otimes_{\mathbb C}\eta.
$$
Suppose that
$$
\omega=\sum_{i=1}^{g_1}\omega_i(e_i\otimes_{\overline{\mathbb Q}}1_{\mathbb C}),\quad \omega_i\in{\mathbb C},\qquad \eta=\sum_{j=1}^{g_2}\eta_j(f_j\otimes_{\overline{\mathbb Q}}1_{\mathbb C}),\quad \eta_j\in{\mathbb C}.
$$
We can rewrite $w\otimes_{\overline{\mathbb Q}}1_{\mathbb C}$ (uniquely) in the form
$$
w\otimes_{\overline{\mathbb Q}}1_{\mathbb C}=(\omega_1+\eta_1)((\Omega_1\otimes_{\overline{\mathbb Q}}\Omega_2)\otimes_{\overline{\mathbb Q}}1_{\mathbb C})
$$
$$
+(\omega-\omega_1(\Omega_1\otimes_{\overline{\mathbb Q}}1_{\mathbb C}))\otimes_{\mathbb C}(\Omega_2\otimes_{\overline{\mathbb Q}}1_{\mathbb C})+(\Omega_1\otimes_{\overline{\mathbb Q}}1_{\mathbb C})\otimes_{\mathbb C}(\eta-\eta_1(\Omega_2\otimes_{\overline{\mathbb Q}}1_{\mathbb C}))
$$
Notice that,
$$
\omega-\omega_1(\Omega_1\otimes_{\overline{\mathbb Q}}1_{\mathbb C})\in {\mathbb C}(\Omega_1\otimes_{\overline{\mathbb Q}}1_{\mathbb C})\oplus W_1,
$$
and
$$
\eta-\eta_1(\Omega_2\otimes_{\overline{\mathbb Q}}1_{\mathbb C})\in {\mathbb C}(\Omega_2\otimes_{\overline{\mathbb Q}}1_{\mathbb C})\oplus W_2.
$$
Using the expansion of $\Omega_i$ in terms of the chosen basis for $V_{i,{\overline{\mathbb Q}}}$, $i=1,2$, we see that
$$
w\otimes_{\overline{\mathbb Q}}1_{\mathbb C}=\delta_0((\Omega_1\otimes_{\overline{\mathbb Q}}\Omega_2)\otimes_{\overline{\mathbb Q}}1_{\mathbb C})
$$
$$
+(\alpha_2(e_2\otimes_{\overline{\mathbb Q}}1_{\mathbb C})+\ldots+\alpha_{g_1}(e_{g_1}\otimes_{\overline{\mathbb Q}}1_{\mathbb C}))\otimes_{\mathbb C}((f_1+b_2f_2+\ldots+b_{g_2}f_2)\otimes_{\overline{\mathbb C}}1_{\mathbb C})
$$
$$
+((e_1+a_2e_2+\ldots+a_{g_1}e_{g_1})\otimes_{\overline{\mathbb C}}1_{\mathbb C})\otimes_{\mathbb C}(\beta_2(f_2\otimes_{\overline{\mathbb Q}}1_{\mathbb C})+\ldots+\beta_{g_2}(f_{g_2}\otimes_{\overline{\mathbb Q}}1_{\mathbb C})),
$$
where $\delta_0=\omega_1+\eta_1$, $\alpha_i, \beta_j\in{\mathbb C}, i=2,\ldots,g_1;\,j=2,\ldots,g_2$. Moreover,
$$
\alpha_2(e_2\otimes_{\overline{\mathbb Q}}1_{\mathbb C})+\ldots+\alpha_{g_1}(e_{g_1}\otimes_{\overline{\mathbb Q}}1_{\mathbb C})\in {\mathbb C}(\Omega_1\otimes_{\overline{\mathbb Q}}1_{\mathbb C})\oplus W_1
$$
equals $\omega-\omega_1(\Omega_1\otimes_{\overline{\mathbb Q}}1_{\mathbb C})$, and
$$
\beta_2(f_2\otimes_{\overline{\mathbb Q}}1_{\mathbb C})+\ldots+\beta_{g_2}(f_{g_2}\otimes_{\overline{\mathbb Q}}1_{\mathbb C})\in {\mathbb C}(\Omega_2\otimes_{\overline{\mathbb Q}}1_{\mathbb C})\oplus W_2
$$
equals $\eta-\eta_1(\Omega_2\otimes_{\overline{\mathbb Q}}1_{\mathbb C})$.
Now, by assumption, we have $w\in V_{\overline{\mathbb Q}}$ and therefore $w\otimes_{\overline{\mathbb Q}}1_{\mathbb C}$ is of the form,
$$
w\otimes_{\overline{\mathbb Q}}1_{\mathbb C}=\sum_{i=1}^{g_1}\sum_{j=1}^{g_2}(w_{ij}(e_i\otimes_{\overline{\mathbb Q}}f_j))\otimes_{\overline{\mathbb Q}}1_{\mathbb C},
$$
where $w_{ij}\in{\overline{\mathbb Q}}$,  and the $e_i\otimes_{\overline{\mathbb Q}}f_j$, $i=1,\ldots,g_1$, $j=1,\ldots,g_2$, are ${\overline{\mathbb Q}}$-linearly independent and form a basis of $V_{\overline{\mathbb Q}}$.

\noindent Comparing the coefficient of $(e_1\otimes_{\overline{\mathbb Q}}f_1)\otimes_{\overline{\mathbb Q}}1_{\mathbb C}$ in the last two expressions for $w\otimes_{\overline{\mathbb Q}}1_{\mathbb C}$, we see that $\delta_0=w_{11}\in{\overline{\mathbb Q}}$, on comparing that of $(e_i\otimes_{\overline{\mathbb Q}}f_1)\otimes_{\overline{\mathbb Q}}1_{\mathbb C}$, we see that $(\alpha_i+\delta_0a_ib_1)=w_{i1}\in{\overline{\mathbb Q}}$, so that $\alpha_i\in{\overline{\mathbb Q}}$, $i=2,\ldots,g_1$. Finally, comparing the coefficient of $(e_1\otimes_{\overline{\mathbb Q}}f_j)\otimes_{\overline{\mathbb Q}}1_{\mathbb C}$ in the last two expressions for $w\otimes_{\overline{\mathbb Q}}1_{\mathbb C}$, we see that $(\beta_j+\delta_0a_1b_j)=w_{1j}\in{\overline{\mathbb Q}}$, so that $\beta_j\in{\overline{\mathbb Q}}$, $j=2,\ldots,g_2$. Therefore, when we write $w\otimes_{\overline{\mathbb Q}}1_{\mathbb C}$ (uniquely) in this form, the coefficients $\delta_0$, $\alpha_i$, $i=2,\ldots,g_1$, $\beta_j$, $j=2,\ldots,g_2$, are all in ${\overline{\mathbb Q}}$. Of course, by assumption, the $a_i$, $i=1,\ldots, g_1$, and the $b_j$, $j=1,\ldots,g_2$, are all in ${\overline{\mathbb Q}}$.

\noindent Let $h=\dim_{\mathbb C}(W_{\mathbb C})$. By the hypotheses of the lemma, there are $w_k\otimes_{\overline{\mathbb Q}}1_{\mathbb C}\in W_{\mathbb C}$ with $w_k\in V_{\overline{\mathbb Q}}$, $k=1,\ldots,h$, that form a basis of $W_{\mathbb C}$. By our discussion, we can write $w_k$ (uniquely) in the form
$$
w_k=\delta_{k0}(\Omega_1\otimes_{\overline{\mathbb Q}}\Omega_2)
$$
$$
+(\alpha_{k2}e_2+\ldots+\alpha_{kg_1}e_{g_1})\otimes_{\overline{\mathbb Q}}(f_1+b_{2}f_2+\ldots+b_{g_2}f_2)
$$
$$
+(e_1+a_2e_2+\ldots+a_{g_2}e_2)\otimes_{\overline{\mathbb Q}}(\beta_{k2}f_2+\ldots+\beta_{kg_2}f_{g_2}),
$$
where $\delta_{k0}$, $\alpha_{ki}$, $\beta_{kj}\in{\overline{\mathbb Q}}$, $i=2,\ldots,g_1$, $j=2,\ldots,g_2$, and
$$
\sum_{i=2}^{g_1}\alpha_{ki}(e_i\otimes_{\overline{\mathbb Q}}1_{\mathbb C})\in ({\mathbb C}(\Omega_1\otimes_{\overline{\mathbb Q}}1_{\mathbb C})\oplus W_1)\cap\{v_1\otimes_{\overline{\mathbb Q}}1_{\mathbb C}: v_1\in V_{1,{\overline{\mathbb Q}}}\},
$$
and
$$
\sum_{j=2}^{g_2}\beta_{kj}(f_j\otimes_{\overline{\mathbb Q}}1_{\mathbb C})\in ({\mathbb C}(\Omega_2\otimes_{\overline{\mathbb Q}}1_{\mathbb C})\oplus W_2)\cap\{v_2\otimes_{\overline{\mathbb Q}}1_{\mathbb C}: v_2\in V_{2,{\overline{\mathbb Q}}}\}.
$$

\noindent Again, let $\omega\in W_1$ with expansion
$$
\omega=\sum_{i=1}^{g_1}\omega_i(e_i\otimes_{\overline{\mathbb Q}}1_{\mathbb C}),\qquad \omega_i\in{\mathbb C},
$$
in terms of the chosen basis $\{e_i\otimes_{\overline{\mathbb Q}}1_{\mathbb C}\}_{i=1}^{g_1}$ of $V_{1,{\mathbb C}}$.
Then  $\omega\otimes_{\mathbb C}(\Omega_2\otimes_{\overline{\mathbb Q}}1_{\mathbb C})\in W_{\mathbb C}$.
Writing $\omega\otimes_{\mathbb C}(\Omega_2\otimes_{\overline{\mathbb Q}}1_{\mathbb C})$ in terms of the basis $w_k\otimes_{\overline{\mathbb Q}}1_{\mathbb C}$, $k=1,\ldots,h$, of $W_{\mathbb C}$ we have, for complex numbers $c_k$, $k=1,\ldots,h$,
$$
\omega\otimes_{\mathbb C}(\Omega_2\otimes_{\overline{\mathbb Q}}1_{\mathbb C})=\left(\sum_{i=1}^{g_1}\omega_i(e_i\otimes_{\overline{\mathbb Q}}1_{\mathbb C})\right)\otimes_{\mathbb C}(\Omega_2\otimes_{\overline{\mathbb Q}}1_{\mathbb C})
$$
$$
=\sum_{k=1}^hc_kw_k=\left(\sum_{k=1}^hc_k\delta_{k0}\right)((\Omega_1\otimes\Omega_2)\otimes_{\overline{\mathbb Q}}1_{\mathbb C})
$$
$$
+\left(\sum_{k=1}^hc_k\left(\sum_{i=2}^{g_1}\alpha_{ki}(e_i\otimes_{\overline{\mathbb Q}}1_{\mathbb C})\right)\right)\otimes_{\mathbb C}((f_1+b_{2}f_2+\ldots+b_{g_2}f_{g_2})\otimes_{\overline{\mathbb Q}}1_{\mathbb C})
$$
$$
+((e_1+a_{2}e_2+\ldots+a_{g_1}e_{g_1})\otimes_{\overline{\mathbb Q}}1_{\mathbb C})\otimes_{\mathbb C}\left(\sum_{k=1}^hc_k\left(\sum_{j=2}^{g_2}\beta_{kj}(f_j\otimes_{\overline{\mathbb Q}}1_{\mathbb C})\right)\right).
$$
Looking at the coefficient of $(e_1\otimes_{\overline{\mathbb Q}}f_1)\otimes_{\overline{\mathbb Q}}1_{\mathbb C}$ on both sides of the above expression, we have $\sum_{k=1}^hc_k\delta_{k0}=\omega_1$. For $j=2,\ldots,g_2$, looking at the coefficient of $(e_1\otimes_{\overline{\mathbb Q}}f_j)\otimes_{\overline{\mathbb Q}}1_{\mathbb C}$ we have
$$
\left(\sum_{k=1}^hc_k\delta_{k0}\right)b_j+\left(\sum_{k=1}^hc_k\beta_{kj}\right)=\omega_1b_j+\left(\sum_{k=1}^hc_k\beta_{kj}\right)=\omega_1b_j.
$$
Therefore, for $j=2,\ldots,g_2$, we have $\sum_{k=1}^hc_k\beta_{kj}=0$. By the way, this is clear {\it a priori} by the uniqueness of the decomposition
$$
\omega\otimes_{\mathbb C}(\Omega_2\otimes_{\overline{\mathbb Q}}1_{\mathbb C})
$$
$$
=\omega_1((\Omega_1\otimes_{\overline{\mathbb Q}}\Omega_2)\otimes_{\overline{\mathbb Q}}1_{\mathbb C})+(\omega-\omega_1(\Omega_1\otimes_{\overline{\mathbb Q}}1_{\mathbb C}))\otimes_{\mathbb C}(\Omega_2\otimes_{\overline{\mathbb Q}}1_{\mathbb C}).
$$
Therefore,
$$
\omega\otimes_{\mathbb C}(\Omega_2\otimes_{\overline{\mathbb Q}}1_{\mathbb C})=\left(\sum_{k=1}^hc_k\delta_{k0}\right)((\Omega_1\otimes\Omega_2)\otimes_{\overline{\mathbb Q}}1_{\mathbb C})
$$
$$
+\left(\sum_{k=1}^hc_k\left(\sum_{i=2}^{g_1}\alpha_{ki}(e_i\otimes_{\overline{\mathbb Q}}1_{\mathbb C})\right)\right)\otimes_{\mathbb C}(\Omega_2\otimes_{\overline{\mathbb Q}}1_{\mathbb C})
$$
$$
=\left(\omega_1(\Omega_1\otimes_{\overline{\mathbb Q}}1_{\mathbb C})+\sum_{k=1}^hc_k\left(\sum_{i=2}^{g_1}(\alpha_{ki}e_i)\otimes_{\overline{\mathbb Q}}1_{\mathbb C}\right)\right)\otimes_{\mathbb C}(\Omega_2\otimes_{\overline{\mathbb Q}}1_{\mathbb C}).
$$
It follows that,
$$
\omega=\omega_1(\Omega_1\otimes_{\overline{\mathbb Q}}1_{\mathbb C})+(\omega-\omega_1(\Omega_1\otimes_{\overline{\mathbb Q}}1_{\mathbb C}))
$$
$$
=\omega_1(\Omega_1\otimes_{\overline{\mathbb Q}}1_{\mathbb C})+\sum_{k=1}^hc_k\left(\sum_{i=2}^{g_1}(\alpha_{ki}e_i)\otimes_{\overline{\mathbb Q}}1_{\mathbb C}\right).
$$
In other words, if $\omega\in W_1$, then
$$
\omega-\omega_1(\Omega_1\otimes_{\overline{\mathbb Q}}1_{\mathbb C})=\sum_{k=1}^hc_k\left(\sum_{i=2}^{g_1}((\alpha_{ki}e_i)\otimes_{\overline{\mathbb Q}}1_{\mathbb C})\right)
$$
We have already observed that, for all $k=1,\ldots,h$,
$$
\sum_{i=2}^{g_1}((\alpha_{ki}e_i)\otimes_{\overline{\mathbb Q}}1_{\mathbb C})\in ({\mathbb C}(\Omega_1\otimes_{\overline{\mathbb Q}}1_{\mathbb C})\oplus W_1)\cap \{v_1\otimes_{\overline{\mathbb Q}}1_{\mathbb C}: v_1\in V_{1,{\overline{\mathbb Q}}}\}.
$$
Let
$$
W_1'={\mathbb C}(\Omega_1\otimes_{\overline{\mathbb Q}}1_{\mathbb C})\oplus W_1.
$$
Then $\omega$ is in the span of the set ${\mathcal S}_1$ with elements
$$
\Omega_1\otimes_{\overline{\mathbb Q}}1_{\mathbb C},\quad \left(\sum_{i=2}^{g_1}\alpha_{ki}e_i\right)\otimes_{\overline{\mathbb Q}}1_{\mathbb C}, \quad k=1,\ldots,h,
$$
all of which lie in $W_1'\cap \{v_1\otimes_{\overline{\mathbb Q}}1_{\mathbb C}: v_1\in V_{1,{\overline{\mathbb Q}}}\}$. But $\omega$ was an arbitrary element of $W_1$, so it follows that every element in $W_1$ is in the span of ${\mathcal S}_1$. Clearly ${\mathcal S}_1$ also spans $W_1'$ and so some subset of ${\mathcal S}_1$ will be a basis of $W_1'$ in $W_1'\cap\{v_1\otimes_{\overline{\mathbb Q}}1_{\mathbb C}: v_1\in V_{1,{\overline{\mathbb Q}}}\}$, as required. Notice that we may really need to consider $W_1'$ and not just $W_1$, as there is no reason for the coefficient $\omega_1$ of $e_1\otimes_{\overline{\mathbb Q}}1_{\mathbb C}$ in the expansion of $\omega$ to be algebraic. We can apply analogous arguments to $W_2$ and $W_2'$. This completes the proof of the lemma.

\medskip

\noindent Let $k\subseteq k'$ be two fields and $V_k$ a $k$-vector space. We denote by $V_{k'}$ the $k'$-vector space $V_k\otimes_kk'$.  Let $K={\overline{\mathbb Q}}$ or ${\mathbb C}$ and view ${\overline{\mathbb Q}}$ as a subfield of ${\mathbb C}$ by fixing an embedding of ${\overline{\mathbb Q}}$ into ${\mathbb C}$. A hermitian form $h_K$ on a $K$-vector space $V_K$ is a map
$$
h_K:V_K\times V_K\rightarrow K
$$
$$
(z,w)\mapsto h_K(z,w)
$$
which is linear in the first variable, anti-linear in the second variable, and satisfies $h_K(w,z)=\overline{h_K(z,w)}$. In particular, the number $h(z,z)$ is in $K\cap{\mathbb R}$. We say that $h_K$ is definite on $V_K$ if $h(z,z)\not=0$ for all $z\in V_K$ with $z\not=0$. Two elements $z,w\in V_K$ are said to be orthogonal with respect to $h_K$ if $h_K(z,w)=0$. An orthogonal set in $V_K$ is a subset of $V_K$ whose elements are mutually orthogonal with respect to $h_K$. A direct sum of two vector subspaces $U_1$ and $U_2$ of $V_K$ is an orthogonal direct sum if every element of $U_1$ is orthogonal, with respect to $h_K$, to every element of $U_2$.

\medskip

\noindent{\bf Lemma 8:}  Let $U_K$ be a finite dimensional $K$-vector space equipped with a hermitian form $h_K$ which is definite on $U_K$. Then $U_K$ has a basis which is an orthogonal set with respect to $h_K$.

\medskip

\noindent {\bf Proof:}  This follows from the Gram-Schmidt orthogonalization process. Let $\{x_1,\ldots,x_h\}$, $h=\dim_K(U_K)$, be any basis of $U_K$. As the result is trivial for $h=1$, we can suppose that $h\ge 2$. Let
$$
v_1=x_1
$$
$$
v_k=x_k-\sum_{j=1}^{k-1}\frac{h_K(x_k,v_j)}{h_K(v_j,v_j)}\,v_j,\qquad k=2,\ldots,h,
$$
which are well-defined as, for $k=1,\ldots,h-1$, we have $v_k\not=0$, by the linear independence of the $x_i, i=1,\ldots,h$, and therefore $h_K(v_k,v_k)\not=0$, as $h_K$ is definite on $U_K$.
We verify inductively that $h_K(v_i,v_k)=0$, for $i\not=k$, and $i,k=1,\ldots,h$. Indeed,
$$
h_K(v_1,v_2)=h_K(x_1,v_2)
$$
$$
=h_K(x_1,x_2)-h_K\left(x_1,\frac{h_K(x_2,x_1)}{h_K(x_1,x_1)}x_1\right)
$$
$$
=h_K(x_1,x_2)-{\overline{h_K(x_2,x_1)}}\frac{h_K(x_1,x_1)}{h_K(x_1,x_1)}
$$
$$
=h_K(x_1,x_2)-{\overline{h_K(x_2,x_1)}}=0.
$$
The result of the lemma therefore follows for $h=2$. Suppose that $h\ge 3$ and that the set $\{v_1,\ldots,v_{k-1}\}$ is orthogonal with respect to $h_K$ for some $k$ in the range $2\le k\le h$. We show that this implies that $\{v_1,\ldots, v_k\}$ is orthogonal. It suffices to check that $h(v_i,v_k)=0$ for $i=1,\ldots,k-1$. We have, for $i=1,\ldots,k-1$,
$$
h_K(v_i,v_k)=h(v_i,x_k)-h\left(v_i,\sum_{j=1}^{k-1}\frac{h_K(x_k,v_j)}{h_K(v_j,v_j)}v_j\right)
$$
$$
=h_K(v_i,x_k)-\sum_{j=1}^{k-1}\frac{{\overline{h_K(x_k,v_j)}}}{h_K(v_j,v_j)}h_K(v_i,v_j).
$$
By the induction hypothesis $h(v_i,v_j)=0$ for $i\not=j$ and $i,j=1,\ldots,k-1$. Therefore, for $i=1,\ldots,k-1$,
$$
h_K(v_i,v_k)=h_K(v_i,x_k)-\frac{{\overline{h_K(x_k,v_i)}}}{h_K(v_i,v_i)}h_K(v_i,v_i)
$$
$$
=h_K(v_i,x_k)-{\overline{h_K(x_k,v_i)}}=0,
$$
as required. Therefore, by induction, the set $\{v_1,\ldots,v_h\}$ is orthogonal with respect to $h_K$, and so its elements are also linearly independent. The lemma follows.

\medskip

\noindent {\bf Lemma 9:} We choose an embedding of ${\overline{\mathbb Q}}$ into ${\mathbb C}$, and thereby identify ${\overline{\mathbb Q}}$ with a subfield of ${\mathbb C}$. Let $V_{\overline{\mathbb Q}}$ be a finite dimensional ${\overline{\mathbb Q}}$-vector space equipped with a hermitian form $h_{\overline{\mathbb Q}}$. Let $h_{\mathbb C}$ be the extension of $h_{\overline{\mathbb Q}}$ to $V_{\mathbb C}=V_{\overline{\mathbb Q}}\otimes_{\overline{\mathbb Q}}{\mathbb C}$ by linearity in the first variable and anti-linearity in the second variable. Namely, for $e_i, f_j\in V_{\overline{\mathbb Q}}$ and $\lambda_i,\mu_j\in{\mathbb C}$, we set
$$
h_{\mathbb C}\left(\sum_{i}e_i\otimes_{\overline{\mathbb Q}}\lambda_i,\sum_{j}f_j\otimes_{\overline{\mathbb Q}}\mu_j\right)=\sum_{i,j}\lambda_i{\overline{\mu_j}}\,h_{\overline{\mathbb Q}}(e_i,f_j).
$$
Let $U_{1,{\overline{\mathbb Q}}}$ be a ${\overline{\mathbb Q}}$-vector subspace of $V_{\overline{\mathbb Q}}$ and let $U_{1,{\mathbb C}}=U_{1,{\overline{\mathbb Q}}}\otimes_{\overline{\mathbb Q}}{\mathbb C}$. Suppose that $U_{1,{\overline{\mathbb Q}}}$ decomposes as an orthogonal direct sum, with respect to $h_{\overline{\mathbb Q}}$, of ${\overline{\mathbb Q}}$-vector subspaces and that $h_{\overline{\mathbb Q}}$ is definite on each summand. Suppose, in addition, that there is a ${\mathbb C}$-vector subspace $U_2$ of $V_{\mathbb C}$ such that
$$
V_{\mathbb C}=U_{1,{\mathbb C}}\oplus U_2
$$
is an orthogonal direct sum with respect to $h_{\mathbb C}$. Then $U_2$ has a basis in $U_2\cap\{v\otimes_{\overline{\mathbb Q}}1_{\mathbb C}: v\in V_{\overline{\mathbb Q}}\}$.

\medskip

\noindent{\bf Proof:} Let $h_1=\dim_{\overline{\mathbb Q}}(U_{1,{\overline{\mathbb Q}}})$ and
$$
{\mathcal B}_1=\{v_1^{(1)},\ldots,v_{h_1}^{(1)}\}
$$
be a basis of $U_{1,{\overline{\mathbb Q}}}$. By assumption $U_{1,{\overline{\mathbb Q}}}$ is an orthogonal direct sum of subspaces on which $h_{\overline{\mathbb Q}}$ is definite. Applying Lemma 8 to each summand, we can, and will, assume that ${\mathcal B}_1$ is an orthogonal set with respect to $h_{\overline{\mathbb Q}}$. Let
$$
h=\dim_{\mathbb C}U_2=\dim_{\overline{\mathbb Q}}(V_{\overline{\mathbb Q}})-\dim_{\overline{\mathbb Q}}(U_{1,{\overline{\mathbb Q}}})
$$
and let ${\mathcal B}_2=\{v_1,\ldots,v_h\}$ be a set of linearly independent vectors in $V_{\overline{\mathbb Q}}$ such that ${\mathcal B}_2\cup{\mathcal B}_1$ is a basis of $V_{\overline{\mathbb Q}}$. Every $v_i\otimes_{\overline{\mathbb Q}}1_{\mathbb C}$, $v_i\in {\mathcal B}_2$, can be decomposed uniquely as a sum of the form
$$
v_i\otimes_{\overline{\mathbb Q}}1_{\mathbb C}=v_{i,1} + v_{i,2},\qquad v_{i,1}\in U_{1,{\mathbb C}},\;v_{i,2}\in U_2.
$$
There are $\lambda_{ij}\in {\mathbb C}$ such that
$$
v_{i,1}=\sum_{j=1}^{h_1}\lambda_{ij}(v_j^{(1)}\otimes_{\overline{\mathbb Q}}1_{\mathbb C}).
$$
As the $v_k^{(1)}$ are mutually orthogonal with respect to $h_{\overline{\mathbb Q}}$ and $v_{i,2}\in U_2$ is orthogonal to each $v_k^{(1)}\otimes_{\overline{\mathbb Q}}1_{\mathbb C}$ with respect to $h_{\mathbb C}$, we have
$$
h_{\overline{\mathbb Q}}(v_i,v_k^{(1)})=h_{\mathbb C}(v_i\otimes_{\overline{\mathbb Q}}1_{\mathbb C}, v_k^{(1)}\otimes_{\overline{\mathbb Q}}1_{\mathbb C})
$$
$$
=h_{\mathbb C}(v_{i,1}+v_{i,2},v_k^{(1)}\otimes_{\overline{\mathbb Q}}1_{\mathbb C})
$$
$$
=
h_{\mathbb C}(v_{i,1},v_k^{(1)}\otimes_{\overline{\mathbb Q}}1_{\mathbb C})
$$
$$
=\lambda_{ik}h_{\overline{\mathbb Q}}(v_k^{(1)},v_k^{(1)}).
$$
Now $h_{\overline{\mathbb Q}}(v_i,v_k^{(1)})$ and $h_{\overline{\mathbb Q}}(v_k^{(1)},v_k^{(1)})\in{\overline{\mathbb Q}}$. Moreover $h_{\overline{\mathbb Q}}(v_k^{(1)},v_k^{(1)})\not=0$ as $h_{\overline{\mathbb Q}}$ is definite on $U_{1,{\overline{\mathbb Q}}}$. Therefore $\lambda_{ik}\in{\overline{\mathbb Q}}$ for all $i,k$, and it follows that we have $v_{i,1}\in \{u\otimes_{\overline{\mathbb Q}}1_{\mathbb C}: u\in U_{1,{\overline{\mathbb Q}}}\}$. Therefore
$$
v_{i,2}=v_i\otimes_{\overline{\mathbb Q}}1_{\mathbb C}-\sum_{j=1}^{h_1}(\lambda_{ij}v_j^{(1)})\otimes_{\overline{\mathbb Q}}1_{\mathbb C}, \qquad i=1,\ldots,h
$$
is in $U_2\cap\{v\otimes_{\overline{\mathbb Q}}1_{\mathbb C}:v\in V_{\overline{\mathbb Q}}\}$. Let ${\mathcal B}_2'=\{v_{1,2},v_{2,2},\ldots,v_{h,2}\}$. It remains to check that ${\mathcal B}_2'$ is a ${\mathbb C}$-linearly independent set. Suppose there are $c_i\in{\mathbb C}$, $i=1,\ldots,h$, such that
$$
\sum_{i=1}^hc_iv_{i,2}=0.
$$
Then
$$
\sum_{i=1}^hc_i(v_i\otimes_{\overline{\mathbb Q}}1_{\mathbb C})-\sum_{j=1}^{h_1}\left(\sum_{i=1}^hc_i\lambda_{ij}\right)(v_j\otimes_{\overline{\mathbb Q}}1_{\mathbb C})=0.
$$
But the elements of the set $\{v\otimes_{\overline{\mathbb Q}}1_{\mathbb C}: v\in{\mathcal B}_2\cup{\mathcal B}_1\}$ are linearly independent over ${\mathbb C}$. In fact, they form a basis of $V_{\mathbb C}$. In particular, in the above expression, the coefficients of the elements $v_i\otimes_{\overline{\mathbb Q}}1_{\mathbb C}$, $v_i\in{\mathcal B}_2$, must all be zero, namely $c_i=0$ for $i=1,\ldots,h$. Therefore ${\mathcal B}_2'$ is a ${\mathbb C}$-linearly independent set and is therefore a basis of $U_2$, as required.

\bigskip

\noindent{\bf Lemma 10:} Let $(V,{\widetilde{\varphi}})$ be a rational Hodge structure of weight $k\ge0$ and $Q$ a non-degenerate bilinear form on $V$. Suppose that
$$
Q(u,v)=(-1)^kQ(v,u),\qquad u,v\in V.
$$
Let $h_{\mathbb C}$ be the hermitian form on $V_{\mathbb C}$ given by
$$
h_{\mathbb C}(u,v)=(\sqrt{-1})^kQ(u,{\overline{v}}),
$$
where $Q$ also denotes the extension of $Q$ to $V_{\mathbb C}$. Suppose that the Hodge decomposition of $V_{\mathbb C}$,
$$
V_{\mathbb C}=\oplus_{p+q=k}V^{p,q},
$$
is orthogonal with respect to $h_{\mathbb C}$ and that $h_{\mathbb C}$ is definite on each $V^{p,q}$. If the Hodge filtration of $(V,{\widetilde{\varphi}})$ is defined over ${\overline{\mathbb Q}}$, then every $V^{p,q}$ has a basis in $V^{p,q}\cap \{v\otimes_{\overline{\mathbb Q}}1_{\mathbb C}: v\in V_{\overline{\mathbb Q}}\}$.

\medskip

\noindent{\bf Proof:} We proceed by induction. The Hodge filtration is given by
$$
F^{p,k}=\oplus_{p'\ge p}V^{p',k-p'},\qquad p=0,\ldots,k.
$$
By the assumptions of the lemma, every ${\mathbb C}$-vector space $F^{p,k}$, $p=0,\ldots,k$, has a basis in $F^{p,k}\cap \{v\otimes_{\overline{\mathbb Q}}1_{\mathbb C}: v\in V_{\overline{\mathbb Q}}\}$. For $p=k$, we have $F^{k,k}=V^{k,0}$, so that $V^{k,0}$ has a basis in $V^{k,0}\cap \{v\otimes_{\overline{\mathbb Q}}1_{\mathbb C}: v\in V_{\overline{\mathbb Q}}\}$. Suppose that $V^{p,q}$ has a basis in $V^{p,q}\cap \{v\otimes_{\overline{\mathbb Q}}1_{\mathbb C}: v\in V_{\overline{\mathbb Q}}\}$ for $p=k,k-1,\ldots,k-i$. We have
$$
F^{k-i-1,i+1}=\oplus_{j\le i+1}V^{k-j,j}=(\oplus_{j\le i}V^{k-j,j})\oplus V^{k-i-1,i+1}=F^{k-i,i}\oplus V^{k-i-1,i+1}.
$$
By the assumptions of the lemma, the ${\mathbb C}$-vector space $F^{k-i,i}$ is orthogonal to $V^{k-i-1,i+1}$ and $h_{\mathbb C}$ is definite on each $V^{k-j,j}$, $j\le i$. By the induction hypothesis, each $V^{k-j,j}$, $j=0,\ldots,i$, has a basis in
$$
V^{k-j,j}\cap \{v\otimes_{\overline{\mathbb Q}}1_{\mathbb C}: v\in V_{\overline{\mathbb Q}}\}.
$$
Therefore, by Lemma 9, the ${\mathbb C}$-vector space $V^{k-i-1,i+1}$ has a basis in
$$
V^{k-i-1,i+1}\cap \{v\otimes_{\overline{\mathbb Q}}1_{\mathbb C}: v\in V_{\overline{\mathbb Q}}\}.$$
The lemma follows by induction.

\bigskip

\section{Some examples of Borcea-Voisin towers}\label{s:exampleslow}

In this section, we give details of the discussion of \S\ref{s:BVtower} for some examples in low dimension. These examples will show that the assumptions that we make in \S\ref{s:BVtower}, Proposition 2,
 are necessary. Some of the arguments from previous sections will be repeated, in order to make this section self-contained. We first look at this construction in the simplest case, when $(A_1,I_1)$ and $(A_2,I_2)$ are elliptic curves with the natural involution.
 The Problem of \S\ref{s:VHS} was also considered for this example in \cite{TCY}, but here we give more details. For $i=1,2$, let $A_i$ be the complex elliptic curve with equation
$$
y_i^2=x_i(x_i-1)(x_i-\lambda_i),\qquad \lambda_i\not=0,1,\infty,
$$
and $I_i$ the involution $(x_i,y_i)\mapsto(x_i,-y_i)$. Then $I_i$ has the 4 fixed points with $y_i=0$, corresponding to the 4 torsion points of order 2 on $A_i$. Let $I_{1,2}$ be the involution $I_1\times I_2$ on $A_1\times A_2$. We have
$$
A_1\times A_2\rightarrow^{I_{1,2}} (A_1\times A_2)/I_{1,2}\rightarrow^\alpha(A_1/I_1)\times(A_2/I_2)=(A_1\times A_2)/G,
$$
where both the maps $I_{1,2}$ and $\alpha$ are generically $2$ to $1$. The composite $\alpha\circ I_{1,2}$ is generically 4 to 1, giving the quotient of $A_1\times A_2$ by the group $G$ of order 4 generated by $I_1\times 1$ and $1\times I_2$. It is 2 to 1 on $A_1\times P_2$, which is fixed by $1\times I_2$, and on $P_1\times A_2$, which is fixed by $I_1\times 1$, where, for $i=1,2$, the point $P_i$ is any of the 4 fixed points of $I_i$ on $A_i$. There are 8 such curves. It is 1 to 1 on the 16 fixed points of $I_{1,2}$ on $A_1\times A_2$. These are the points with $y_1=y_2=0$, that is, the 16 torsion points of order 2 on the abelian surface $A_1\times A_2$. Let $\widehat{A_1\times A_2}$ be the blow-up of $A_1\times A_2$ at these 16 points, and ${\widehat{I}}_{1,2}$ the induced involution on $\widehat{A_1\times A_2}$. The ramification locus of ${\widehat{I}}_{1,2}$ is the union of the 16 distinct rational curves given by the blow-ups of the 16 fixed points of $I_{1,2}$. The Calabi-Yau variety $B=\widehat{A_1\times A_2}/{\widehat{I}}_{1,2}$ is called the Kummer surface associated to $A_1\times A_2$, and gives a minimal resolution of the singular surface $(A_1\times A_2)/I_{1,2}$. The involution $\alpha$ on this last quotient induces an involution $I={\widehat{\alpha}}$ on $B$. We have, for $k=0,1,2$,
$$
H^k(B,{\mathbb C})=H^k(\widehat{A_1\times A_2},{\mathbb C})^{{\widehat{I}}_{1,2}},
$$
where the superscript on the right hand side denotes the ${\widehat{I}}_{1,2}$-invariant part of $H^k(\widehat{A_1\times A_2},{\mathbb C})$. From Lemma 3, \S\ref{s:lemmas}, we have
$$
H^2(\widehat{A_1\times A_2},{\mathbb C})=H^2(A_1\times A_2,{\mathbb C})\oplus W_{\mathbb C},
$$
where $W$ is the underlying ${\mathbb Q}$-vector space of a rational Hodge structure ${\mathbb W}$, independent of $\lambda_i$, $i=1,2$, concentrated in bi-degree $(1,1)$, and stable under ${\widehat{I}}_{1,2}$. The ${\mathbb C}$-vector space $W_{\mathbb C}$ has dimension 16. From Lemma 2, \S\ref{s:lemmas}, we have
$$
H^{2,0}(A_1\times A_2)=H^{1,0}(A_1)\otimes H^{1,0}(A_2),\quad H^{0,2}(A_1\times A_2)= H^{0,1}(A_1)\otimes H^{0,1}(A_2),
$$
$$
H^{1,1}(A_1\times A_2)= (H^{1,0}(A_1)\otimes H^{0,1}(A_2))\oplus (H^{0,1}(A_1)\otimes H^{1,0}(A_2))
$$
$$
\oplus(H^{1,1}(A_1)\otimes H^{0,0}(A_2)\oplus (H^{0,0}(A_1)\otimes H^{1,1}(A_2)).
$$
Therefore $h^{2,0}=h^{0,2}=1$, and $h^{1,1}=4$. For $i=1,2$, the rational Hodge structure $H^0(A_i,{\mathbb Q}_{A_i})$ is concentrated in bi-degree $(0,0)$, and the rational Hodge structure $H^2(A_i,{\mathbb Q}_{A_i})$ is concentrated in bi-degree $(1,1)$. Therefore, the term
$$
(H^{1,1}(A_1)\otimes H^{0,0}(A_2))\oplus (H^{0,0}(A_1)\otimes H^{1,1}(A_2))
$$
is the complexification $U_{\mathbb C}$ of the underlying 2 dimensional ${\mathbb Q}$-vector space $U$ of a rational CM Hodge structure ${\mathbb U}$, by Lemma 1, \S\ref{s:lemmas}. This Hodge structure is independent of $\lambda_1$, $\lambda_2$ and is invariant under $I_{1,2}$. It follows that
$$
H^2(B,{\mathbb C})=H^2(\widehat{A_1\times A_2},{\mathbb C})^{{\widehat{I}}_{1,2}}
$$
$$
=H^2(A_1\times A_2,{\mathbb C})^{I_{1,2}}\oplus W_{\mathbb C}=H^2(A_1\times A_2,{\mathbb C})\oplus W_{\mathbb C}
$$
$$
=(H^{1,0}(A_1)\otimes H^{1,0}(A_2))\oplus (H^{0,1}(A_1)\otimes H^{0,1}(A_2))
$$
$$
\oplus(H^{1,0}(A_1)\otimes H^{0,1}(A_2))\oplus (H^{0,1}(A_1)\otimes H^{1,0}(A_2))\oplus U_{\mathbb C}\oplus W_{\mathbb C}
$$
Observing that $h^{1,1}(B)=20$, we have $h^{2,0}+h^{1,1}+h^{0,2}=22$, which is consistent with our claim that $B$ is a $K3$-surface. The $F^{2,2}$-part of the Hodge filtration $F^{\ast,2}$ of $H^2(B,{\mathbb C})$ is given by
$$
(H^{1,0}(A_1)\otimes_{\mathbb C} H^{1,0}(A_2))\subset H^2(B,{\mathbb C}).
$$
Assume that $F^{2,2}$ has a basis in $H^2(A_1\times A_2,{\overline{\mathbb Q}})$. The ${\mathbb C}$-vector spaces $H^{1,0}(A_1)$ and $H^{1,0}(A_2)$ are 1-dimensional. Let $e_1$ span $H^{1,0}(A_1)$ and $e_2$ span $H^{1,0}(A_2)$. Then $e_1\otimes_{\mathbb C}e_2$ spans $H^{1,0}(A_1)\otimes_{\mathbb C} H^{1,0}(A_2)$. If we assume there is a basis of $H^{1,0}(A_1)\otimes_{\mathbb C} H^{1,0}(A_2)$ in
$$
H^2(A_1\times A_2,{\overline{\mathbb Q}})=(H^0(A_1,{\overline{\mathbb Q}})\otimes_{\overline{\mathbb Q}}H^2(A_2,{\overline{\mathbb Q}}))\oplus (H^2(A_1,{\overline{\mathbb Q}})\otimes_{\overline{\mathbb Q}}H^0(A_2,{\overline{\mathbb Q}}))
$$
$$
\oplus (H^1(A_1,{\overline{\mathbb Q}})\otimes_{\overline{\mathbb Q}}H^1(A_2,{\overline{\mathbb Q}})),
$$
then there is a complex number $\lambda\not=0$ such that $\lambda(e_1\otimes_{\mathbb C} e_2)=e_{1,2}\otimes_{\overline{\mathbb Q}} 1_{\mathbb C}$, where $e_{1,2}$ lies in the 4-dimensional ${\overline{\mathbb Q}}$-vector space $H^1(A_1,{\overline{\mathbb Q}})\otimes_{\overline{\mathbb Q}} H^1(A_2,{\overline{\mathbb Q}})$. Therefore, there is a non-zero complex number $\alpha$ such that $\lambda\alpha e_1\in H^1(A_1,{\overline{\mathbb Q}})$ and $\alpha^{-1}e_2\in H^1(A_2,{\overline{\mathbb Q}})$. It follows that, for $i=1,2$, the Hodge filtration $F^{\ast,1}(A_i)$ of $A_i$ is defined over ${\overline{\mathbb Q}}$. If we also have $\lambda_i\in{\overline{\mathbb Q}}$, for $i=1,2$, then by Schneider's Theorem, the filtration $F^{\ast,1}(A_i)$ has CM, so that $F^{\ast,2}$ has CM  by Lemma 1, \S\ref{s:lemmas}, since the  rational Hodge structure $H^2(B,{\mathbb Q}_B)$ is the tensor product of $H^1(A_1,{\mathbb Q}_{A_1})$ and $H^1(A_2,{\mathbb Q}_{A_2})$, up to a rational Hodge structure concentrated in bi-degree $(1,1)$ and independent of $\lambda_1$, $\lambda_2$. For an argument using a direct computation of normalized periods, see \S\ref{s:examples}, Example 1. One can also check directly that $H^1(B,{\mathbb C})=\{0\}$, which is consistent with $B$ being Calabi-Yau. Therefore $B$ has CM if and only if $A_1$ and $A_2$ both have CM. The ``if'' part again uses Lemma 1, \S\ref{s:lemmas}. The CM property and strong CM property are equivalent in this case. We note in passing that, although we do not use this fact in this example, once we deduce from the assumption that $F^{2,2}(B)$ has a basis in $H^2(B,{\overline{\mathbb Q}})$ that $F^{\ast,1}(A_i)$, $i=1,2$, is defined over ${\overline{\mathbb Q}}$, we can \emph{a fortiori} deduce directly from the above arguments that the whole Hodge filtration $F^{\ast,2}(B)$ of $H^2(B,{\mathbb C})$ is defined over ${\overline{\mathbb Q}}$, without invoking the CM property of $B$. Notice that if the $F^{\ast,1}(A_i)$, $i=1,2$ are not defined over ${\overline{\mathbb Q}}$, it may well happen that $H^{1,1}(B)$ has no basis in $H^2(B,{\overline{\mathbb Q}})$, even though the Mumford-Tate group acts via scalars on this piece of $H^2(B,{\mathbb C})$.

\noindent We now have a new Calabi-Yau variety $B$ with involution $I={\widehat{\alpha}}$. The ramification locus $R$ of $I$ has 8 connected components consisting of smooth rational curves $C_s$, $s=1,\ldots,8$, given by the union of the image under the degree 2 rational map $A_1\times A_2\dashrightarrow B$ of $P_1\times A_2$ and $A_1\times P_2$, where, for $i=1,2$, the point $P_i$ runs over the four 2-torsion on $A_i$. As in \cite{Bor}, we repeat the Borcea-Voisin construction with $(A_3,I_3)=(B,I)$ and $(A_4,I_4)$ the elliptic curve
$$
y^2=x(x-1)(x-\lambda),\qquad \lambda\not=0,1,\infty,
$$
with involution $I_4:(x,y)\mapsto(x,-y)$. We proved an analogue of Schneider's Theorem for this example in \cite{TCY}, but here we, again, give more details. Let $I_{3,4}$ be the involution $I_3\times I_4$ on $B\times A_4$. It has ramification locus the 32 rational curves
$$
Z=\cup_{s=1}^8\cup_{k=1}^4 (C_s\times\{P_k\})
$$
where $P_k$ runs over the four 2-torsion points of $A_4$. Let ${\widehat{B\times A_4}}$ be the blow-up of $B\times A_4$ along $Z$, and ${\widehat{I}}_{3,4}$ the involution induced on it by $I_{3,4}$. The quotient $Y={\widehat{B\times A_4}}/{\widehat{I}}_{3,4}$ is a Calabi-Yau 3-fold with involution $\gamma$ induced by $I_3\times 1$, or equivalently by $1\times I_4$, on $B\times A_4$. The variety $Y$ is also a minimal resolution of
 $$
 A_1\times A_2\times A_4/H
 $$
 where $H$ is the group of order 4 generated by $I_1\times I_2\times 1$ and $1\times I_2\times I_4$. The singularities of this last quotient lie along a configuration of $48=16+32$ rational curves with $4^3$ intersection points. The ramification locus $R$ of $\gamma$ consists of the image under the degree 4 rational map
 $$
A_1\times A_2\times A_4\rightarrow Y
 $$
 of the union of
 $$
P_i\times A_2\times A_4,\quad A_1\times P_j\times A_4,\quad A_1\times A_2\times P_k,
 $$
 where $P_i, P_j, P_k$ run over the 2-torsion points of $A_1, A_2, A_4$, respectively. As $Y$ is a Calabi-Yau 3-fold, or by direct computation, we have that $H^0(Y,{\mathbb C})={\mathbb C}$, $H^{1,0}(Y)=H^{0,1}(Y)=\{0\}$, $H^{2,0}(Y)=H^{0,2}(Y)=\{0\}$. Therefore the rational Hodge structure $H^2(Y,{\mathbb Q}_Y)$ is concentrated in bi-degree (1,1) and therefore has CM by Lemma 1, \S\ref{s:lemmas}. We have, with similar notation to that used above,
 $$
 H^3(Y,{\mathbb C})=H^3({\widehat{B\times A_4}},{\mathbb C})^{{\widehat{I}}_{3,4}}.
 $$
 By Lemma 3, \S\ref{s:lemmas}, as $Z$ is a union of rational curves, we have
 $$
 H^3({\widehat{B\times A_4}},{\mathbb C})=H^3(B\times A_4,{\mathbb C})\oplus H^1(Z,{\mathbb C})(-1)=H^3(B\times A_4,{\mathbb C}).
 $$
 Therefore
 $$
 H^3(Y,{\mathbb C})=H^3(B\times A_4,{\mathbb C})^{I_{3,4}}.
 $$
 Denote by $H^{r,s}(B)^\pm$ the $\pm$-eigenspace for the induced action of $I_3$ on $H^{r,s}(B)$, and by $H^{r,s}(A_4)^\pm$ the $\pm$-eigenspace for the induced action of $I_4$ on $H^{r,s}(A_4)$. By Lemma 2, \S\ref{s:lemmas}, we have, for all $r,s$ with $r+s=3$,
 $$
H^{r,s}(B\times A_4)^{I_{3,4}}=\oplus_{p+p'=r,q+q'=s}(H^{p,q}(B)^+\otimes H^{p',q'}(A_4)^+)
 $$
 $$
\oplus_{p+p'=r,q+q'=s}(H^{p,q}(B)^-\otimes H^{p',q'}(A_4)^-).
 $$
 Since the CM-nature of $B\times A_4$ is inductively determined by the CM nature of three elliptic curves used in its construction, it is easy to see that one only needs to relate $H^{3,0}(Y)$ to the $H^{1,0}(A_i)$, $i=1,2,4$. This is explained in more detail in \S\ref{s:examples}, Example 2. For the sake of demonstrating the general case, we proceed more formally. The $H^{p,q}(B)^+\otimes H^{p',q'}(A_4)^+$ terms in the above sum equal $\{0\}$. Indeed, as $B$ has real dimension 4, we have $H^1(B,{\mathbb C})\simeq H^3(B,{\mathbb C})$, and as it is Calabi-Yau, there are no $H^{3,0}$ and $H^{0,3}$ terms. As $H^{1,0}(A_4)=H^{1,0}(A_4)^-$, and $H^{0,1}(A_4)^-=H^{0,1}(A_4)$ there are no $H^{1,0}(A_4)^+$ and $H^{0,1}(A_4)^+$ terms, so no $H^{p,q}(B)^+\otimes H^{p',q'}(A_4)^+$ terms with $p+q=2$. One therefore has
 $$
 H^{3,0}(B\times A_4)^{I_{3,4}}=H^{2,0}(B)^-\otimes H^{1,0}(A_4)^-=H^{2,0}(B)\otimes H^{1,0}(A_4),
 $$
 and
 $$
 H^{2,1}(B\times A_4)^{I_{3,4}}=(H^{2,0}(B)^-\otimes H^{0,1}(A_4)^-)\oplus(H^{1,1}(B)^-\otimes H^{1,0}(A_4)^-)
 $$
 $$
 =(H^{2,0}(B)\otimes H^{0,1}(A_4))\oplus(H^{1,1}(B)^-\otimes H^{1,0}(A_4)).
 $$
 Therefore
 $$
 H^{3,0}(Y)=H^{2,0}(B)\otimes_{\mathbb C} H^{1,0}(A_4),
 $$
 and
 $$
 H^{2,1}(Y)=(H^{2,0}(B)\otimes_{\mathbb C} H^{0,1}(A_4))\oplus(H^{1,1}(B)^-\otimes_{\mathbb C} H^{1,0}(A_4)).
 $$
 We notice that, by our discussion of $B$ in the preceding paragraph, and the fact that $I_3$ is induced by $I_1\times 1$ (and by $1\times I_2$), we have
 $$
 H^{1,1}(B)^-=(H^{1,0}(A_1)\otimes H^{0,1}(A_2))\oplus (H^{0,1}(A_1)\otimes H^{1,0}(A_2)).
 $$
 Now, if the Hodge filtration $F^{\ast,3}$ of $H^3(Y,{\mathbb C})$ is defined over ${\overline{\mathbb Q}}$, we have in particular that the 1-dimensional ${\mathbb C}$-vector space $H^{3,0}(Y)$ has a generator in $H^2(B,{\overline{\mathbb Q}})\otimes_{\overline{\mathbb Q}}H^1(A_4,{\overline{\mathbb Q}})$. In exactly the same way as we did for the previous step in the tower, we deduce that the 1-dimensional vector spaces $H^{2,0}(B)$ and $H^{1,0}(A_4)$ have generators in $H^2(B,{\overline{\mathbb Q}})$ and $H^1(A_4,{\overline{\mathbb Q}})$, respectively. Since $H^{2,0}(B)=H^{1,0}(A_1)\otimes_{\mathbb C}H^{1,0}(A_2)$, we again deduce that the 1-dimensional vector space $H^{1,0}(A_i)$ has a generator in $H^1(A_i,{\overline{\mathbb Q}})$, $i=1,2$. Although we do not need it here, it follows from our discussion that $H^{2,1}(Y)$ also has a basis in $H^3(Y,{\overline{\mathbb Q}})$. Therefore $F^{\ast,3}(Y)$ is defined over ${\overline{\mathbb Q}}$. By Schneider's Theorem, if we also assume that $\lambda_1$, $\lambda_2$, and $\lambda_4$ are algebraic numbers, the elliptic curves $A_1$, $A_2$, $A_4$ have CM, and therefore $Y$ has CM by Lemma 1, \S\ref{s:lemmas}, see also \cite{Bor}. In fact, it is easy to see, again by Lemma 1, \S\ref{s:lemmas}, that $Y$ has CM if and only if $A_1$, $A_2$, $A_3$ have CM. By the Calabi-Yau property, the CM and strong CM property are equivalent for $Y$.

 \noindent We now consider the Borcea-Voisin construction in the case where the initial varieties are any algebraic $K3$-surface with an involution which is non-trivial on its $(2,0)$-forms and an elliptic curve with the usual involution considered above.  For some general facts about the Borcea-Voisin construction and CM we use the references \cite{Bor}, \cite{Roh}, Chapter 7, \cite{Thac}, \cite{Voi} several times without citing them each time. Suppose that we have a smooth projective family ${\mathcal B}\rightarrow S$ of $K3$-surfaces with involution, defined over ${\overline{\mathbb Q}}$, such that the involution leaves $S$ invariant and is non-trivial on $H^{2,0}({\mathcal B}_s)$, for every fiber ${\mathcal B}_s$ of ${\mathcal B}$. Suppose that, for all $s\in S({\overline{\mathbb Q}})$, if the Hodge filtration $F^{\ast,2}$ of $H^2({\mathcal B}_s,{\mathbb C})$ is defined over ${\overline{\mathbb Q}}$, then ${\mathcal B}_s$ has CM. We fix $s\in S$ and $\lambda_4\in{\mathbb C}\setminus\{0,1,\infty\}$. Let $B={\mathcal B}_s$, with involution $I_3$, and let $A_4=A_4(\lambda_4)$ be the elliptic curve with involution $I_4$ as above. Let $Y={\widehat{B\times A_4}}/{\widehat{I}}_{3,4}$ be the next step in the Borcea-Voisin construction, where ${\widehat{B\times A_4}}$ is the blow-up of $B\times A_4$ along the ramification
 divisor $Z$ of $I_{3,4}=I_3\times I_4$ and ${\widehat{I}}_{3,4}$ is the involution induced on this blow-up by $I_{3,4}$. As $Y$ is a Calabi-Yau 3-fold, we have $H^1(Y,{\mathbb C})=\{0\}$ and $H^2(Y,{\mathbb C})=H^{1,1}(Y)$. Therefore, the CM and strong CM properties are equivalent for $Y$. Now, with an analogous notation to that used in the previous paragraph, we have
 $$
 H^{1,1}(B\times A_4)^{I_{3,4}}=\oplus_{p+p'=1,q+q'=1}(H^{p,q}(B)^+\otimes H^{p',q'}(A_4)^+)
 $$
 $$
\oplus_{p+p'=1,q+q'=1}(H^{p,q}(B)^-\otimes H^{p',q'}(A_4)^-).
 $$
 $$
 =(H^{1,1}(B)^+\otimes H^{0,0}(A_4))\oplus (H^{0,0}(B)^+\otimes H^{1,1}(A_4)).
 $$
 Let $Z$ be the ramification locus of $I_{3,4}$. Then, by Lemma 3, \S\ref{s:lemmas}, we have,
 $$
 H^{1,1}(Y)=(H^{1,1}(B)^+\otimes H^{0,0}(A_4))\oplus(H^{0,0}(B)\otimes H^{1,1}(A_4))\oplus H^{0,0}(Z)(-1).
 $$
 As $H^2(Y,{\mathbb Q}_Y)$ is concentrated in bi-degree $(1,1)$, it has CM, and therefore the corresponding Hodge filtration $F^{\ast,2}$ is defined over ${\overline{\mathbb Q}}$. The first two summands in the above formula clearly share these properties. Let $H^k(B,{\mathbb Q}_B)^\pm$ be the rational Hodge substructure of $H^k(B,{\mathbb Q}_B)$  corresponding to the $\pm$-eigenspace of the induced action of $I_3$. The Hodge structure $H^2(B,{\mathbb Q}_B)^+$ is concentrated in bi-degree $(1,1)$, so has CM and Hodge filtration defined over ${\overline{\mathbb Q}}$, which is consistent with the above formula.

 \noindent For $r,s$ with $r+s=3$,
 $$
H^{r,s}(B\times A_4)^{I_{3,4}}=\oplus_{p+p'=r,q+q'=s}(H^{p,q}(B)^+\otimes H^{p',q'}(A_4)^+)
 $$
 $$
\oplus_{p+p'=r,q+q'=s}(H^{p,q}(B)^-\otimes H^{p',q'}(A_4)^-).
 $$
 and similar arguments show that, in fact,
 $$
 H^{r,s}(B\times A_4)^{I_{3,4}}=\oplus_{p+p'=r,q+q'=s}(H^{p,q}(B)^-\otimes H^{p',q'}(A_4)^-).
 $$
 Again, repeating the reasoning of the previous paragraph,
 $$
 H^{3,0}(B\times A_4)^{I_{3,4}}=H^{2,0}(B)\otimes H^{1,0}(A_4),
 $$
 and
 $$
 H^{2,1}(B\times A_4)^{I_{3,4}}=(H^{2,0}(B)\otimes H^{0,1}(A_4))\oplus(H^{1,1}(B)^-\otimes_{\mathbb C} H^{1,0}(A_4)).
 $$
 Moreover,
 $$
 H^{1,1}(B\times A_4)^{I_{3,4}}=(H^{1,1}(B)^+\otimes H^{0,0}(A_4))\oplus (H^{0,0}(B)\otimes H^{1,1}(A_4)).
 $$
By Lemma 3, \S\ref{s:lemmas}, we have
 $$
 H^3(Y,{\mathbb C})=H^3(B\times A_4)^{I_{3,4}}\oplus H^1(Z,{\mathbb C})(-1),
 $$
 and, as $H^{2,0}(Y)=H^{0,2}(Y)=\{0\}$, and $H^{1,0}(B)=H^{0,1}(B)=\{0\}$,
 $$
 H^2(Y,{\mathbb C})=H^{1,1}(Y)=H^{1,1}(B\times A_4)^{I_{3,4}}\oplus H^{0,0}(Z,{\mathbb C})(-1)
 $$
 $$
 =(H^{1,1}(B)^+\otimes H^{0,0}(A_4))\oplus (H^{0,0}(B)\otimes H^{1,1}(A_4))\oplus H^{0,0}(Z,{\mathbb C})(-1).
 $$
 Unlike the special situation in the preceding paragraph, the $H^1(Z,{\mathbb C})$ term may no longer be trivial. The 1-dimensional $Z$ is the union of $R\times P_i$, where $R$ is the ramification locus of $I_3$ on $B$ and $P_i$ runs over the four 2-torsion points of $A_4$. The fixed locus of $I_3$ on $B$ can be either empty, consist of a finite number of rational curves and at most one curve of positive genus, or consist of two elliptic curves, see \cite{Nik}, \cite{VoiK3}, \cite{Zha}. We therefore make the additional assumption that the ramification loci of the fibers of ${\mathcal B}\rightarrow S$ form a smooth projective family ${\mathcal R}\rightarrow S$ such that, if $s\in S({\overline{\mathbb Q}})$, and the Hodge filtration $F^{\ast,1}$ of $H^1({\mathcal R}_s,{\mathbb C})$ is defined over ${\overline{\mathbb Q}}$, then ${\mathcal R}_s$ has CM (and therefore automatically strong CM). We have
 $$
 H^{3,0}(Y)=H^{2,0}(B)\otimes H^{1,0}(A_4),\quad  H^{0,3}(Y)=H^{0,2}(B)\otimes H^{0,1}(A_4)
 $$
 and
 $$
H^{2,1}(Y)=(H^{2,0}(B)\otimes H^{0,1}(A_4))\oplus(H^{1,1}(B)^-\otimes H^{1,0}(A_4))\oplus H^{1,0}(Z)(-1).
$$
Suppose now that the Hodge filtration $F^{\ast,3}(Y)$ of $H^3(Y,{\mathbb Q}_Y)_{\mathbb C}$ is defined over ${\overline{\mathbb Q}}$. In particular, the 1-dimensional complex vector space $F^{3,3}(Y)=H^{3,0}(Y)$ is therefore spanned by an element of $H^3(Y,{\overline{\mathbb Q}})$.
We deduce, using our previous arguments, that $H^{2,0}(B)$ is spanned by an element of $H^2(B,{\overline{\mathbb Q}})$, and $H^{1,0}(A_4)$ is spanned by an element of $H^1(A,{\overline{\mathbb Q}})$. Suppose, in addition, that $s\in S({\overline{\mathbb Q}})$ and $\lambda_4$ is an algebraic number. Then, by Schneider's Theorem, the elliptic curve $A_4$ has CM, that is, the rational Hodge structure $H^1(A_4,{\mathbb Q}_{A_4})$ has CM. We have the Hodge decomposition of $H^2(B,{\mathbb Q}_B)_{\mathbb C}$ given by
$$
H^2(B,{\mathbb C})=H^{2,0}(B)\oplus H^{1,1}(B)\oplus H^{0,2}(B).
$$
The complex vector space $H^{2,0}(B)$ is spanned by $e\otimes_{\overline{\mathbb Q}} 1_{\mathbb C}$, for some element $e$ of $V=H^2(B,{\overline{\mathbb Q}})$. As ${\overline{H^{0,2}}}(B)=H^{2,0}(B)$, it follows that $H^{0,2}(B)$ is spanned by an element ${\overline{e}}\otimes_{\overline{\mathbb Q}} 1_{\mathbb C}$ where $\overline{e}\in V$. Consider the non-degenerate symmetric pairing on $H^2(B,{\mathbb Q})$ of \S\ref{s:VHS}, given (up to a constant) by
$$
q(\alpha,\beta)=\,-\int_B\alpha\wedge\beta.
$$
It extends by bi-linearity to $V=H^2(B,{\overline{\mathbb Q}})$. The following hermitian pairing on $H^2(B,{\mathbb C})$ restricts to a well-defined hermitian pairing on $H^2(B,{\overline{\mathbb Q}})$,
$$
h(\alpha,\beta)=\,q(\alpha,\overline{\beta}).
$$
The Hodge decomposition of $H^2(B,{\mathbb C})$ is orthogonal (with respect to $h(\cdot,\cdot)$). Therefore, the subspace $H^{1,1}(B)$ of $V_{\mathbb C}$ is the orthogonal complement in $V_{\mathbb C}$ of $H^{2,0}(B)\oplus H^{0,2}(B)$. Let $U$ be the ${\overline{\mathbb Q}}$-vector subspace of $V$ spanned by $e$ and ${\overline{e}}$, and $U^\perp$ the ${\overline{\mathbb Q}}$-vector subspace of $V$ given by the orthogonal complement in $V$ of $U$. Then $V=U\oplus U^\perp$, where $\dim_{\overline{\mathbb Q}}(U)=2$ and $\dim_{\overline{\mathbb Q}}(U^\perp)=h$, where $h=\dim_{\mathbb C}H^{1,1}(B)$. Therefore $V_{\mathbb C}=U_{\mathbb C}\oplus (U^\perp)_{\mathbb C}$. Clearly $(U^\perp)_{\mathbb C}$ is a ${\mathbb C}$-vector subspace of $(U_{\mathbb C})^\perp=H^{1,1}(B)$. But $\dim_{\mathbb C}((U^\perp)_{\mathbb C})=h$, so $(U^\perp)_{\mathbb C}$ must be all of $H^{1,1}(B)$. Therefore, the complex vector space $H^{1,1}(B)$ is the complexification of a ${\overline{\mathbb Q}}$-vector subspace of $H^2(B,{\overline{\mathbb Q}})$. Combining these facts, it follows that the Hodge filtration of $H^2(B,{\mathbb Q}_B)_{\mathbb C}$ given by
$$
H^{2,0}(B)\subset H^{2,0}(B)\oplus H^{1,1}(B)\subset H^2(B,{\mathbb C})
$$
is defined over ${\overline{\mathbb Q}}$ (see also Lemma 8 and Lemma 9 of \S\ref{s:lemmas}). As $B={\mathcal B}_s$, with $s\in S({\overline{\mathbb Q}})$, our assumptions on the family ${\mathcal B}
\rightarrow S$ imply that $B$ has CM (and hence automatically strong CM). For alternate arguments, see \S\ref{s:lemmas}.

\noindent Therefore, we now know that, if $s\in S({\overline{\mathbb Q}})$ and $\lambda_4\in{\overline{\mathbb Q}}$, $\lambda_4\not=0,1$, then $B={\mathcal B}_s$ and $A_4=A_4(\lambda_4)$ have (strong) CM. The (polarized) rational Hodge structure $H^2(B,{\mathbb Q}_B)$ is the direct sum of the rational Hodge structures $H^2(B,{\mathbb Q}_B)^+$ and $H^2(B,{\mathbb Q}_B)^-$, so that $H^2(B,{\mathbb Q}_B)^\pm$ has CM and the Hodge filtration of $H^2(B,{\mathbb Q}_B)^\pm_{\mathbb C}$ is defined over ${\overline{\mathbb Q}}$. Notice that $H^2(B,{\mathbb C})^+=H^{1,1}(B)^+$, so that $H^2(B,{\mathbb Q}_B)^+$ is concentrated in bi-degree $(1,1)$, which again shows it has CM. Recall that $H^1(B,{\mathbb Q})=\{0\}$, and $H^1(A_4,{\mathbb Q}_{A_4})=H^1(A_4,{\mathbb Q}_{A_4})^-$ has CM and is defined over ${\overline{\mathbb Q}}$. We have
$$
H^2(B,{\mathbb C})^-=H^{2,0}(B)\oplus H^{1,1}(B)^-\oplus H^{0,2}(B),
$$
and
$$
H^3(B\times A_4,{\mathbb Q}_{B\times A_4})^{I_{3,4}}=\bigoplus_{a+b=3}(H^a(B,{\mathbb Q}_B)^+\otimes H^b(A_4,{\mathbb Q}_{A_4})^+)
$$
$$
\oplus \bigoplus_{a+b=3} (H^a(B,{\mathbb Q}_B)^-\otimes H^b(A_4,{\mathbb Q}_{A_4})^-),
$$
so that by Lemma 1, \S\ref{s:lemmas}, the rational Hodge structure $H^3(B\times A_4,{\mathbb Q}_{B\times A_4})^{I_{3,4}}$ has CM and Hodge filtration defined over ${\overline{\mathbb Q}}$.

\noindent It remains to look at the ramification locus $Z$ of $I_{3,4}$. If $Z$ is empty, or consists only of rational curves, then $Z$ clearly has CM. We recall again that by Lemma 3, \S\ref{s:lemmas}, we have the relation between rational Hodge structures
 $$
 H^3(Y,{\mathbb Q}_Y)=H^3(B\times A_4,{\mathbb Q}_{B\times A_4})^{I_{3,4}}\oplus H^1(Z,{\mathbb Q}_Z)(-1).
 $$
From our initial assumptions, we have deduced that $H^3(B\times A_4,{\mathbb Q}_{B\times A_4})^{I_{3,4}}$ has CM. Therefore, since the Mumford-Tate group does not change by the twist $(-1)$, if the rational Hodge structure $H^3(Y,{\mathbb Q}_Y)$ fails to have CM, it will be because there is a component of $Z$ which is a curve of positive genus without CM. By Lemma 4, \S\ref{s:lemmas}, if the Hodge filtration of $H^3(Y,{\mathbb Q}_Y)_{\mathbb C}$ is defined over ${\overline{\mathbb Q}}$, then so is that of $H^1(Z,{\mathbb Q}_Z)_{\mathbb C}$. Therefore, if in addition the curve $Z$ is defined over ${\overline{\mathbb Q}}$, then $H^1(Z,{\mathbb Q}_Z)$ has CM by \cite{Co1}, \cite{SW}. Therefore, if the Hodge filtration of $H^3(Y,{\mathbb Q}_Y)_{\mathbb C}$ is defined over ${\overline{\mathbb Q}}$, $B={\mathcal B}_s$, $s\in{\overline{\mathbb Q}}$, $A_4=A_4(\lambda_4)$, $\lambda_4\in{\overline{\mathbb Q}}$, $\lambda_4\not=0,1$, and $Z$ is defined over ${\overline{\mathbb Q}}$, then $H^3(Y,{\mathbb Q}_Y)$ has CM.

\end{document}